THEORY BUILDING IN GEOMETRY EDUCATION:
DESIGNING AND IMPLEMENTING A COURSE FOR 12-15 YEAR OLD STUDENTS

by

Madhav Kaushish



A Dissertation Submitted to the Faculty of the

DEPARTMENT OF MATHEMATICS

In Partial Fulfillment of the Requirements

For the Degree of

DOCTOR OF PHILOSOPHY

In the Graduate College

THE UNIVERSITY OF ARIZONA

2021



THE UNIVERSITY OF ARIZONA
GRADUATE COLLEGE

As members of the Dissertation Committee, we certify that we have read the dissertation prepared by: Madhav Kaushish

titled: Theory Building in Geometry Education: Designing and implementing a course for 12-15 year old students

and recommend that it be accepted as fulfilling the dissertation requirement for the Degree of Doctor of Philosophy.

| | |
|---|---|
| *Rebecca H McGraw* | Date: Jul 19, 2021 |
| Rebecca H McGraw | |
| *Deborah Hughes Hallett* | Date: Jul 19, 2021 |
| Deborah Hughes Hallett | |
| *Bill McCallum* | Date: Jul 19, 2021 |
| Bill McCallum | |
| *Massimo Piattelli-Palmarini* | Date: Aug 7, 2021 |
| Massimo Piattelli-Palmarini | |

Final approval and acceptance of this dissertation is contingent upon the candidate's submission of the final copies of the dissertation to the Graduate College.

I hereby certify that I have read this dissertation prepared under my direction and recommend that it be accepted as fulfilling the dissertation requirement.

| | |
|---|---|
| *Rebecca H McGraw* | Date: Jul 19, 2021 |
| Rebecca H McGraw | |
| Mathematics | |



## Acknowledgments

This dissertation would not have been possible without the involvement of many other people. Hence, the blame for problems with the dissertation should also be shared!

I would first like to thank all the students who agreed to be part of the course implementation, along with Malati Kalmadi, Pallavi Naik, Lakshmi Gandhi, Anjali Gulanikar and Shital Sonawane, who gave me access to their schools. They also made arrangements to reschedule parts of the course when I fell sick for a couple of weeks. Teaching this course was one of my most enjoyable teaching experiences and I hope to work with these students again, sometime soon.

My advisor, Rebecca McGraw, has been extremely helpful during the entire process. Her comments and feedback have been invaluable in getting the dissertation to where it is. I have also learnt a lot from the rest of my committee, Deborah Hughes-Hallett, Bill McCallum and Massimo Piattelli-Palmarini. Their inputs, especially at the outset, have shaped the direction I took.

I would also like to thank the Mathematics Department at The University of Arizona, especially Tina Deemer and Dave Glickenstein, who made arrangements so that I could do my fieldwork in India in early 2020.

KP Mohanan and Tara Mohanan, along with other members of the ThinQ team, were always available to discuss various ideas. Many of the ideas in this dissertation are the result of conversations with them.

Finally, and most importantly, I would like to thank all my friends and family.



Table of Contents

































**Abstract**

If we want mathematics education to be valuable to everybody, including those who do not pursue mathematics-related careers, we need to use mathematics as a training ground for certain ways of thinking. This dissertation focuses on one aspect of mathematical thinking, namely theory building. It has three parts – unpacking what constitutes theory building, the development of a course aimed at theory building, and a qualitative description of an implementation of that course with 12 to 15-year-old students in two schools in Pune, India.

The description aims to understand some of how students engaged with the course and the role of the instructor in the implementation of the course. The data used in this description includes video recordings of the sessions, teacher reflections at the end of each session, and student reflections and feedback. The ultimate goal of the description is to suggest considerations we should keep in mind when designing and implementing theory-building courses.

The focus of the dissertation is within mathematics, and more specifically, within geometry. However, this is hopefully just a first step towards a transdisciplinary framework aimed at developing thinking abilities useful to everybody in their personal, public, and professional lives.



## Chapter 1 – Introduction

In *How We Think*, Dewey says:

> *While it is not the business of education to prove every statement made, any more than to teach every possible item of information, it is its business to cultivate deep-seated and effective habits of discriminating tested beliefs from mere assertions, guesses, and opinions; to develop a lively, sincere, and open-minded preference for conclusions that are properly grounded, and to ingrain into the individual's working habits methods of inquiry and reasoning appropriate to the various problems that present themselves. (Dewey, 1997)*

If we accept Dewey's vision of education, our role as educators is to develop these habits and the abilities which reside beneath them. Since most of this thesis focuses on mathematics education, we need to ask what role learning mathematics can play.

One of the possible goals of mathematics education is to help students think like mathematicians. However, to be valuable to those who do not end up in mathematics or a mathematics-related discipline, there must be transfer of learning certain ways of thinking in mathematics to other disciplines and one's public, professional and personal lives. There are other important goals of mathematics education, such as learning mathematical content to become a mathematician or to apply mathematics to other disciplines and the real world. Those are not the focus of this thesis.

Building on my Master's Thesis (Kaushish, 2019), this dissertation focuses on one aspect of mathematical thinking: theory building. As mentioned in that thesis, W.T. Gowers distinguishes theory-building and problem solving (Gowers, 2000). Gowers points to areas like algebraic number theory and geometry as the abodes of theory builders, where progress "is often



a result of clever combinations of a wide range of existing results" (pp. 3). On the other hand, the study of graph theory involves many interesting problems, most of whose solutions are not all that dependent on previous results.

Hyman Bass defines theory-building as "creative acts of recognizing, articulating, and naming a mathematical concept or construct that is demonstrably common to a variety of apparently different mathematical situations, a concept or construct that, at least for those engaged in the work, might have had no prior conceptual existence" (Bass, 2017). Naming mathematical objects and constructs is an essential part of theory building. However, I use theory-building more broadly to include activities like Hilbert's axiomatization of Euclidean Geometry (Hilbert, 2014). Hilbert did not generate any new results or new objects but just added more rigor to the theory. I also include activities like the development of alternative geometries, which mathematicians came up with while attempting to prove the parallel postulate.

In this thesis, I will be using theory building to refer to a concern for the relationships between claims and the basic objects of a theory. By basic objects of a theory, I mean axioms and undefined entities, definitions, conjectures/theorems, classificatory systems, and proofs. I call these the 'concepts of theory building' and will discuss them in the literature review in relation to education and clarify how I am using these terms in Chapter 3. One of the main differences between theory building and regular K-12 and undergraduate mathematics is that theory building is concerned about things for which there may not be deductive proofs. For instance, why should mathematicians study parallelograms and not quadrilaterals with two equal sides? Why should squares be types of rectangles rather than separating them? Which definition of prime numbers should be chosen?



Related to each of the objects of a theory are questions of this type. Here is a listing of some questions related to each of the concepts:

- Definitions

    o  Should an object X be included in the theory?

    o  Why should we choose definition A over B?

    o  What is the best way to extend a definition outside its initial context?

- Classificatory systems:

    o  Should A be a type of B?

- Conjectures and theorems:

    o  Is this conjecture plausible?

    o  Is this theorem valuable?

- Proofs

    o  Is this proof rigorous enough?

    o  Is this proof explanatory?

    o  What is the role of diagrams and other visual elements in a proof?

- Axioms and undefined entities

    o  What should our axioms be, and what should we leave undefined?

    o  Can we get more foundational axioms?

This list is not comprehensive. I call these questions, and others like these, 'questions of theory building.' I will unpack some of these questions in Chapter 3. Returning to Dewey's stated goal of education, many of these questions are applicable outside of mathematics, to other academic domains, and other aspects of our lives. I will discuss this further when I discuss transdisciplinarity in this chapter and the discussion chapter.



I have developed a course on theory building and implemented it in two schools in Pune, India. This thesis involves both the laying out of the course and a qualitative description of the implementation. The questions of theory building I will unpack in Chapter 3 are those which are salient to that course.

**Theory Building and Mathematical Practice**

Many of the questions of theory building I described above are pursued as questions about mathematical practice. Investigation into the nature of mathematical practice is a relatively new area of inquiry in the philosophy of mathematics. Traditionally, the discipline has focused on topics like the ontology of mathematical objects and how we have access to abstract mathematical objects (Mancosu, 2008). Inspired by Polya's investigation into problem-solving heuristics, Lakatos appears to have been one of the first philosophers to investigate mathematical practice. However, as Mancosu says, Lakatos's work did not bring this into the mainstream amongst philosophers of mathematics. That has happened much more recently.

Some aspects of mathematical practice have been studied more than others. For instance, as Tappenden (2008a) says, the concept of a 'natural proof' has not been studied in detail while we do have some understanding of 'mathematical explanation.' One way to approach this in education is to wait for future insights of philosophers and then apply them to education. However, that is not the approach I will be taking in this thesis, and that is not something we can afford to do if we wish to bring these learning outcomes into education in a reasonable time frame. So instead, I take insights from this research but build on them.

**Transdisciplinarity**

As I mentioned above, during general education, learning in mathematics ought to be valuable even to those who do not pursue mathematics. Some aspects of mathematics are useful



for application purposes. For instance, arithmetic is helpful for basic money management. An understanding of percentages and probability is helpful for many practical reasons and inquiry in other disciplines.

However, the type of disciplinary transfer I am interested in here is thinking tools or aspects of mathematical practice that can be used at some level of abstraction in other disciplines. For instance, classification in mathematics shares some similarity with classification in biology, reasoning about scientific theories is related to reasoning in mathematics, and deducing consequences of ethical principles is related to deducing consequences of axioms. I am not claiming that these are the same. Rather, there are significant enough similarities that there can hopefully be some valuable transfer of learning. This learning transfer is not just useful for inquiry in other disciplines. For instance, many 'disagreements' we have boil down to different uses of the same word. As an example, there is an often-heard debate amongst children regarding how many fingers they have. The 'disagreement' doesn't have any substance – it depends on the meaning of the term 'finger.' If you consider the thumb to be a finger, most humans have five fingers. Otherwise, they have four. This type of definition-based 'disagreement' comes up again and again in conversation and political debates. Theory building, amongst other outcomes, allows us to see the relationship between definitions and their consequences.

This transfer of thinking tools across disciplines, by abstracting away from particular fields, is what I am terming transdisciplinarity. The course I have designed for this dissertation is not transdisciplinary in nature. It is situated within mathematics and within geometry. Hence, I am not expecting transfer of learning to occur. However, the ultimate goal is to develop a framework and courses which aim at transdisciplinary thinking. This work within mathematics is



intended to be the first step towards that. However, that would require the same unpacking of thinking tools used in other disciplines.

**Geometry Education and Theory Building**

As Gowers (2000) said, geometry is an ideal space for theory builders. The existence of a variety of geometries makes it even more appealing for education in theory building. Additionally, Geometry is one of the few parts of the K-12 curriculum in which students are expected to write proofs. However, many existing geometry courses involve quasi-algorithmic proofs such as 2-column proofs, which require minimal creativity and focus on form rather than substance.

All the modules in the course I have developed involve geometry. One of the modules, called Triangle Theory Building, engages with students' existing understanding of Euclidean geometry as they construct a rigorous theory of Euclidean triangles. Another module, called Discrete Geometry, involves students creating many geometries involving a finite number of points.

**Productive Struggle and Theory Building**

A positive notion of intellectual struggle is something that has existed throughout written history. This is visible in various traditions, including ancient India and Greece. Socrates is the canonical example of a person who saw value in intellectual struggle. The Buddha also discusses the value of self-inquiry and internal intellectual struggle in the Kalama Sutta.

The idea first appears in the education literature, seemingly with Dewey and later with Polya. As Dewey says, the process begins with "some perplexity, confusion or doubt" and continues as students try to fit things together to make sense of them, to work out methods for



resolving the dilemma (Heibert & Grouws, 2007). However, the crucial thing about struggle is that it must be productive.

**Structure and Goals of the Thesis**

There are broadly three aspects to this thesis. The first is the unpacking of the various aspects of theory building, including addressing the questions of theory building raised above, to list some of the learning outcomes of a theory-building curriculum. This will be a part of Chapter 3. Concepts involved in theory building such as proofs, definitions, theories, and conjectures will be the focus of my literature review in Chapter 2, along with a discussion of the literature on productive struggle.

The second aspect of the thesis involved the creation of an approximately 25-hour course on theory building based on the framework. The course structure and the outline of the teaching-learning materials will also be in Chapter 3, with the more detailed materials in the Appendices. I will end Chapter 3 with a list of the research questions for the empirical part of the thesis.

The final aspect is the implementation of the course and the analysis of the results of that implementation. The implementation was done in two schools in Pune, India, which I have worked with before. More details about the schools can be found in Chapter 4, dealing with methods. The two primary goals of the analysis of the implementation are to understand how students engage with such a course and what the role of the instructor is in that process. The third goal is to suggest revisions to the course materials considering how the actual sessions went. Methods of Data Collection and Analysis are in Chapter 4. Chapter 5-10 will contain the qualitative study results, and the final chapter will put things together and discuss the future direction of the broader research program. The data analysis is meant to be illustrative and not



grounds for arriving at conclusions about the goals or the materials and methods used to achieve those goals.



## Chapter 2 – Literature review

In this literature review, I will focus on the education literature relating to Conjecturing & Proving, Defining, and Axiomatic Systems, followed by a summary of the literature pertaining to Productive Struggle. While I will be interested in these concepts generally, I will focus on the literature on these in Geometry Education wherever possible.

**Conjecturing and Proving**

In his book *How to Solve It*, Pólya (2004), referring to a proof that the sum of the angles of a triangle is two right angles, writes:

> *If a student has gone through his mathematics classes without having really understood a few proofs like the foregoing one, he is entitled to address a scorching reproach to his school and to his teachers. In fact, we should distinguish between things of more and less importance. If the student failed to get acquainted with this or that particular geometric fact, he did not miss as much; he may have little use for such facts in later life. But if he failed to get acquainted with geometric proofs, he missed the best and simplest examples of evidence and he missed the best opportunity to acquire the idea of strict reasoning. (pp. 216-7)*

High School Euclidean Geometry focusing on proofs began in the United States in the late 19$^{\text{th}}$ Century (Herbst, 2002). The idea behind it was similar to that in Polya's quote above – to develop certain abilities in students rather than just focusing on knowledge transfer. Herbst shows how this ideal, due to practical considerations, eventually resulted in a disassociation between proving and knowledge construction with the advent of 2 column proofs.



As Weiss & Herbst (2015) point out, the current High School geometry course is often students' first introduction to proof, and sometimes their only example of proof in High School. Proof is rarely touched before High School. As many have pointed out, at least in the United States, the Geometry course is a caricature of actual mathematics, where form triumphs over substance, there are too many postulates, and there is a lack of clarity in the meanings of words (Christofferson, 1930; Weiss et al., 2009; Weiss & Herbst, 2015).

There are different ways in which proof in education has been thought about. Harel & Sowder (1998) introduce the notion of a proof scheme. A person's proof scheme in a particular context is the type of argument they find convincing. They suggest that there are broadly three proof schemes: external-conviction, empirical and deductive. An example of conviction through an external proof scheme is accepting a claim because it is mentioned in a textbook. Accepting a claim because it works on a few examples would constitute being convinced by an empirical proof scheme. Both of these types of proof schemes have been sub-categorized.

Harel & Sowder see two major sub-categories within deductive proof schemes: transformational proof schemes and axiomatic proof schemes. Transformational proof schemes consist of three shared characteristics: generality, operational thought, and logical inference (Harel & Sowder, 2007). Generality is to do with students understanding the class of objects for which they need to prove the claim and not proving it for examples or sub-classes. Operational thought is to do with individuals forming goals and sub-goals along the path to proving the claim. Finally, the proof must follow the rules of logical inference. Axiomatic proof schemes are transformational proof schemes with the understanding that proofs must ultimately start from axioms. While the idea of proof schemes does not contribute to the design of the course, they do play a role in evaluating how students engaged with the course.



It is also important to point out that mathematicians themselves do not necessarily gain conviction about results through deductive proofs (Weber et al., 2014). Many mathematicians will often use authority and/or empirical means to get convinced about results before using them since it would be impossible to work through detailed proofs for each result. Hence, the goal of proof education need not necessarily be that students prove every result but that they see the value of deductive proof.

Weber & Alcock (2004) differentiate between semantic and syntactic proof production. Syntactic proofs are created by correctly stated definitions and related facts in a logically permissible manner. On the other hand, semantic proofs refer to situations in which the prover uses instantiations of a mathematical object the statement is about. Instantiations here refer to systematically repeatable ways in which provers can think about mathematical objects. For instance, an instantiation of a general convergent sequence may be a graph of a 'prototypical' convergent sequence.

### Student Understanding of Proof.

Schoenfeld (2013) shows an interesting result regarding students' understanding of proof after their geometry course – students saw empirical methods as determining truth, while deductive arguments were just exercises teachers gave them. This has been shown in other studies such as Martinez & Pedemonte (2014).

In another study by Martin & Harel (1989), fifty-two percent of student teachers accepted an incorrect deductive argument as proof for an unfamiliar statement. Even after accepting a deductive argument, high school students saw room for potential counter-examples (Fischbein & Kedem, 1982; Harel & Sowder, 2007), while Galbraith (1981) found that over a third of the students they studied did not understand the concept of counter-examples and 18% of them



thought a single counter-example was insufficient to disprove a claim (Battista & Clements, 1995).

Harel & Sowder (2007) found something similar when working with university students. They did not have an axiomatic proof scheme and relied on empirical and authoritative proofs. Even when reasoning deductively, students make inferences that do not follow from their premises. For instance, they conclude $Q \Rightarrow P$ from $P \Rightarrow Q$ . Students also struggle to read proofs and judge their suitability (Hoyles & Healy, 2007; Lin & Yang, 2007).

As mentioned above, students tend to believe that the reason you prove is to complete a task assigned by a teacher or verify something you already know to be true. There is no discovery associated with proof. This is an unfortunate feature of many types of existing mathematics education. One of the goals of the course I have designed for this dissertation is to get students to see the inherent value in proving a claim to gain an understanding of the landscape you are working in.

**Development of activities and models for proof education.**

Given that students conceive of proof as verification or as doing an assigned task, it is clear that proof education is severely lacking. However, there have been efforts to improve this over the years. Polya and Fawcett are two early examples. Polya conceived of proof as problem-solving (Pólya, 2004). He focused on techniques and heuristics to conjecture and prove, such as using empirical means to come up with conjectures, generalizing conjectures to prove them more easily, and so on.

Fawcett (1938) provides a fascinating example of mathematics education research. The research touches on aspects of theory building beyond proof, and I will also return to this work in other sections. Fawcett created a course that was aimed at students constructing Euclidean



Geometry. Students engaged, not just in coming up with conjectures and proofs but also in laying out the axiomatic system on which their proofs were based. The results of the experiment seem to have been universally positive. Not just did students learn how to prove, but they also learned the same amount of geometry that other students knew, as was demonstrated through tests after the experiment. Also, the students who went through the course judged it to have had a significant impact on their lives many years after the fact (Flener, 2009).

While both Polya and Fawcett achieved seeming success, the impact of their work on the actual geometry course has been minimal (Herbst, 2002). More recent research into proof has seen a reconceptualization of the notion of proof, as mentioned above. Proof as a convincing argument and the introduction of levels of proof and proof schemes give new tools to proof educators. Rather than thinking of students as misconceived, these concepts allow us to locate students' thoughts, reasoning, and motivation.

Moving on to some specific areas of interest related to proof production in Geometry, generating examples of concepts has been hypothesized to help in proof production. However, there has not been any significant effect found to this practice – indeed, it seems to be as effective as students studying given examples (Iannone et al., 2011). However, in the same study, they point out that the methodology of example generation might affect their result.

There has also been some research on empirical verification after proof, especially in the case of 'proof problems with diagrams' (Komatsu, 2017), where claims are made in reference to a diagram. In such cases, a claim may be false or require additional specifications since the diagram may contain certain hidden assumptions not made in the original statement. It also may be the case that examining examples after a proof may allow one to prove something more general. The paper also suggests the following roles of the teacher:



1. Prompting students to draw diagrams different from the given diagram or presenting such diagrams

2. Posing questions that move students to either revise the statement or proof

3. Selecting students with ideas worth examining and getting them to present to the class

Concerning another class of proofs, namely existence proofs, Samper et al. (2016) suggest that such proofs are not intuitive. Usually, a student's move is to impose the conditions on a randomly chosen object. The paper highlights the need for the teacher to play a role in mediation and suggests that we cannot expect students to work through such proofs autonomously without the help of a 'more competent doer,' a teacher who has worked through the details themselves.

The educator needs to play a significant role in proof education, especially early on. As Komatsu (2017) suggests, they need to pose the types of questions which will get students to see potential flaws in their proofs or places where they need to be more rigorous. They also need to be able to spot interesting arguments students have made, whether flawed or not, which can deepen the understanding of proof.

**Defining**

Closely related to proof is definition. As DeVilliers (1998) puts it, there are two ways of dealing with definitions: to teach definitions or to teach students to define. The latter does not imply that students ought to come up with every definition. As Freudenthal (1973) says, instructors deny a learning opportunity for students by giving them definitions. He says that most definitions are not preconceived. They are the finishing touches of the organizing activity, and hence students should have experience working with them.



**A Good Definition.**

What makes a definition 'natural' or 'good' is an ongoing area of research. Edwards & Ward (2008), drawing upon van Dormolen & Zaslavsky (2003), give two sets of criteria for good mathematical definitions: necessary and preferred criteria. The preferred criteria include minimality, elegance, and exclusion of degenerate cases. The necessary criteria come primarily from Aristotle:

1. Hierarchy: Objects must be special cases of other objects

2. Existence: The object must be instantiated at least once

3. Equivalence: Multiple definitions must be shown to be equivalent

4. Acclimatization: The definition must fit into a deductive system

Another idea discussed by Tappenden (2008b) is fruitfulness. The consequences of a definition need to be 'significant.' Significant is not a synonym for 'large in number' since if P is true, then (P & P) is also true. Rather, evaluating significance involves considering the explanatory value of the conclusions and the role the definition plays in arguing for those conclusions. Tappenden gives the example of quadratic reciprocity and the Legendre symbol. Before the introduction of the Legendre symbol, the proofs for quadratic reciprocity were tedious and involved many cases. Dirichlet reduced the number of cases in Gauss' proof from eight to two using the Legendre symbol.

In the same article, Tappenden talks about another consideration for a 'natural' definition – that it shows up in different places, especially when not expected. It results in deep connectivity between different parts of mathematics. The Legendre symbol once again is an example of this – it shows up in various number fields and not just where it was defined.



**Definitions and Concepts.**

As Vergnaud (1991) points out, a definition on its own will not enable a learner to apprehend and comprehend a concept (as cited in Ouvrier-Buffet, 2006). Rather, situations and problem solving give a concept meaning (Ouvrier-Buffet, 2006). Hence, the construction of definitions requires students to play with a concept in various situations, extract the important aspects of that concept, and work with multiple representations of that concept before defining.

Defining in geometry involves many different aspects. Unlike in other areas of mathematics, students have access to mental pictures and diagrams of the objects they are defining. However, this can also make things more difficult. Mariotti and Fishbein (1997) talk about harmonizing the conceptual and figural aspects of geometric objects. A paradigmatic example of this is that of a square and rectangle. Students, especially in elementary school (Bartolini Bussi & Baccaglini-Frank, 2015; Kaur, 2015; Tsamir et al., 2015), tend to see squares and rectangles as different entities rather than squares being a subset of the class of rectangles. There is a conflict between their perceptual experiences, the figural aspects, and the need to unify and generalize the conceptual aspects (Bartolini Bussi & Baccaglini-Frank, 2015; Mariotti & Fischbein, 1997).

Tall & Vinner's distinction between Concept Image and the Concept Definition (Tall & Vinner, 1981) has had a significant impact on the field of Mathematics Education Research. The Concept Image is the set of pictures, representations, properties, and statements associated with the concept in the learner's mind. The Concept Definition is the actual definition of the object. The most significant impact is that it has shown that learning the definition is not enough – students need to understand the object the definition is referring to.



Zandieh & Rassmussen (2010) use the concept Image-Definition framework along with Gravemeijer's (1999) RME activity framework to create a framework for definition as a Mathematical Activity (DMA). They use this framework to construct a series of activities that transition students from triangles on the plane to creating and inquiring into the concept of spherical triangles. The activity they use is an example of extending a definition outside of its initial scope. This will form an important part of the course I have designed, and I will discuss this idea further in Chapter 3.

Mariotti & Fischbein (1997) discuss two types of definitions: the basic objects of the theory and new elements within the theory defined in terms of the basic objects. These basic entities have a close relationship with the axioms of the theory. Using this and Fischbein's conceptual-figural distinction discussed above, the paper proposes a pedagogy for coming up with definitions via a classroom discussion. This involves:

1. Observing

2. Identifying the main characteristics

3. Stating properties based on them

4. Returning to observations to check

In a teaching experiment on defining and classifying quadrilaterals, Fujita et al. (2019) found that through a semiotic/dialogic process, students were able to transform their intuitions of what parallelograms were to a collective notion. Even though this notion did not necessarily agree with the conventional definition, students were able to use their definitions to solve other problems. Coming up with definitions based on ideas and examples is the basis of one of the modules of the course I have designed. Also, one of the tasks in another module involves classifying triangles via a similar process to the one in Fujita et al. (2019).



**Axiomatic Systems**

Fawcett (1938) is one of the first examples of students constructing an axiomatic system in a high school course. Students listed words they had heard of in relation to geometry and named properties related to them. They then named some statements as axioms, some words as undefined, defined others, and wrote proofs to statements deduced from this system. Students then used the same thinking tools they had developed to evaluate things outside of mathematics, including legislation.

Healy's Build-a-book Geometry (1993) is another example of working with axiomatic systems. Unlike Fawcett, where students started with ideas and observations, Healy began by giving students three axioms: parallel lines never meet, triangles contain 180 degrees, and circles contain 360 degrees. From this, students worked cooperatively in groups to solve problems they posed on their own. They would then present their findings to the class. The students performed about as well on standardized exams as more traditional classes. However, they also learned other skills like self-confidence, teamwork, creativity, and responsibility (Ridgway & Healy, 1997).

Apart from Fawcett and Healy, there has been some other work on different axiomatic systems, mainly on Spherical Geometry and Discrete Geometry (Lenart, 1996; Junius, 2008; Zandieh & Rassmussen, 2010; Ada & Kurtulus, 2012).

**Spherical and discrete geometry.**

Spherical Geometry refers to a 2D geometry where the surface we are concerned with is a sphere. We are concerned with similar things as we would be in Euclidean Geometry, but we are constrained to the surface of a sphere. The advent of Spherical and other geometries resulted in a significant change in the nature of mathematics from viewing Euclidean geometry as the 'only true geometry' to a large number of geometries coming into existence.



Lenart (1996), making a case for teaching different types of geometry alongside Euclidean Geometry, suggests that Spherical Geometry is especially useful since the sphere shape is prevalent naturally, most obviously in the shape of the Earth (Sinclair et al., 2016). Lenart suggests introducing what he calls 'comparative geometry' where students work on spherical and plane geometry together. He says that if one of the goals of mathematics education is that students learn axiomatization, a single context is not sufficient.

While exploring different geometries, there are different views we can take – we can take an extrinsic view by thinking of a sphere being embedded in Euclidean 3-space, or we could think of the sphere are the space itself. Junius (2008) worked with students to move them from an extrinsic view of straightness of lines to an intrinsic view. It requires us to take the intrinsic view to see great circles as equivalent to straight lines. Exploring straight lines in various worlds is a part of the discrete geometry module in the course I have designed. Seeing discrete worlds as the space we are working in instead of points in Euclidean space is also an important aspect of the same module.

Junius takes the concept of straight line from Euclidean geometry and extends it to the sphere. Zandieh & Rassmussen (2010), as discussed in a previous section, worked on the notion of spherical triangles with students by extending the notion of planar triangles, which students were familiar with. Extending definitions outside their initial domain of definition is something I will be discussing in the next chapter.

Unlike Spherical Geometry, Discrete Geometry doesn't refer to a single axiomatic system. Rather, it refers to a collection of geometries that are non-gradient. Some of these geometries have points with internal structure, such as geometry related to simplices or some pixel geometries. Others, such as geometric graph theory, do not. Taxi-cab geometry is one such



geometry that has been used at the K-12 level (Ada & Kurtulus, 2012). Taxi-cab geometry is

modeled on a grid city. You can only travel along the streets and not through the blocks. So,

there can be multiple straight lines between two given points. Discrete geometry will be an

important part of the course I have designed. However, taxi-cab geometry is not part of that

course.

**Productive Struggle**

Heibert and Grouws (2007) define struggle as students expending effort to make sense of

mathematics, to figure out something which is not immediately apparent. As they write, struggle

does not mean needless frustration or extreme levels of challenge created by nonsensical or

overly difficult problems.

This positive notion of struggle has been a part of education literature since at least

Dewey and Polya. As Dewey (1926, from Heibert and Grouws, 2007) claims, deep knowledge

"is the fruit of the undertakings that transform a problematic situation into a resolved one."

For the instructor/course designer, this implies a tension between making a task

challenging enough so that it results in struggle but not so challenging that it makes students

disengage. This is related to Vygotsky's zone of proximal development (ZPD). If the goal is to

induce struggle in students, tasks need to be such that they are on the boundary of a student's

ZPD. Polya (2004) mentions something similar. He wrote, "Helping the student. One of the most

important tasks of the teacher is to help his students. This task is not quite easy; it demands time,

practice, devotion, and sound principles. The student should acquire as much experience of

independent work as possible. But if he is left alone with his problem without any help or with

insufficient help, he may make no progress at all. If the teacher helps too much, nothing is left to



the student. The teacher should help, but not too much and not too little, so that the student shall have a reasonable share of work."

Hence, ensuring productive struggle is a matter of both appropriate and challenging task design and the art of deciding when the right time is for the instructor to intervene. The role of the instructor in ensuring productive struggle is similar to the role of the instructor in a problem-based learning (PBL) course. Problem-based learning courses usually consist of open problems (Barrows, 2000). Hence, the role of the instructor in a PBL course is not to lead students to the 'correct answer' but to model good thinking practices (Hmelo-Silver & Barrows, 2006) and engage with responses. This is how I see the role of the instructor in the course I have designed.

**Summary**

While theory building has not been emphasized in the education literature, different aspects of it have been studied and have played a role in unpacking theory building, in the development of the course, and in describing the implementation of the course in this dissertation. The idea of proof schemes comes up when describing the implementation of the course. Fawcett's (1938) course on creating Euclidean geometry has helped me shape the triangle theory-building module. The literature on productive struggle and problem-based learning gives a useful pedagogical approach for the course. The literature on mathematical practice in philosophy has also played an important role in various aspects of this thesis. Tappenden's (2008a; 2008b) work on definitions has been especially useful in thinking about choosing between definitions and extending definitions outside their original context. I will discuss these further in the next chapter, in which I will unpack the learning outcomes associated with theory building and lay out the course I have designed based on those learning outcomes.



## Chapter 3 – A course on theory building

As mentioned in the introduction, this thesis involves developing and implementing a course on Theory Building. This chapter consists of three parts. The first part unpacks learning outcomes associated with Theory Building. The second part lays out the modules of the course which aim at these learning outcomes. The last part involves listing the research questions I will address in the thesis's empirical part.

The goal of the course is to equip students with the ability to construct mathematical knowledge with a focus on theory building. This involves understanding various concepts of mathematical theories like definitions, theorems, proofs, axioms, and undefined entities. It also involves an understanding of heuristics and considerations mathematicians use when constructing knowledge in mathematics – the answers to the questions of theory building I mentioned in the introduction.

### Learning outcomes associated with theory building

To unpack the learning outcomes associated with Theory Building, it would be useful to first answer the question – what is the difference between Theory Building and 'Doing Mathematics'? Bass' notion of theory building includes activities that involve a concern for the structure of mathematics. However, Bass doesn't clarify the notion of structure. As mentioned in the introduction, I will be using structure to mean the relationships between statements of a theory and the basic objects of the body of knowledge. This focus is different from Polya's work on problem-solving, where the goals are conjecturing, solving problems, and proving. Of course, any theory-building exercise will involve conjecturing and proving, and Polya's work does involve some concern for the theory. However, the distinction here is a matter of focus.



This section will have two parts. In the first part, I will discuss the concepts of theory building, which are touched upon in the literature review, including axioms, definitions, conjectures, proofs, etc. The second part is to do with the questions of theory building, which include thinking about what we should consider when choosing one definition over another, deciding upon which classificatory system to use, what counts as proof, what we should take as axioms, and so on.

Before jumping into that, it is useful to present a high-level course structure so that I can refer to the activities in the course. I will discuss these activities in more detail later in this chapter.

The course is designed for a minimum of five sessions of ninety minutes each and consists of three modules. The first session is dedicated to the first module. At a minimum, two sessions need to be spent on each of the other two modules, but each of these modules could go on for weeks. Each of the three modules involves different activity types. I am calling these activity types:

1. Definition guessing

2. Assumption digging and integration

3. Creating worlds by changing assumptions

The particular modules related to these activity types are:

a. Podgons (Definition guessing)

b. Triangle Theory Building (Assumption digging and integration)

c. Discrete Geometry (Creating worlds by changing assumptions)

I will be laying out and unpacking these activities in detail later in this chapter.



**Concepts of theory building.**

The concepts involved in mathematical theory building include theorems, conjectures, proofs, definitions, classificatory systems, undefined entities, and axioms. Many of these have already been discussed in the context of the literature in mathematics education. So, this section is just focused on clarifying how I am using these terms and the role they will play in the course.

*Conjectures and theorems.*

A conjecture is a claim we have not proved or shown to be false. If we find a proof to the conjecture, it becomes a theorem. If we show it is false, its negation is a theorem, at least when working with two-valued logic. Conjecturing, or coming up with conjectures, will play a role in the course modules, mainly in Discrete Geometry and Podgons. I will discuss the art of conjecturing and plausible conjectures while discussing the questions of theory building.

*Proofs.*

A proof is an argument that shows that a particular claim is true. In other words, it converts a conjecture into a theorem. A proof can rely on statements previously proved or on axioms of the theory. Proving will form a significant part of the Triangle Theory Building module. Finding counter-examples, which prove the negation of a claim, forms an important part of the Podgon module. Further along in this chapter, I will discuss what constitutes a proof in terms of rigor and what makes a proof a good proof.

*Definitions.*

Definitions in mathematics take the form of if and only if statements. They will play an important role in all of the three modules. I will discuss various aspects of definitions when discussing the questions of theory building.



### *Undefined entities and axioms.*

Every theory consists of undefined entities which are constrained by axioms. For instance, in the usual telling of Euclidean Geometry, points and lines are undefined but are constrained by axioms like 'given two points, there exists exactly one straight line which passes through it.' These come up in Triangle Theory Building and, to some extent, in Discrete Geometry. I will discuss what entities ought to count as undefined entities and what ought to count as axioms further along in this chapter.

### *Classificatory systems.*

Classificatory systems are closely related to definitions. Whether a square is a type of rectangle is determined by the definition of square and rectangle we choose. Classificatory systems give us logical inheritance, which means that statements which are justified for the super-category are automatically true about the sub-category. Classificatory systems play a role in Triangle Theory Building.

### **Questions of theory building.**

In the previous section, I gave short explanations for the different concepts involved in theory building. However, I have not touched upon how they are used and what role they play in mathematics. This section will ask various questions related to these concepts and discuss the considerations that need to be kept in mind while addressing them.

There is an infinite range of mathematical theories and objects we could be interested in. Within Euclidean Geometry, we study objects like triangles but not closed curves whose boundary consists of exactly two straight lines but could have non-straight line parts. Also, we study Euclidean Geometry in great detail, but many other types of geometry are not given as much importance.



We also see that some definitions of the same object are more revealing and natural than others. For instance, we could define even numbers as those that end in 0, 2, 4, 6, or 8 in base 10, or we could define them as numbers divisible by 2. The latter definition feels more natural and seemingly gives more insight. Similarly, we use a classificatory system with squares as types of rectangles but not one that sees them as distinct categories.

We accept certain proofs as being sufficiently rigorous but not others. We accept certain systems of axioms but not others. For instance, we could just state every result in Euclidean Geometry as an axiom. Since they are not inconsistent, that doesn't result in a contradiction. Yet, we prefer a small set of axioms.

Some of these decisions may have had a sociological, applicational, or contingent historical basis. Maybe we give importance to Euclidean Geometry over other geometries because it happened to be the first type of geometry that was developed. In this section, I will largely ignore the sociological, applicational, and historical bases for these decisions. Rather, my focus will be on unpacking the heuristics and considerations we ought to be keeping in mind while making such decisions.

Polya's discussions on heuristics for doing mathematics would fall into this category. Lakatos' Proofs and Refutations is another example of work that concentrates on these 'softer' aspects of mathematical thinking. Apart from Lakatos, people like Kitcher (1984) and more recently Mancosu (2008) and others have made contributions to this area in the Philosophy of Mathematics literature. Some of what follows has a basis in this work.

As I mentioned in the introduction, I will only be unpacking a few of the questions of theory building that you could ask. I will be concentrating on:

- Should X be in our inventory?



- How do we choose a definition?

- Which definition should we choose to extend?

- Should A be a type of B?

- Is this conjecture plausible?

- What should our axioms be, and what should we leave undefined?

- When is a proof?

***Should X be in our inventory?.***

As mentioned above, there are certain objects we decide to study and certain objects we do not. We tend to give names to the objects we do study, like triangles, circles, squares, and so on. These objects can be thought of as entries in an inventory. In what follows, I will attempt to discuss some of the considerations to keep in mind when deciding to include a particular object in the inventory.

Tappenden (2008a) discusses the idea of addition, comparing the 'plus' function to a 'quus' function. The 'plus' function is regular addition, while 'quus' represents an infinite class of functions that agree with the 'plus' function for any conceivable computation humans would want to do. However, it can differ on other pairs of numbers, for instance, numbers that are sufficiently large. It does seem that the 'plus' function is more natural than any possible 'quus' function. One argument given for that is 'plus' is more intuitively natural and simple. This may work for something like 'plus,' which is a concept available to most people. However, when dealing with more complex areas of mathematics that require training, it isn't clear whether our intuitions are a good guide.



It is useful here to compare two objects: parallelograms and quadrilaterals with two sides equal. Parallelograms have a lot of theorems about them. For instance, the following are equivalent definitions of parallelograms:

1. Quadrilaterals with opposite sides equal

2. Quadrilaterals with opposite sides parallel

3. Quadrilaterals whose diagonals bisect each other

On the other hand, while there may be something to say about quadrilaterals with any two sides equal, there almost certainly are not a lot of things you can say about them which are not true about quadrilaterals generally. So, a possible consideration when deciding whether to introduce a concept to our inventory is whether there are interesting things to say about that concept. Behind this is the idea of 'fruitfulness' (Tappenden, 2008b). We include those objects into our inventory which result in significant consequences. Significant here is not synonymous with 'large in number' since if P is true, then (P & P) is also true. It is to do with the value of the consequences in 'explaining' the mathematics. For the purposes of this thesis, I do not intend to go deeper here and will assume an intuitive understanding of 'explanation' in mathematics.

To take another example, right-angled triangles are objects we study in textbooks. However, triangles with one 47.3 degree angle is not a class of objects widely studied. On the face of it, there doesn't appear to be a difference – in both cases, we are specifying one angle of a triangle. There are a few possible reasons for this. One reason, contingent on the reality we live in, is that we use right angles in areas like engineering. Another reason is that in Euclidean Geometry, there is less arbitrariness when dealing with right angles than with other angles. For instance, given a point and a line, there is exactly one perpendicular that can be drawn passing through that point to the line. For any other angle, that is not the case. A third reason is that with



right triangles, equations like the one given by Pythagoras' theorem are cleaner. You could use the law of cosines to write out the theorem for other angles, but there would be extra coefficients in those equations.

So, when adding entries to our inventory, some considerations to keep in mind are the significance and range of consequences of the objects, the lack of arbitrariness in those properties, and the simplicity when writing out theorems about those objects. One thing to note here is that definitions often do not come first as they do in textbooks. It is often conclusions that arise first, and definitions are ways in which to systematize them. An example I will be discussing in the next subsection in a slightly different context is that of the Legendre symbol. Various proofs of quadratic reciprocity existed well before the Legendre symbol was first written down. Amongst other things, what the symbol did was to make those proofs much shorter.

### How do we choose a definition?.

Given an object, we may not yet have a well-articulated definition of that object. Rather, we may have a general sense of what the object should be along with some conclusions we wish to be true about that object. That results in a choice of possible definitions. Some of these definitions may be equivalent in the theory we are working in – for instance, the parallelogram definitions mentioned above. Others may have different consequences.

Consider the example of a triangle. How do we define a triangle in Euclidean Geometry such that we get the consequences that we want? Clearly, three straight line segments is not sufficient since that may not even give us the concept of angles. What about a closed shape made of three straight line segments? The consequences of this depend on what we mean by 'made of.' By one notion of this term, it is impossible to have a shape that consists of exactly three line segments since every line segment can be broken into smaller line segments, and hence every



object 'made up' of three line segments is made up of an infinite number of line segments. The other notion would be along the lines of, 'there are three line segments which make up the entire shape such that every point on the shape is on at least one of the line segments and only the end points can be shared by any two of the line segments.' One potential issue here is that if we take a straight line segment and two sub-line segments of that line segment which together make up the line segment, that would be a triangle. It may be possible to reject this by saying that the concept of a closed shape requires a non-empty inside, which would then require us to define inside.

Even if we do this, the problem remains when we generalize. Let us use the same definition we used for triangle but for quadrilateral – we just replace three with four. Now, take something which looks like a triangle. That is also a quadrilateral since we can break it into four line segments that share only end points. It also doesn't have a non-empty inside. So, now we have a choice: we can define triangles in such a way that what looks like a triangle is also a quadrilateral, a pentagon, and so on, and that what looks like a straight line could be a triangle, or we could exclude such possibilities. We can do the latter by stating something like, 'three consecutive vertices should not be collinear.'

What are the considerations we should keep in mind while choosing between these two? One is the idea of maximal generality, while the other is the properties we want the object to have. Clearly, the first way of defining triangles and other polygons is more general – every example of the second type is an example of the first type. However, by using the first definition, we lose out on certain properties like the strong triangle inequality. Depending on the use cases, it might be better to use one definition over the other.



Moving on to a case where definitions are equivalent. Even in this case, there may be reason to prefer one definition over another as our main definition. Take, for instance, the concept of even numbers. We could define them to be natural numbers that end in 0, 2, 4, 6, or 8 when written in base 10. Alternatively, we could define them as natural numbers which are divisible by two. The second definition seems to give us a better understanding of what even numbers are. Also, it isn't dependent on something as arbitrary as the base we are working in. The concept of 'two' is more general than how it is represented in a particular base. So, while we may want to use the other definition for particular purposes, such as checking whether a particular number in base 10 is even, the definition to do with divisibility is more central to the actual concept of evenness.

Taking the case of three equivalent definitions of a parallelogram:

1.  Quadrilateral with opposite sides parallel

2.  Quadrilateral with opposite sides equal

3.  Quadrilateral whose diagonals bisect each other

All three of these entail the other two. Hence, these are equivalent definitions. There is no obvious reason that one of these is a better definition than the others. However, while I have not been able to find any literature on this, it appears hardest to mentally visualize how the third definition entails the other two. If that is the case for most people, that is also a useful consideration when deciding upon which definition to pick.

To summarize, in this subsection, the focus has been on what considerations to keep in mind when deciding upon a definition. Some of those highlighted here are:

1.  Matching the definition to the mental picture we have as in the case of triangles



2. Defining something in order to get consequences we desire, as was the case with the polyhedra and the triangles

3. The level of generality we define objects at as in the case of triangles

4. The centrality of the definition to the concept as in the case of even numbers

5. The ease of using, understanding, and manipulating the definition as in the case of parallelograms or even numbers

Some of these considerations may come into conflict at times, as shown in some of the examples. At that point, the mathematician would have to weigh the different considerations and make a judgment. One other consideration worth thinking about is to do with extending a definition outside its initial domain of applicability. You may find that some definitions do better than others in such cases. That is what I will discuss in the next subsection.

***What definition should we choose to extend?.***

Given a definition in a particular theory, we may want to use the concept underlying the definition in another theory. This new theory may be a generalization of the original theory, as in the case of Rings being a generalization of the Integers, or it may be a different theory with contradictory axioms, such as spherical vs. Euclidean geometry.

When we decide to extend a concept outside the theory it exists in, we often have various choices we can make since it may be the case that certain concepts do not exist in the new theory.

Tappenden (2008a) gives the example of the definition of prime numbers. The elementary school definition of prime numbers is those numbers that are only divisible by one and themselves. Another definition of prime numbers is that a number is prime if whenever it divides a product of two numbers, it at least divides one of them. In the case of natural numbers,



these are equivalent definitions. However, when we generalize, this is no longer the case. In certain integer rings, numbers are prime in the first sense but not in the second. Hence, we have a choice here as to which of these we want to take as the definition of prime numbers to extend. Mathematicians choose the second definition and not the one we learn in elementary school because many of the things we consider significant regarding prime numbers rely on that.

As another example, take the case of Taxi-cab geometry. For intuition about what Taxi-cab geometry is, think about a grid city where you are driving a taxi along the streets. Suppose we want to borrow the concept of straight-line segment from Euclidean Geometry. Straight-line segments in Euclidean geometry have many properties. For instance:

1. given the concept of angle, we can say straight line segments are paths in the same direction

2. straight line segments are the shortest paths between two points

3. given two points, there exists a straight-line segment between them

4. given two points, the straight-line segments between them is unique

In Taxi-cab geometry, it isn't totally clear what 'in the same direction' would mean. Even if we do have a notion of direction, that would be to do with paths along a particular street. Hence, there would not be 'straight line segments' between any two given points and only between those on the same street if we wanted to use direction to define them. However, if there existed a straight-line segment, that would be unique. If we use shortest path, on the other hand, there are always shortest paths given any two points. However, paths are almost never unique. So, even though the two conceptions of straight-line segments give us the same consequences in Euclidean Geometry, they do not in Taxi-cab geometry. So, if we want to extend the concept



from Euclidean Geometry to Taxi-cab geometry, we need to make a choice as to which of the properties we want to preserve.

### Should A be a type of B?.

Closely related to definitions is the notion of classification. Should a square be a rectangle or not? For descriptive purposes, when talking about the real world, it is reasonable to argue that squares should not be rectangles. However, for academic purposes, how do we make this choice? If squares are not rectangles, then the definition of rectangles has to contain something like: 'adjacent sides cannot be equal.'

What are we giving up if we do include that in our definition – if we treat squares as non-rectangles? We lose at least two things. Firstly, we lose logical inheritance of properties. When treating squares as rectangles, anything we prove about rectangles is automatically proved about squares. This is no longer the case if squares are not rectangles, and we will have to prove a lot more statements in this case.

Secondly, even stating definitions is harder. To illustrate this in a convincing manner, it would be useful to try and compare two different classificatory systems for quadrilaterals. The first is where we treat squares, rectangles, parallelograms, trapeziums, and 'other' quadrilaterals as separate categories with no intersections. The second is where trapeziums are types of quadrilaterals, parallelograms are types of trapeziums, rectangles are types of parallelograms, and squares are types of rectangles.



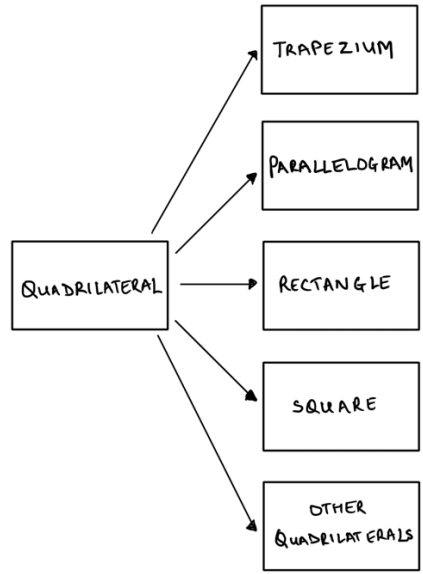

*Figure 1: Classificatory System A*

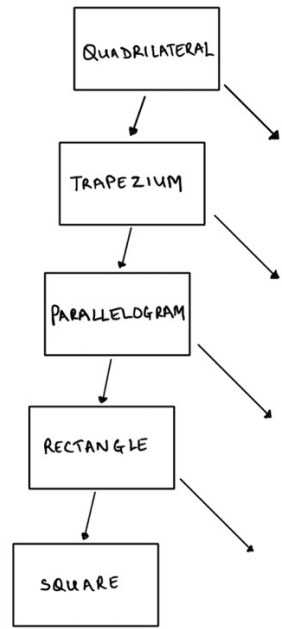

*Figure 2 Classificatory System B*

Classificatory System B gives us the following definitions:

*Trapezium*: Quadrilateral with one pair of parallel sides

*Parallelogram*: Trapezium with the parallel sides being equal

*Rectangle*: Parallelogram with one right angle



*Square*: Rectangle with adjacent sides equal

Classificatory System A gives the following definitions:

*Square*: Quadrilateral with all sides equal and one right angle

*Rectangle*: Quadrilateral with opposite sides equal, one right angle, and adjacent

sides not equal

*Parallelogram*: Quadrilateral with opposite sides parallel and no right angle

*Trapezium*: Quadrilateral with one pair of opposite sides parallel and the other

pair not parallel

*Other Quadrilaterals*: Quadrilateral with no pair of opposite sides parallel

Notice that if we want to include rhombuses and other shapes in our system as non-intersecting categories, the definitions get even harder to state.

You see something similar in the case of biological taxonomy. We could have humans and animals as distinct. In other words, humans not being animals. As an alternative, we could say humans are primates which are mammals, which are vertebrates, which are animals. In the latter case, things we say about the general categories are inherited by the more specified categories. Also, to specify what a human is, all you need to do is to state some set of properties that distinguish humans from other primates.

So, some of the important considerations while classifying are logical inheritance and the simplicity of definitions.

### *What makes a conjecture plausible?.*

Plausible conjecturing is an art that Polya discusses at length in Volume 1 of Mathematics and Plausible Reasoning. In that book, Polya makes a case for induction and analogy as being useful tools in the process of discovery in mathematics. Induction here refers to



empirical induction as opposed to mathematical induction. There are a few ways he discusses

using induction. One is trying a sample of examples of what the claim is saying. Another is

looking at the consequences of a claim. If those consequences happen to be true, the claim is

more likely to be true. Similarly, if something which implies another claim turns out to be false,

our degree of belief in the plausibility of that other claim should decrease.

Analogy is another way to establish plausibility. Polya discusses an example of the

infinite series made of the sum of reciprocals of squares. Euler saw the relation between this and

the roots of a polynomial. He used this to develop a technique to compute the sum of the series.

This was before he rigorously showed that this relationship existed.

Given that there are infinite possible conjectures we could be interested in, establishing

the plausibility of a conjecture is valuable before spending time attempting to justify it.

### *What should our axioms be, and what should we leave undefined?.*

At the base of any mathematical theory lie undefined entities and axioms which constrain

them. For example, in Euclidean Geometry by Hilbert's axiomatization, points and lines are

undefined, but there are various axioms, such as there being exactly one line segment between

two points, which constrain how these entities can behave. In set theory, the concept 'in' is

undefined. In Peano Arithmetic, the successor function is undefined. Defined entities are defined

from these undefined entities and other defined entities.

In principle, we could treat every entity in Euclidean Geometry as undefined and treat

every claim about those entities as axioms. Clearly, that goes against one of the main purposes of

knowledge construction – to unify disparate conclusions so that we don't have to memorize

individual statements and that we get a more abstract understanding. Also, given a large number

of axioms, it would be very hard to tell whether they contradict each other.



For a particular group of practitioners, especially those learning how to do mathematics, this requirement of minimizing the axioms needs to be balanced by the understandability of the system and ease of computation. Hence the answer to 'when is an axiom?' in a classroom is something that will have to be negotiated between the students and instructors.

The idea of undefined entities is present in academia outside of mathematics. For instance, Newton's notion of mass is undefined. The three laws are the equivalent of axioms which constrain what sorts of concepts can apply to the word mass. Similarly, the concept of syllable in theoretical linguistics doesn't have an if and only if definition. Rather, the concept is constrained by a large number of statements, the equivalent of axioms, about the concept.

### *When is a proof?.*

This is a similar question to 'when is an axiom?'. However, the concern, in this case, is with the level of rigor. This is a question Keith Devlin asked in a blog post (Devlin, 2003). He says there are two ways of looking at what a proof is. The first is to think of a proof as a logically correct argument that establishes the truth of a proposition, while the other is an argument that convinces a typical mathematician of the truth of a given statement. He shows that the prior descends into the latter if we look carefully at it since it is human mathematicians who have to evaluate the correctness of the proof.

Hence, a better question than 'what is a proof?' is 'when is a proof?'. Devlin illustrates this with examples from research mathematics. For the purpose I am asking the question, it might be more useful to illustrate this with proofs you might encounter in a classroom.

Take the claim: equilateral triangles in Euclidean geometry are equiangular. At one level, that could be an 'obviously true' assertion. However, if that is not acceptable, you may attempt to justify the claim by appealing to some notion of congruence – you can say that if you drop a



median from one vertex, the two triangles will be congruent by SSS, and hence the angles will be equal. You can do this for two vertices giving you that all the angles are equal. While the congruence assumption is explicitly stated, there are many other assumptions here that are not explicit. For instance, the proof does not explicitly state that it is always possible to draw a median. It also does not explicitly state that equality of angles has transitivity. A proof with those assumptions explicitly stated would be more rigorous than one with those not stated. You could go even deeper and appeal to logical assumptions such as modus ponens or the law of the excluded middle. If you explicitly state that you are accepting modus ponens, that makes your proof even more rigorous.

Of course, explicitly stating these many assumptions in every proof would be an extremely cumbersome task. Hence, the level of rigor expected needs to be negotiated in the classroom and will have to vary from task to task depending on the particular goals for that task.

**Summary.**

These learning outcomes form the heart of the course I have designed. Words like 'good,' 'natural,' 'rigorous,' 'explanatory,' 'structure,' come up again and again in the language mathematicians use. However, students, especially in K-12, almost never see this side of mathematics. They see mathematics as solving problems or at best as conjecturing and proving within a theory. While developing an interest in mathematics is not one of the main goals of this course, these 'softer' aspects of mathematics could potentially attract those students who were turned off by the lack of creativity in what they saw as the discipline.

In the next section, I will lay out the modules in the course designed for this thesis. The course consists of three modules, one on definitions, the second focused on rigorous reasoning, and the third on an exploration of various worlds.



## Course Activities and Structure

The course is designed for a minimum of five sessions of ninety minutes each and consists of three modules. In the actual implementation for this thesis, the course was run for nine sessions. The first session is dedicated to the first module. At a minimum, two sessions need to be spent on each of the other two modules, but each of these modules could go on for weeks. In the implementation, three sessions were spent on each of the two modules, starting from the third session. The second session was dedicated to an example of assumption digging, and the final session was spent on a short definition guessing example run by the students, while the rest of the time was spent on the post-course survey and the post-test.

In this section, for each of the three modules, I will first describe the general activity and then discuss the specific module. More details about the specific modules are in the appendices.

### Podgons - Definition guessing.

The goals of this activity type involve understanding definitions and looking for counter-examples, which is a type of proof. In the course, it is meant as an introductory ice-breaker activity that is not intended to be very rigorous. However, it can be done as a rigorous activity.

Definition guessing involves figuring out a definition when given objects which are and are not members of the class being defined. This is done in the form of a game. Two or more people are given definitions of the same word. I will refer to them as 'definers.' These definitions have at least some different consequences. Each of the definers are asked to give an example of an object which is a member of the class by their definition and an example that is not. The students then attempt to guess the definitions. There will of course be many possible definitions people come up with. They can check whether a particular definition works for a particular



definer by asking whether the examples they come up with belong to that class or not. The role

of the definers can be played by the instructor(s) or by multiple students.

As an example, let us take the activity which is used in the course designed for this thesis.

The instructor takes the role of two definers. Let's call them A and B. The word they are defining

is podgon. They both give the following examples:

| Podgon | Non-Podgon |
|---|---|
| 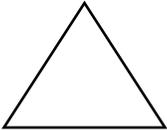 | 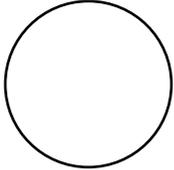 |

*Figure 3 Podgon or not?*

The students come up with the following definitions:

Def 1: An equilateral triangle

Def 2: Any shape consisting of straight lines

In order to test these definitions, students have to come up with examples that test the

definitions. For example, they may posit an example of a triangle whose sides are different

lengths. If that turned out to be a podgon for A, Def 1 would not be a possible definition for A.

So, students have to search for potential counter-examples to a given definition – examples that

fit the definition but not the judgment of a definer or examples which fit the judgment of a

definer but not the definition.

**Triangle Theory Building - Assumption digging and integration.**

The goals of this module involve an understanding of all of the concepts of theory

building mentioned in the previous section, with the focus being on proofs. In terms of the

practice-related learning outcomes, this module aims at all of those mentioned above apart from

extending definitions and judging the plausibility of conjectures.



To explain what this module is about, I will first unpack the notion of assumption digging and then discuss integration. For a more detailed unpacking of assumption digging, see my Master's Thesis (Kaushish, 2019), which was focused entirely on it.

### Assumption digging.

As I said in my Master's thesis, the metaphor of assumption digging comes from archaeology. Imagine being an archaeologist who has discovered what looks like the top of an ancient ruin. You would like to uncover that ruin. In order to do that, you will have to dig. As you dig, you realize that parts of the building are not in great shape. You have to repair them before you dig further in order to preserve the structure. In this analogy, the top of the building represents conclusions, intermediate floors of the building represent intermediate conclusions, the walls and pillars represent justification, and the foundation of the building represents the axioms, undefined entities, and definitions the theory is built from. Parts of the building in bad shape represent hidden assumptions or bad arguments.

An example of this type of activity is what Hilbert engaged in with his Grundlagen – finding flaws in Euclid's geometry and fixing them by adding assumptions. The motivation here is to place a theory on firmer footing. So, given a conclusion and an argument purporting to prove that conclusion, Assumption Digging is the process of clarifying terms and statements, extracting hidden premises, and clearly spelling out the steps of reasoning to a certain degree of rigor. In a classroom, the decision on what constitutes a rigorous enough justification will have to be decided upon collectively.

To unpack this further, it would be useful to use a diagrammatic representation of Assumption Digging. I am using the same example here as was used in my Master's thesis.



The conclusion (referred to as C) we began with was: Triangles with the same base, and between the same pair of parallel lines, have the same area. The proof for this relies on the following three premises:

Premise 1 (P1): The height of two triangles between the same pair of parallel lines is the same

Premise 2 (P2): The area of a triangle is ½ base x height

Premise 3 (P3): The bases of the two triangles are of the same length

Figure 4 is a diagrammatic representation of this argument.

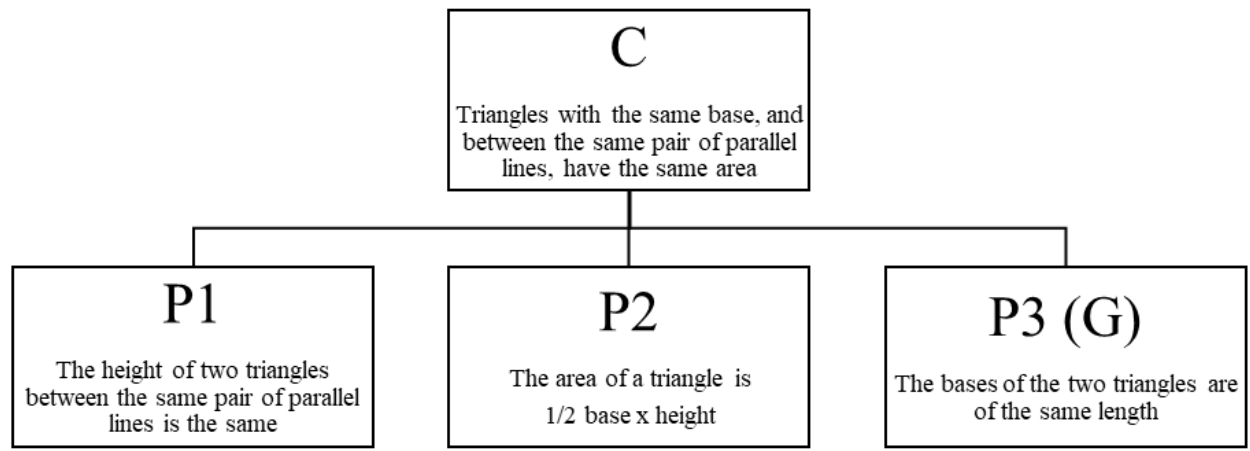

*Figure 4 Tree representation of the argument above.*

At the top of the tree is a conclusion we have begun with. We could begin with more than one conclusion, but the representation would be much messier. Since the particular example used in this thesis starts with just one conclusion, the above representation is sufficient. The proof given for the conclusion (C) is based on the three premises at the second level of the tree. There are three possible types of premises: Definitions, Axioms, Given, and Claims (which may be theorems – those require a proof).

P3 is given in the statement of the claim we are interested in, and that is indicated by a G. Taking Premise 2, we could start with the area of a triangle being ½ base x height as the



definition of area of a triangle. However, we could also define the area of a triangle in terms of the area of a parallelogram. Here are the premises involved in that argument.

Premise 2-1 (P2-1): Area of a parallelogram is base x height

Premise 2-2 (P2-2): Given a triangle, there exists a parallelogram with the same base and height and double the area.

Figure 4 is a tree representation of the entire argument so far.

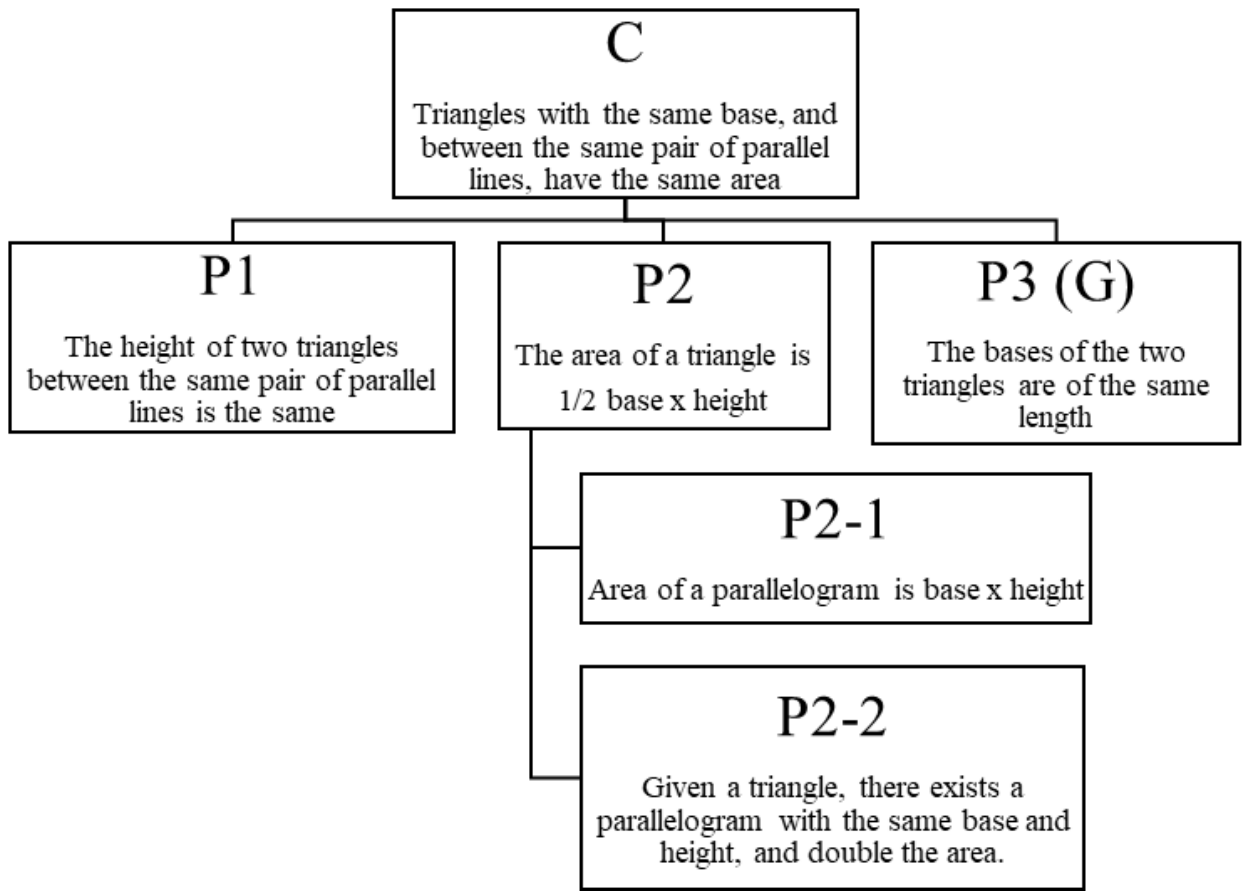

*Figure 5 Tree representation of the argument above.*

P2-2 will take some work to prove rigorously – it requires the notion of congruence. With P2-1, we once again have the choice of using that as a definition, or we could define the area of a parallelogram in terms of the area of a rectangle.



Taking P1, it will eventually require some form of the parallel postulate, even though we

do not get to that in the actual sessions. Figure 6 represents this.

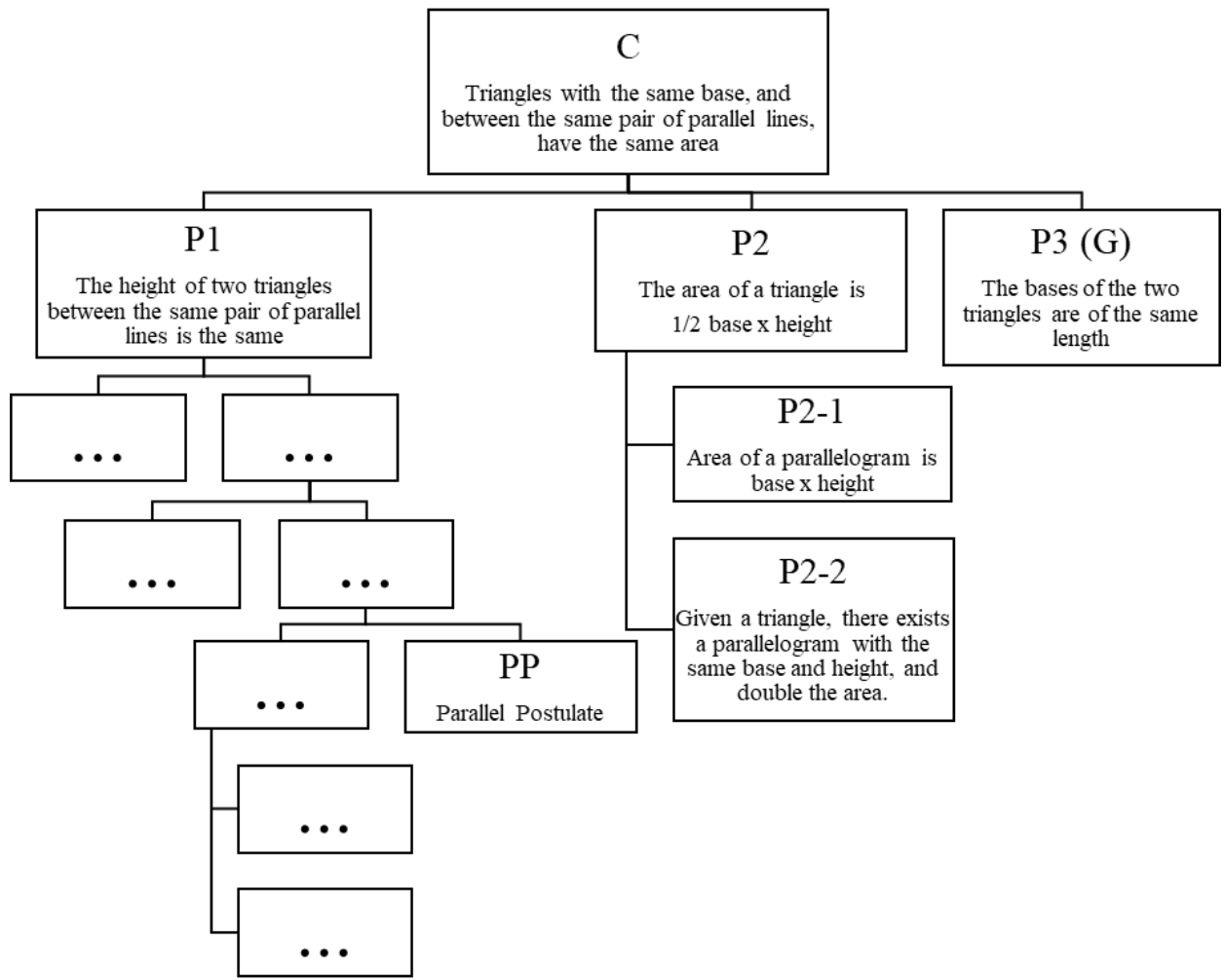

*Figure 6 Tree representation of an argument going all the way down to the Parallel Postulate.*

PP represents the parallel postulate in some form. It may well be the case that the same

axiom/definition is used at various different points in the diagram – to identify these pictorially

would make it extremely messy, and the tree structure would be lost. What is also missing in the

diagram is an indication as to what the undefined entities are. While they will be mentioned in

the axioms, they also need to be explicitly mentioned. An example of an undefined entity in



Euclid's version of his geometry is a straight line. It is undefined but is constrained by certain axioms.

To summarize, the vertices of the tree represent propositions. The root of the tree is a conclusion, while the other vertices are premises on which the conclusion is based. Once the digging is complete, the leaves represent axioms and definitions. The other premises are intermediate conclusions and could be considered theorems in their own right. The edges of the tree represent the steps of reasoning from each layer of premises to the conclusion one step above.

An important part of assumption digging, as mentioned above, is to fix flaws in the structure. These flaws usually come in two flavors, one related to the edges and the other to the vertices:

1. leaves without justification which have not been stated as axioms or definitions

2. flaws in the edges – in the steps of reasoning

The assumption digging task, then, becomes an exercise in fixing these flaws. When we see a leaf without justification from a lower layer, there are three possibilities – the leaf contradicts other statements in the theory and should be dropped, the leaf is justified by a lower layer of premises, or the leaf is stated as an axiom or a definition. In order to show the first case, we look for counter-examples that satisfy other claims we wish to accept. It may well be the case that we find a counter-example while attempting a justification.

***Integration.***

This is the heart of what Theory Building aims at – putting together different things into a coherent whole. An example of integration outside of mathematics is Newtonian gravity. Newton managed to explain a large number of seemingly unrelated observational generalizations using a



small set of statements. In the history of mathematics, it seems like Euclidean Geometry was the first attempt at doing something similar – using a small number of axioms and definitions, it attempts to justify a large number of claims.

Writing proofs and assumption digging result in a certain level of integration. They establish relationships between certain statements. I will refer to that type of integration as vertical integration since it establishes a partial order amongst statements.

The other type of integration involves noticing similar objects and statements which occur throughout the theory and abstracting them away from the specific context. I will refer to this as horizontal integration. A couple of examples of horizontal integration are:

1. Noticing that objects with four straight-line sides with opposite sides parallel show up in various theorems, giving that object a label (parallelogram) and studying that object.

2. Noticing that certain statements, for which you do not have justification, show up repeatedly. You can then choose to state those as axioms.

The triangle theory-building module of this course works as follows: Students list out a large number of things they know about triangles. Then, they pick one claim and engage in Assumption Digging. Ideally, their proofs would involve other claims they have listed. If not, they add those claims to their list. Eventually, the goal is to justify the large list of facts they know about triangles in terms of a smaller set of axioms and definitions. Given the time constraints in this course, the amount of horizontal integration students would get the chance to do is minimal. So, the idea will be to just sensitize students to the idea.



**Discrete Geometry - Creating worlds by changing assumptions.**

The goals of this activity type involve an understanding of all the concepts of theory building apart from classification. While proving is involved, that is not the primary focus. In terms of the practice-related learning outcomes, Discrete Geometry has a primary focus on extending definitions, choosing definitions, and conjecturing.

In this activity type, we are looking to explore how changing the axioms and definitions of our theory changes the consequences. The canonical example of this is dropping the parallel postulate from Euclidean Geometry. By doing that, we end up creating multiple new geometries. There are two types of tasks that I class under this type. The first is taking a set of clearly articulated axioms and definitions, changing a subset of them, and following the consequences. Since this is meant to be an introductory course, this type of task is not a part of it. The other type of task is to take a more intuitive understanding of a theory, change certain assumptions within that theory and formalize the consequences. This is the type of task that will form a large part of this course.

The example used in this course involves asking the question, 'In a world with exactly 6 points, can every straight line segment be bisected?'. The words used give the impression that we are dealing with geometry. However, the sort of Geometry students have dealt with is Euclidean Geometry. In that situation, they know that the answer to the question is yes – every Euclidean line segment can be bisected. However, this question is extremely ambiguous since it requires us to understand what a world with 6 points looks like, what the rules of such a world are, what straight lines in that world are, and what bisection could mean.

This requires students to define these concepts and set up rules so that they can answer the question. This is not to be done arbitrarily – for instance, we don't want to define bisection as



a triangle with two equal sides. However, for many concepts, it may not be possible to preserve all the characteristics when we are trying to construct a new theory. So, the goal is to preserve aspects of the object which we consider to be valuable, as mentioned when discussing extending definitions in the previous section. As an example, we may create a world without a useful concept of direction but with a useful concept of distance. In that world, defining straight lines in terms of direction isn't useful, but we can define straight lines as shortest paths.

Once we have set up at least some of the rules of the world, the next step of this activity is to explore how other objects behave in this world. For instance, taking the definition of circle from Euclidean Geometry as the set of points equidistant from a center, we can explore what circles look like in the different worlds that we create. We can then generalize from particular worlds with 6 points to different types of worlds with n points.

**Summary.**

Each of the modules is looking to give students some insight into different aspects of theory building. Podgons aims at giving students an understanding of definitions and of justifying through counter-examples. Triangle Theory Building and Discrete Geometry look to give students an experience of engaging with some of the questions of theory building as well as giving them an understanding of the various concepts of theory building.

The Triangle Theory Building and Discrete Geometry modules are relatively open in terms of the number of sessions a course could dedicate to them. Each of them could range from two sessions to an entire semester, and maybe even a whole year. In this implementation, I chose to spend three sessions on each of them. The Podgon module was used in this implementation as an ice-breaker on the first day. Hence, it was not implemented as rigorously as it could be. Activities of the definition guessing type could be done independently in a very rigorous manner



and could involve students being the definers. In fact, in this implementation, the last session did involve a definition guessing activity in which students took these roles. However, I have not included it in my analysis since the video footage is unusable.

Detailed plans for each of the modules are available in Appendix A. The next chapter will be dedicated to the methodology for the implementation of the course as well as to the methodology for the analysis and reporting of the data gleaned from that. The next section explores the research questions for the empirical part of this dissertation, which the methodology looks to address.

**Empirical Research Questions**

The goal of the empirical part of this thesis is to explore what a course on theory building looks like when implemented. The study is qualitative and will involve describing the classroom interactions between the students and me, the instructor, in detail. The questions I will be guided by while doing this are:

1. How do 12-15-year-old students engage with a course on mathematical theory building?

2. What is the role a facilitator can play in guiding the direction of such a course?

3. What are the various things to keep in mind when designing a course of this type?

While these are the guiding research questions, the empirical part of this thesis is focused on a single course implemented in two schools that are part of the same school group and had a single facilitator. Hence, I make no claim to address these questions in any generality. Rather, the empirical work should be seen as addressing the following set of narrower questions:

1. How did the students in the two schools engage with different aspects of this course?

2. What are the various ways in which the facilitator guided the direction of the course?



Using the results of the first two questions, what should be kept in mind when attempting future iterations of this course?



## Chapter 4 – Methods

**Background and Logistics**

  **Course implementation.**

   I taught the course to two groups of students in two schools in Pune, Maharashtra. I will be referring to these schools as Indus and Ganga. They are relatively low-cost private schools, owned by the same non-profit.

   The course consisted of three modules relating to the three activities discussed in Chapter 4: Podgons, Triangle Theory Building and Discrete Geometry. Apart from that, there was a session specifically on Assumption Digging about circles in the second session of the course which I do not have video recordings of due to issues with my equipment.

   Podgon was a single session long in both the schools while the other two modules spanned three sessions each. A session ranged from 1 hour 30 minutes to 2 hours. The podgon session was on the first day following by the assumption digging session. After that, I alternated between the other two modules after that starting with Triangle Theory Building.

   Each of the schools had 15 students. On the first day in Ganga, half the students were missing due to something they had to attend outside of school. Apart from that day, there were at most two students missing on a given day in Ganga. In Indus, at most one student was missing on any given day. In Ganga, I was joined by a friend, Mohanan, to help out during the first session while another friend, Sabareesh, was present for multiple sessions in both the schools.

   The following is the schedule for the two schools. P refers to Podgons, C to circle assumption digging, TTB to triangle theory building and DG to discrete geometry. The session lengths are rounded to the closest half hour.



Ganga:

|         | Day 1 | Day 2 | Day 3 | Day 4 | Day 5 | Day 6 | Day 7 | Day 8 | Day 9 |
|---------|-------|-------|-------|-------|-------|-------|-------|-------|-------|
| Module  | P     | C     | TTB   | DG    | TTB   | DG    | DG    | TTB   | Wrap  |
| Length (Min) | 120 | 90 | 90 | 90 | 90 | 90 | 120 | 90 | 90 |

Indus:

|         | Day 1 | Day 2 | Day 3 | Day 4 | Day 5 | Day 6 | Day 7 | Day 8 | Day 9 |
|---------|-------|-------|-------|-------|-------|-------|-------|-------|-------|
| Module  | P     | C     | TTB   | DG    | TTB   | DG    | TTB   | DG    | Wrap  |
| Length (Min) | 120 | 90 | 90 | 90 | 90 | 120 | 120 | 90 | 90 |

**My history with the schools.**

Through an organization I co-founded called ThinQ, I have worked with the two schools in the past. ThinQ set up a course for 6[th] graders at the two schools aimed at developing Academic Inquiry abilities across disciplines, including mathematics. I have also done various workshops in the school on Mathematical Thinking, including one which turned into my Masters Thesis. Most of the students in the workshop in Ganga went through the 6[th] grade course, and some of the students in both the schools attended previous workshops of mine. Hence, to give context, it would be useful to elaborate on these.

***Sixth grade course.***

I had some involvement with creating the materials for this course and was involved during the original implementation, but most of the work was done by other members of ThinQ



since I was just getting ready to start graduate school when the course began. The goal of the course was to develop in students what we call Academic Inquiry abilities and mindset by which we mean the abilities and mindsets which go into being a good mathematician, scientist, philosopher, historian, and so on. We aim at developing these abilities through students actually creating new knowledge in these fields. New knowledge here doesn't necessarily refer to something new to everybody, but something new to the students.

An example of an extended scientific inquiry activity we did with students was related to the movement of the Earth and Sun. Given phenomena, students had to construct alternative models which explain those phenomena. We started with the existence of day and night. For the purposes of this exercise, we assumed the Earth and Sun were spherical, but you could imagine a session where you add that as another variable trait. Here are some of the models students came up with:

1. Earth rotating in a 24 hour period with the sun stationery

2. Earth revolving around the sun in 24 hours

3. Sun revolving around the Earth in 24 hours

4. Earth and Sun revolving around each other

We then added in further phenomena like properties of shadows, and eventually moved to seasons along with other phenomena. Students were able to reject some of these models since they didn't explain particular phenomena. Rather than students being told about the relationships between the Earth and the Sun, they had to figure them out for themselves.

Another example is to do with the cross sections of solid shapes including spheres, cylinders and cubes. We asked students to imagine what the cross sections could be and what sort of cut would result in what cross section. We then cut objects in order to check whether their



conjectures worked. The value of this exercise is the development of visual imagination, which would not have happened if they were given the objects to cut initially.

### *Workshops.*

I did workshops in the summer of 2017, 2018 and 2019. The 2017 and 2019 workshops were single day long workshops and were only in Ganga, while the 2018 workshop was four days long in each of the schools. The 2018 workshop formed the data I eventually used for my Masters Thesis. The focus of the 2017 workshop was on Mathematical Thinking more generally, while the 2018 and 2019 workshops were on Theory Building.

*Theory building workshop 2018.*

The following are some of the activities involved in the workshop:

- Straight Lines and Intersections: This is an example of long-form problem solving which can result in Theory Building. We only spent a couple of hours on this and, hence, didn't get very far. This activity works far better in the context of a course than in an intensive workshop since we can come back to it again and again. I had intended on doing something like this in this course, but I decided against it due to lack of time.

- Coloring Problems: Given various map coloring problems students had to guess at the minimum number of colors required to color such maps such that no two adjacent regions had the same color. Some examples were maps where regions were bounded by straight lines, by circles, and so on. They then had to classify maps by the number of colors required to color them.

- Discrete Geometry: I did a very short version of what did in this course. I only did it in Ganga.



- Assumption Digging: The claim I used for assumption digging was – Two triangles with the same base and between two parallel lines have the same area.

**Student Background.**

The specific background of the students related to my interactions with them is given in the previous section. In this section, I will be concerned with the more general context of education in the country and the school.

### *About mathematics in the schools.*

The schools are relatively low-priced private schools affiliated to the CBSE board for the purposes of 10th and 12th grade examinations. Almost half of all students in India go to private schools. Within the schools which are affiliated with CBSE, these schools are somewhere in the middle in terms of their annual fees.

For mathematics, they use textbooks from NCERT, and the medium of instruction is English. Euclidean geometry is introduced formally in the 9th grade in the NCERT textbooks. Prior to that, the geometry chapters in the 7th and 8th grades are concerned with compass and straight edge constructions, and with understanding polygons and their properties.

The Euclidean Geometry course begins with an introduction to Euclid's axioms, postulates, and definitions. At the end of that section, it states that while what comes next is Euclidean Geometry, it isn't based on the same axiomatic system as Euclid's due to it not being satisfactory. In the Euclidean Geometry chapters, students are given definitions, axioms, and proofs to all the important theorems. They have to come up with proofs themselves, but only to less crucial theorems. Conjecturing is not explicitly encouraged. This is even more the case outside of geometry where proving is not the focus.



### *Student demographic information.*

I have used pseudonyms for specific students in this thesis including in this section. This information is from a survey I asked students to fill in after the first session of the course. The survey is in Appendix C.

*Ganga.*

In Ganga, there were seven students who listed their gender as female and eight who listed their gender as male. Three students were in the eighth grade (Vivaan, Anya and Vivek), one in the seventh grade (Imran), and the rest in the ninth grade. Apart from two of the students (Anuj and Jaya), all the rest had been part of the sixth-grade course described above, and three of the students had been a part of my Master's Thesis workshop in 2018 (Vivaan, Gauri and Jaya).

*Indus.*

In Indus, six students identified as male while nine identified as female. Arvind, Anil and Devika were eighth graders while the rest were in the ninth grade. Eight students had been through the sixth grade course and all but five students (Serena, Anil, Arun, Sanya and Meghna) were part of the workshop for my Master's Thesis.

## Data collection

The goal of collecting data was to provide insight in order to address the research questions. The first two research questions involve student engagement with the course and the role of the instructor. The main mode of engagement with the course were class discussions. Hence, a two-camera video of those class discussions, one pointed at the instructor and the other at the students, was the most important part of the data.



Supplementing that was video of group discussions, my reflections at the end of each day, student written reflections, surveys at the beginning and end of the course, and tests at the beginning and end of the course.

The following expands on each of these data sources and their purpose:

**Video recordings (class and group discussions).**

The sessions were recorded with two cameras, one pointed at me and the board, and the other at the students. These cameras were stationary and placed to capture as wide a view as possible. Apart from the mics on the camera, there were also three mics placed amongst the students, one in each group. As mentioned above, the video consisted of class and group discussions. The main object of analysis was the class discussion. The role of the group discussion data was to give context and to provide evidence for or against particular interpretations of what students said.

**Initial survey.**

In the first week of the course, I conducted a survey. It elicited demographic information along with students' prior interactions with mathematics, me or with materials I have been involved with creating. The list of questions in the survey are in Appendix C. The main goal of this survey was to get basic demographic information along with some idea of previous interactions students may have had with me since these interactions could have an impact on how they interact with the course.

**Pre and Post Tests.**

I carried out written assessments before and after the course. The tests aimed to assess some aspects of Theory Building but only a small part of this vast area. Both the tests had the similar questions with specific features changed. The Pretest and Posttest are in Appendix B. At



this point, I am not looking to arrive at any general conclusions from these tests. However, I will

be presenting the data from the tests to give a sense of how a course of this type can be evaluated

and to encourage future research into the instrument.

**Student written reflections.**

Student written reflections potentially help in giving a sense of how they engaged with

the course. Hence, I asked students to write written reflections four times during the course at the

end of a session. I did not ask students to write reflections after each session since, when I have

done so in the past, it tends to become a formulaic exercise rather than something students spend

effort on. The four days I chose were the first and last of the triangle theory building and discrete

geometry modules. This gave students' reactions to the set-up of the module and their reflections

on the module as a whole.

**Post-course survey.**

The post-course survey was intended to give a sense of how students experienced the

course as a whole and how they thought they engaged with the course. The questions asked were

to do with what students enjoyed and did not about the course, what they found valuable and not

valuable during the course, and what changes they would recommend to the course. Each of

these questions asked for an explanation. The goal here was to see how students experienced the

course. This is especially useful as a way to look at those students who did not speak up during

the class discussions.

**Teacher reflections.**

Since one of the research questions is to do with the role of the teacher, my reflections on

the sessions were a source of data to address that. I used these reflections in the analysis in order



to narrow down data which was used and to give context to certain decisions I made during the course. I wrote reflections on the sessions each week. The focus of these reflections was:

1. The differences between what was intended and what occurred

2. The moves I made which seemed to work and those which didn't, giving specific reasons why I thought so

3. What I would have done differently

4. The difference between the sessions in the two schools

The reason for these choices is that they address the second and third research questions from the point of view of the teacher of the course.

**Data analysis and reporting**

Using the data, I intend to address the following research questions related to the course implementation:

1. How do 12-15 year old students engage with a course aimed at developing theory building abilities?

2. What is the role a facilitator can play in guiding the direction of such a course?

3. What are the various things to keep in mind when designing a course of this type?

Since these are very general research questions, I do not have the ability to answer them given the data I have access to. So, as a first step to address the research questions above, I will answer the following related questions:

a. How did students engage with this particular implementation of this course?

b. What role did the facilitator play in guiding the direction of this course?

c. What are some things to keep in mind while implementing and refining this course given questions a and b?



While I will attempt to suggest some more general learnings in relation to these questions, most of the analysis is local to specific parts of the implementation. The local analysis involves particular parts of the implementation of this module in the two schools, addressing students' engagement with particular aspects of the course, to do with the first research question, as well as with how I facilitated the sessions, to do with the second research question. To achieve this, I have attempted to interpret what students have said given the direct context along with their group discussions and written feedback. I have also attempted to suggest alternative ways in which I could have facilitated particular parts of the course.

**Analyzing and reporting the video data.**

The main mode of engagement with the course in terms of time spent was the class discussions. Hence, the focus of my analysis will be the videos, specifically the class discussion parts of the videos.

***Coding the video.***

The first way I divided the class discussion videos was into the set-up and the rest of the module. By set-up I mean the initial question or the first few examples in the module. The reason for concentrating on this is the assumption that the initial framing and understanding of a module or session can have a significant impact on learning. For these parts of the video, I will be reporting the class discussions sequentially from beginning to end.

For the Podgon module, by set-up, I mean the initial question, the first definition and an evaluation of that definition. For the Triangle Theory Building module, the set-up includes the first session and the claims we evaluated fully during that session. For Discrete Geometry, it involves the time students spent addressing the initial question posed.



After that, I went through the class discussion parts of the videos and coded them in three ways. The first is by the mathematical content we are dealing with. In the definition guessing module, this coding is not relevant. However, in the discrete geometry and triangle theory building, it is. In triangle theory building, the mathematical content will consist of various claims students are evaluating. In discrete geometry, the different types of content are the different questions being asked. The reason for this type of coding is that it divides these two modules into natural units which require minimal external context to understand.

The second type of coding is to do with individual student interactions in the class discussions. The goal of this coding was to relate what students said to the learning outcomes of the course. The units of coding are individual student statements since the focus was student engagement. While coding the videos, I was interested in the activities I have mentioned in Chapter 3, i.e., conjecturing, proving, classifying, defining, and assumptions (which could be axioms). The codes I used are below:

| Code | Description | Example |
|---|---|---|
| Asking for Proof | When students ask for the justification of a claim either from other students or from the instructor | Student A says claim P is true since it follows from claim Q. Student B asks why they should believe claim Q |
| Asking for Definition | When students ask for the definition of a word either from other students or from the instructor | Instructor asks whether straight lines can be bisected. Student responds saying, 'that depends on what you mean by straight lines.' |
| Providing Proof | When a student justifies a claim either when asked to do so or spontaneously. While coding, I will include all sorts of justification including | Student A says claim P is true since it follows from claim Q. Student B asks why they should believe claim Q. |



|  | empirical and appeal to authority. This also involves showing that a claim is false. | Student A responds to Student B with an argument for claim Q. |
| --- | --- | --- |
| Evaluating a Proof | When a student points out flaws in a proof given by themselves, another student, or the instructor. | Student A gives a proof for a claim P.<br><br>Student B mentions that the proof only works for specific cases and doesn't work generally |
| Providing Definition | When a student gives a definition for a word. This could be something they remember or something they have come up with. | Instructor says, 'given that the pictures on the board are examples of circumscription or examples of non-circumscription as indicated, what is a possible definition for circumscription?<br><br>Student responds by giving a definition of circumscription. |
| Evaluating a Definition | When a student argues that a definition doesn't capture the intended concept. This includes giving counter-examples: objects which are examples of the concept but do not fit the intended definition or the other way around. | Student A gives a definition of Rectangle which includes Squares.<br><br>Student B suggests that square should not be included in the category of rectangles since in regular language it is not. |
| Defending a Definition | When a potential flaw is pointed out in a definition, a student defends the definition saying that what seems to be a flaw isn't one, and argues for that. | In the above example, Student A retorts that squares should be rectangles since including them doesn't change any theorems, and we should include it since the more general, the better. |
| Asking Questions | This will involve all math related questions students ask which are not asking for justification or definition. | During a session on triangles on a sphere, student asks whether there is a theorem about the sum of their angles like there is on the plane. |



| Coming up with Conjectures | When a student makes a claim which has not been stated before. | Student A defines a triangle to be a polygon with 3 angles. Student B makes the claim that triangles have 3 straight line sides. |
|---|---|---|
| Refining Conjecture | When a student changes a conjecture to add more clarity. | In the previous example, the instructor gives an example of a shape with 3 straight line sides and a few curved sides. The student refines the definition to say exactly 3 straight line sides |
| Stating assumptions | A student makes a claim which they say should be treated like an assumption/axiom. | In order to prove P, student claims Q. When asked to prove Q, student says we should treat it as an axiom. |
| Extracting assumptions | A student identifies an assumption being made in an argument or in order to justify a claim. | Student states that the proof worked only because we assumed the triangle was acute. |

For this type of coding, for the podgon module, I coded each statement made by students if they fit into one of these codes. For the other two modules, I did the same for the set-up of the module. For the rest of the module, I initially coded chunks of mathematical content (the first type of coding) with these codes if there were student statements in those which reflected them. I then narrowed down the data to be analyzed using this initial coding and the mathematical content coding, as I will discuss in the next subsection. For the narrowed down data, I coded them by statement as I had done with the set-up.

The third type of coding was to do with emerging themes. For the definition guessing module, I took individual statements which had been coded by the second type of coding and looked for general themes within those, using my teacher notes as a guide. I chose two themes to



focus on: reasons provided by students for examples they posited and unwarranted conclusions students arrived at. I then picked examples I judged to be representative of those themes.

For the other two modules, I looked for themes across the narrowed down set of data as well as the set-up. These themes formed additional codes to the chunks already coded by mathematical content. These themes are reflected in the general observations and reflections I have made at the end of each chapter of the results.

### *Narrowing the video data to be analyzed and reported.*

For the Podgon module, the data to be analyzed was not narrowed. For the other two modules, due to the large amount of data, I narrowed down data to be analyzed and reported corresponding to the mathematical content codes mentioned above, given the following constraints:

1. Along with the set-up of the module, individual student statements span as diverse a set of the learning-outcome related codes as possible to give a sense of different ways in which students engaged with the course.

2. They include as many of the episodes I have mentioned in my post-session notes. The reason for doing this is based on the assumption that my post-session notes shed insight into the role of the instructor.

3. The episodes are not all in the same session – there should be variance across when the episodes occur. This is so that there is variance in content being covered.

### *Transcribing the video data and images.*

I then transcribed all the classroom discussions of the episodes I had narrowed my search to. The transcription was done to get the essence of what students were saying. While I have tried to be as faithful as possible to their words, I have ignored pauses and 'ums'. In terms of



images of the board-work I have used in the description, the quality of the actual images was not very good. Hence, I have recreated the board work I have used. I have tried to remain as faithful to the original while making the text more legible.

### *Using other data.*

While analyzing the class discussion parts of the video, I also looked at video of the group discussions nearby the class discussion in question and my teacher notes. The goal of using these was to either give more context or to give evidence for or against a particular interpretation of something in the class discussions. As I mentioned above, given that I wanted the results chapter to be self-contained and containing as much context as possible, I have only referred to the group discussion data where it is necessary for choosing between two interpretations.

### Student feedback and tests.

Another type of data I am reporting is students' feedback as well as the results of the pre and post-tests I gave them.

The student feedback involves three aspects: reflections at the end of a session (which students in each of the schools wrote four times), end of course feedback, and the audio interviews which students did at the beginning and end of the course. The student feedback is used in two ways – for general comments students made about the course and for specific comments students made about particular modules.

For the pre and post-tests, I have reported the raw scores of students and the changes to the scores from the pre-test to the post-test. I am not intending to arrive at any conclusions based on that change. I also discuss the ways in which students thought about the questions as they reported to me during the audio interviews.



**Structure of results chapters**

There are six chapters dedicated to results, one the podgon modules, two each for the other two modules, and the last dedicated to the student feedback and tests.

The first results chapter will be on the Podgon module. That is broken into two three parts. The first part focuses on the set-up, the second on student feedback related to this module, and the third on general observations from the implementation of this module.

Triangle Theory Building and Discrete Geometry have two chapters each dedicated to them. They are structured very similarly. The first chapter is on the set-up of the modules, and the second is on selected episodes from the rest of the module, student feedback related to these modules and the last on general observations from the implementation of these modules.

The sixth results chapter is focused on student feedback and the pre and post-tests. The section on the pre and post-tests involves reporting the scores and how they changed. The student feedback part is primarily focused on the feedback which is general to the course and which was not focused on in the earlier chapters.



## Chapter 5 – Podgons: Definition guessing

As I mentioned in Chapter 3, this module in the course was intended to be an ice-breaker activity. However, it resulted in some interesting discussions. The learning goals of this module are relatively narrow. It focuses on defining and on finding counterexamples.

The goal of this chapter is to give a description of the session focusing on how students engaged with the module and the role played by the instructor in steering the module. I will also be reflecting on possible changes to the module for future iterations.

After describing the module, the rest of this chapter is in two parts. The first part focuses on the set-up of the module. In the case of this module, by the set up I mean the first two attempts at defining what podgons are. I will present the class discussion and interpret what students say in light of the goals of the course. I will also use other sources of data such as the group discussions and my teacher notes where they are helpful in interpreting.

The second part focuses on general reflections and observations on the module. Here, I will be using the data I presented in the set-up part along with other instances during the implementation of the module in the two schools.

I did do another definition guessing session on the last day of the course with the students playing the role of the definers. However, the video from that day is not usable.

The following is a table of contents for this chapter:







**Description of the module**

The module involves students figuring out a definition when given objects which are and are not members of the class being defined. This is done in the form of a game. Two or more people are given definitions of the same word. I will refer to them as 'definers.' These definitions have at least some different consequences. Each of the definers are asked to give an example of an object which is a member of the class by their definition and an example that is not. The audience then attempts to guess the definition. There will of course be many possible definitions people come up with. They can check whether a particular definition works for a particular definer by asking whether examples they come up with belong to that class or not. The role of the definers can be played by the instructor(s) or by multiple students.

As an example, the instructor takes the role of two definers. Let's call them A and B. The word they are defining is podgon. They both give the following examples:

| Shape | Podgon A | Podgon B |
|:---:|:---:|:---:|
| 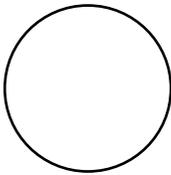 | 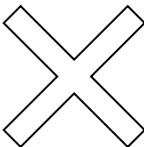 | 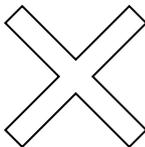 |



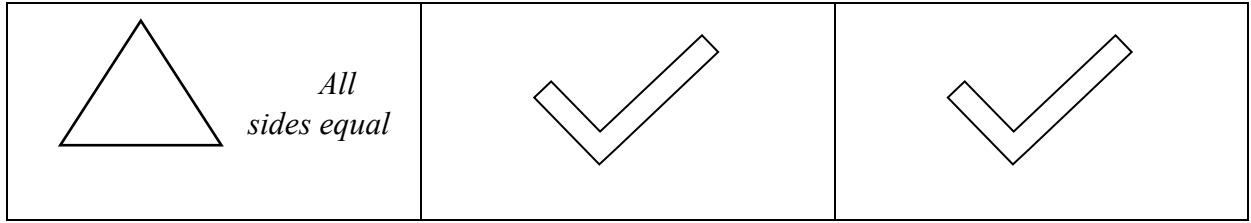

*Figure 7 Podgon or not?*

The students invent the following definitions:

Def 1: A podgon is a closed figure with vertices

Def 2: A podgon is a figure with vertices

In order to test these definitions, students have to invent examples which test the definitions, for instance an open figure with vertices. It might be the case that neither of these two definitions match the concept. In that case, they may need to come up with more definitions.

The definition guessing sessions I am exploring in the thesis were the first sessions in both schools. Each was just under an hour long. These were shorter than other sessions in the course because I spent time on introductions and the pre-test.

**The set-up**

As mentioned in the methods section, I will be starting by exploring the beginnings of the two definition guessing sessions in detail.

**Ganga.**

Students were given initial examples (which can be found in the methods section) and asked to develop possible definitions in the course of a class discussion. Pankaj suggested the following:

*Pankaj: A podgon is a closed figure with vertices.*

*I added in an alternative:*



*"A podgon is a figure with vertices."*

*I then gave the following instructions:*

*"Let us start with this (pointing to the definition I added). In your groups, tell me what question you would ask each of us to see whether this is our definition or not… So, what example would you give us."*

*Group 1 gave the following examples: a five-pointed star and a semicircle.*

*Me: "If you found out that this shape (the star) is or isn't (a podgon), how would that help you with the definition?"*

*Anjali: "It would tell us whether the diagonal was inside the figure or outside."*

*Mohanan intervened shortly after to say:*

*"Hold on. Hold on. When you ask a question, you should know what it is you want to find out. So, when you ask the question, 'is a semicircle a podgon?', what do you want to find out?"*

*Anjali: "Whether a curved side counts in a podgon?" [Evaluating a definition]*

*Mohanan: "This one says nothing about curved sides, so you want to find out whether curved sides are allowed."*

*When asked again why they wanted to know about semi circles:*

*Gauri: "We are asking whether it has one vertex or requires more."*



*Rather than continuing with the discussion, we gave them our responses to the*

*semicircle and then went to the other group.*

The goal of my question about the star, and of Mohanan's intervention, was to get students to think through the examples they suggest. It isn't clear whether the group understood the task they were given. Anjali's mention of diagonals seemingly shows that they didn't, given that this particular consideration would be irrelevant to the task of evaluating the definition I gave them. If we listen to their group discussion, the discussion is mainly about diagonals and equal sides (which they did not bring up in the class discussion). While the equal sides discussion is relevant to the session's broader goals, it is much less relevant to the particular task they had been set.

Of course, the examples they gave could help evaluate the definition since if they were not podgons, the definition would be invalid. An alternative set of instructions would involve asking the students to find examples that allow them to choose between the two definitions. The examples they gave do not achieve that since a rejection of those examples would result in a rejection of both definitions. Moving on to the other group, this is how it seems they interpreted the question.

*Me: "Do you have something which deals with the one (pointing at the*

*definition)?"*

*Imran: "Is an open figure a podgon?"*

*Me: "So, I can't answer that. It is too broad of a question."*

*Imran comes up to the board and draws a 'greater than' shape.*



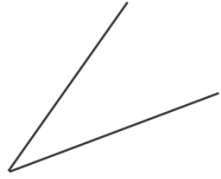

*Figure 8 First counter-example*

*Imran: "Is this a podgon?"*

*I said it isn't for either of the two definitions and asked:*

*"Now that you asked that question, how does that help you with this definition?"*

*Imran: "So, a figure can be open and closed. So, by this, we know it is only a closed figure."*

From this, it seems like at least somebody in the group understood the purpose of the exercise. This conclusion is supported by evidence from the group discussion. This group discussion was relatively short since they spent most of the time figuring out their group roles. Here is a short discussion between Imran and Pankaj towards the end:

*Imran: What I think is we have been given two statements. One is about only closed figures, and the other is about both closed and open. So, we should know about the open.*

*Pankaj: But this isn't correct. We haven't proved that this is correct.*

*Imran: That is only what we have to ask.*

From this, it seems clear that Imran understood the task. He finds the differences between the two definitions, notices that one of them is entirely contained in the other, and looks for an



example outside of the 'closed only' definition but within the more general definition. Given Pankaj's response, it isn't clear what he is referring to by 'this' when he says 'this' has not been proved, and I have not been able to figure that out from context.

Going back to the class discussion part of this, Imran's claim at the end that 'by this, we know it is only a closed figure' doesn't follow from what he had said previously. All the information he had was that there are some non-closed figures which are not podgons. This reasoning from an example to a conclusion is something I will discuss later in the chapter. I will now explore the set-up in Indus.

**Indus.**

Given the two initial examples, the definitions the three groups came up with (in order of reading them out to the class) were:

> *Group 2 (Aditi): Closed figure made up of straight lines which are equal*
>
> *Group 3 (Tarini): A closed figure with finite vertices and congruent straight sides*
>
> *Group 1 (Vandana): Equilateral polygon*

Group 1 and Group 3's definitions are equivalent, assuming by 'made up of' they mean shapes made only of straight lines. This seems to be what they meant since they never brought up non-straight lines in their group discussion. Group 2's is similar but adds that the number of vertices is finite[1]. Listening to their group discussion, they spent a lot of time on this aspect of the definition. I decided to sidestep that point. As I said in my notes:

---

[1] This idea that a circle is an infinite-sided polygon is something I have come across various times amongst students. It is not something I have found in Indian textbooks, but it seems to be commonly and uncritically accepted by



> *"... I did not want to spend time discussing whether a circle is an infinite-sided polygon. While that would have been an interesting discussion, it would have been a distraction from the goals of the session."*

I wrote down group 2's definition:

> *Me: How can we check this definition? What are the various checks we can do?*
>
> *Uday: We can ask Me 1 and Me 2 whether this shape comes under the definition or not.*
>
> *Me: Why are you doing that?*
>
> *Uday: So that we can come up with a modified definition*
>
> *Me: So, you want to modify the definition by seeing things that fit the definition right now, but may not fit the definitions that Me 1 and Me 2 have. So, does anyone have anything?*
>
> *Uday: Can a square be a podgon?*
>
> *Me: Okay. Why are you asking that?*
>
> *Uday: Because it will prove that the equilateral sides*

Asking whether a square is a podgon or not is valuable since if it is not, the definition will need to be modified - not all equilateral polygons would be podgons. However, Uday's reason for asking this seems to be that he wants to 'prove' that podgons are equilateral, which

---

students. I have not been able to find anything on this in the mathematics education literature. While I think that addressing this makes for a good session, and I have done so in other workshops, it takes a lot of time and would have decreased the amount of time I would get to spend on this session.



would be a 'proof' by example. Assuming this is what he meant, Harel and Sowder (2007) would class this under an 'empirical' proof scheme.

I then decided to invert the question and asked:

---

*Me: So, what happens if a square is not a podgon?*

*Tanya: Then it comes down to three sides.*

*I give the judgments on the square. Then, Tarini asks:*

*Tarini: Is a pentagon with equal sides is a podgon or not?*

*Me: Why are you asking that?*

*Tarini: Because what if the sides are odd sides*

*Me: So, first, you are changing the definition right to three, made up of three equal line sides. Now, you are asking if five works. You know that four doesn't work.*

*Tarini: Yes. A pentagon with equal sides.*

*Uday: A regular pentagon.*

*I tell them it is for both Me 1 and Me 2.*

*Tarini (along with others in unison): Odd number of sides.*

*Me: Made up of an odd number of sides (while writing that)*

---

Tanya makes a similar move to what Uday did previously. The claim she is making is that if we find out that four-sided equilateral polygons are not podgons, we can conclude that



only three-sided equilateral polygons are podgons. Tarini notices that while it may be the case that four-sided polygons are not podgons, it doesn't mean that polygons with a larger number of sides are not.

She seems to have in mind a definition that the number of sides is odd and is attempting to evaluate that definition by suggesting an equilateral pentagon. I framed what she was saying as an evaluation of the definition that the number of sides is precisely three. In retrospect, it may have been useful to spend more time on this distinction. It would have allowed me to clarify two strategies to the students:

1. Looking for elements that fit the definition but not my judgment

2. Looking for elements that fit my judgment but do not fit the definition

This ended the first part of the session. We worked through the iterative process of definition guessing for the rest of the session as I will discuss in the next section after comparing the two set-ups.

**Comparing the two set-ups.**

There were apparent pedagogical differences between the way the two sessions began. Two obvious differences were:

1. Starting with a group discussion to search for a definition in Indus but not in Ganga

2. Offering the alternative that the figure might be open in Ganga but not in Indus

While the group discussion does take some time, it would have been valuable. Given that it is the most straightforward exercise in the session, most students have a chance to get involved. Also, in Ganga, it seemed like one of the groups did not initially get a sense of what the session was about. A group discussion followed by a class discussion on the results could have helped clarify that before jumping into more challenging parts of the session.



On the closed figure vs. figure exercise in Ganga, I had intended to start with that in Indus. As I mention in my notes, it completely slipped my mind. It did help clarify the goals of the session in Ganga. However, given that the Indus students had already been through a similar session in a previous workshop, they already seemed to have a better initial grasp of the task.

**General observations and reflections**

In this section, I will focus on some general observations and reflections on the module. I will be using data from the set-up and across the rest of the module. When using data from the rest of the module, I will present it in the required context.

### Reasons provided by students.

One of the critical aspects of the task is that students think through why they are asking about a particular example. I attempted to achieve that in the session by asking them their reasons for positing a particular example. We have seen examples of how students responded to this. When tasked with evaluating the definition that allowed for an open figure, Anjali suggested a star and said it would help her figure out something about the diagonals. She then went on to suggest a semicircle to figure out whether curved sides are allowed. As mentioned above, these examples could help evaluate the definitions but not for the reasons given.

Imran's reason for suggesting an open shape in the same situation and Uday's reason for suggesting a square to evaluate the equilateral definition involved flawed reasoning, which I will explore in the next section, but were at least relevant to the definition they were evaluating.

Imran later asks about whether an ice-cream cone is a podgon. The reason he gives is that rejecting it would be:

> *"a counterexample to Mo's. It has vertices and is a closed figure but is not a podgon."*



Similarly, Tarini suggests a rhombus choose between whether a podgon is equilateral or regular. These are two examples of good reasons given for examples students posited.

However, there were many situations, in fact, a majority, where students gave examples, and I did not ask them to give reasons. This included Tarini giving the following example with three straight lines and a curved line:

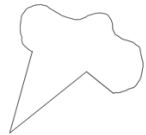

*Figure 9 Three straight lines and a curved line*

This allowed Tarini to come up with the definition:

*"closed figure with an odd number of straight-line sides, which may have a curved side."*

In this case, given that she could come up with this definition given the example, she would probably have come up with a good reason for her example. However, there were other situations where that is less clear.

For instance, just before this, Uday had asked whether an ice-cream cone is a podgon. This was just after I had rejected the following:

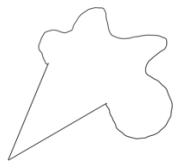

*Figure 10 Two straight lines and a curved line*

I am unable to come up with a coherent reason for what he could have had in mind for a definition that allows us to reject the above, but accept this:



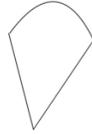

*Figure 11 Ice-cream cone*

It would have been valuable to have asked him.

There seem to be two crucial variables along which reasons students give occur. The first is relevance. I have given examples of situations where students give reasons which are relevant and which are not. The other variable is the validity of the argument.

A pedagogical strategy that could address the first of these (relevance) would be to ask:

a. What is the consequence if X is a podgon?

b. What is the consequence if X is not a podgon?

This can be done in the form of a worksheet, which groups would have to fill every time they posit an example. I will discuss the validity of the reasons given as part of the next subsection.

**Unwarranted conclusions.**

We have already seen a few examples of students arriving at unwarranted conclusions. For instance, Uday concluding that knowing that a square is a podgon would 'prove' that the definition of podgon is 'equilateral polygon.' We also had Tanya seeming to conclude that me rejecting a square as a podgon would limit podgons to having three sides. The third example we have seen is Imran concluding that rejecting a single open figure can allow us to conclude that podgons are closed. I am going to explore that a little further:

*Imran comes up to the board and draws a 'greater than' shape.*



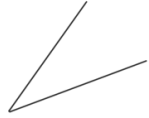

*Figure 12 Greater-than shape*

*Imran: "Is this a podgon?"*

*I said it isn't for either of the two definitions and asked:*

*"Now that you asked that question, how does that help you with this definition?"*

*Imran: "So, a figure can be open and closed. So, by this, we know it is only a closed figure."*

*Me: "Uh. Do you know that? We don't know whether this is a podgon or not (drawing an open arc)."*

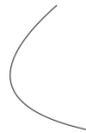

*Figure 13 Open curve*

*"You know that this open figure is not a podgon (pointing at the one he had originally asked about)."*

*Imran: "That one doesn't contain any vertices."*

If the shape I suggested turned out to be a podgon, it would be a counterexample to both the definitions since it doesn't have vertices. However, if it turned out not to be a podgon, it would not be a counterexample to either of the definitions since neither of the definitions



suggests it is a podgon. It seems like this is what Imran was getting at. He was looking for shapes

that fit the general definition but not the specific definition. I do think it would have been

beneficial here to have explicitly stated this in the class.

Instead, what I did do was to offer another figure, which was convincing:

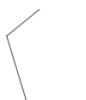

*Figure 14 Open figure with straight sides*

I could not find even one instance of the opposite type of unwarranted conclusion, where

students do not reject a definition even when presented with counterexamples. The closest I have

been able to find is students not spotting a counterexample already on the board. An instance of

this was when I had to point out that a semicircle not being a podgon for me was a

counterexample to the definition that all closed shapes with vertices are podgons. However, once

they saw the counterexample, there was nobody who suggested the same definition again.

**Student feedback.**

In this section, I have attempted to give an overview of the definition guessing session in

the two schools, highlighting the ways in which the students engaged with the session and the

ways in which I (and Mohanan in Ganga) played a role in shaping the direction of the session.

First focusing on the students, from the post-course surveys in the two schools, it seems

like they enjoyed the session. 7 out of 16 Indus students and 4 out of 10 Ganga students

explicitly mentioned this session when asked what their favorite part of the course was. No

student in either school mentioned this session in answering what their least favorite part of the

course was.



However, enjoyment doesn't necessarily correlate with learning. While it is hard to tell whether significant long term learning took place given the data I currently have, I have highlighted instances above where students had to think carefully and express themselves.

**Instructor choices.**

There were many choices made by me during the session which guided what the session looked like. The decision on whether a particular task ought to be a class discussion or a group discussion is something I have highlighted above when comparing the two sessions. Since these sessions are not just driven by the instructor but also by student responses, many of these decisions were made on the fly. Some guidelines for making such decisions would be valuable. However, what those could look like is unclear to me at this point. Further trials of such sessions involving other instructors would be useful to achieve build these guidelines.

Similarly, there were certain paths I chose not to take or to delay even when students were pulling in that direction. One I have highlighted above is avoiding the 'finite vertices' discussion. Another was when I delayed discussing the first definition Pankaj gave in Ganga and asked students to explore my definition. In the first case, the reason I did that was prior experience with similar discussions – I knew how long such discussions can take. In the latter case, I wanted the students to work on a simpler example for which I was quite certain they would easily find a counterexample.

**Adding rigor.**

The possible changes to the sessions I have discussed so far can be summarized as making the session more rigorous. This includes having worksheets where students have to explicitly write out their reasoning, working through the reasoning on the board more carefully,



and carefully taking students through the various types of conclusions they can draw given some information.

There is a slight tension here since one of the goals of this session is as an icebreaker aimed at getting as many students involved. The way I would address these competing values going forward would be to have two facilitated definition guessing sessions, one at the beginning and the other in the middle of the course. The first one would maintain a similar level of rigor to the sessions done as part of this workshop. The second would be more rigorous.



## Chapter 6 – Triangle theory building: Assumption digging & integration

This and the next chapter are focused on the second course module – triangle theory building. A sketch of a lesson plan for this module is available in Appendix A and the first subsection of this chapter describes the lesson plan and the learning outcomes it is aiming at.

After reiterating the mathematical content and the learning outcomes of this module, the bulk of this and the next chapter focuses on a local analysis of the implementation. My analysis involves going deep into particular parts of the implementation of this module in the two schools, dealing both with students' engagement with particular aspects of the course, addressing the first research question, as well as with how I facilitated the sessions, addressing the second research question. To achieve this, I have attempted to interpret what students have said given the direct context along with their group discussions and my teacher notes. I have also attempted to suggest alternative ways in which I could have facilitated particular parts of the course.

In this chapter, I discuss the set-up of the module. In this case, by the set-up, I mean the first day of the module in the two schools. I present the entire set-up in both the schools in a sequential manner.

The next chapter discusses the implementation of the rest of the module and student feedback. It ends with a discussion of general observations and reflections from the implementation in the two schools.

Here is a table of contents for this chapter:







**Description of the module**

The goal of this module is for students to construct a theory of triangles. The idea is to write out a large number of claims we believe to be true about triangles and connect them by justifying some of them using others. So, there are two important aspects of the module: integration and assumption digging. Integration refers to the idea of placing these claims about triangles into a coherent body of knowledge with a small number of starting points from which you can deduce the rest. Assumption digging refers to the process of asking 'why should I believe that?' and 'what does this mean?' in a recursive manner.

The intended learning outcomes of this module include an understanding of assumption digging and integration. Apart from that they include an understanding of:

1. definitions

2. classification

3. undefined entities and axioms

4. proofs



In terms of the questions of theory building, this module involves all the questions of theory building described in Chapter 3 apart from the one to do with extending definitions.

**The set-up**

As I mentioned above, by the set up, I mean the first day of the module in the two schools. All the mathematical interactions in the class discussions on that day are presented. I also attempt to interpret what students have said and try to take learnings for future iterations of the course.

### Eliciting statements in the two schools.

In both the schools, I asked students to come up with 'things they know to be true about triangles' in their groups. I took examples from the groups one by one and ended up with 10 claims from each school.

Ganga:

1. Side opposite to greater angle is greater

2. Area of a triangle = in-radius x semi perimeter

3. Sum of two sides > third side

4. Sum of angles is 180 degrees

5. Pythagoras Theorem

6. Angles opposite equal sides are equal

7. Area = $[s(s-a)(s-b)(s-c)]^{\frac{1}{2}}$

8. If all the sides are equal, then it is an equilateral triangle

9. Circumcenter is the intersection of perpendicular bisectors

10. Angles of an equilateral triangle are equal

Indus:



1.  Sum of angles is 180 degrees

2.  Sum of 2 sides is greater than 3rd side

3.  4 types: equilateral, right angled, scalene and isosceles

4.  Congruence by: SSS, SAS, ASA, SAA, AAS and RHS

5.  Sum of exterior angles is 360 degrees

6.  Angle opposite to greatest side is greatest

7.  Side opposite to greatest angle is greatest

8.  Angles opposite equal sides are equal

9.  Perimeter is sum of sides

10. Area = ½ bh = $[s(s-a)(s-b)(s-c)]^{½}$

While there were a large number of similarities in what the two groups gave, there are also notable differences. In Indus, the mention of the different types of triangles and congruence stand out while in Ganga, there is a mention of equilateral triangles and Pythagoras' theorem.

However, in both the schools, some of the simplest things they know about triangles were not mentioned: that they have three sides, that they are polygons, that they have three angles, and so on. A disadvantage of students giving the initial list of claims is that it may be harder for an instructor to prepare for the future path of the sessions in advance. The list given by the students resulted in two different starting points for the two schools.

**Evaluating claims.**

As I say in my post session notes, I was looking for the simplest claims to start with which would also be telling. So, I decided to choose equilateral triangles in Ganga. In Indus, while I did discuss equilateral triangles later in the first session, I started with the classification of triangles since the students had not mentioned anything about equilateral triangles.



***Ganga.***

The first claim we took up in Ganga was:

Angles of an equilateral triangle are equal

*Equilateral Triangles.*

I asked the students to justify this claim, attempting to use other statements on the board in their justification. This is given as a group discussion task.

I ask the middle group to give their justification and Arnav came to the board. He drew the following and said:

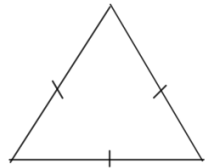

*Figure 15 Equilateral triangle*

*Arnav: We have an equilateral triangle. All the sides are equal and so this is also an isosceles triangle. So, this angle and this angle are equal (pointing at the two bottom angles). If we consider this vertex (pointing to the bottom right vertex), this angle is equal to this angle (pointing at the two angles on the left), and all the angles are equal.*

*Pankaj (interrupting): Angles opposite equal sides are equal.*

*I then write up what Arnav said:*

*Equilateral triangles have equal sides*

*I then label the vertices A, B and C and continue:*



*Angle A = Angle B since CA = BC*

*Angle C = Angle B since AC = AB*

*So, Angle A = Angle B = Angle C*

*Me: Just because Angle A = Angle B and Angle B = Angle C, why does Angle A = Angle C?*

*Vivaan: If we say that if something is equal to something else and something else is equal to something other than that, then the original thing will equal to that. If we work like that, then this follows. [Laying down axioms]*

*Vivek: A is equal to B, so in the place of B you can put A.*

I decided to accept that and move on. It seems like I did not parse what Vivaan had said in the moment. In retrospect, it appears to be quite insightful. At the end he says 'if we work like that, then this follows.' In other words, he seems to be saying that if we assume transitivity of equality of angle measure, then this conclusion follows. This would have been valuable to discuss further as this was the first instance of a student appealing to an axiom.

*An open question.*

I then picked the first statement in the list: sides opposite the greater angle is greater. I asked them to either prove that statement or use the statement to prove something else. The second option I gave them made the task quite open. They could prove anything which involved this statement in the proof. I left this as a group discussion task.

I invited Anya from the group on the right to come up to the board to present her group's work. She draws a circle with its diameter and says:



*Anya: The angle subtended by the diameter of a circle is 90 degrees… Any point on the circumference, the angle made is 90 degrees. But, if we take any random chord, like join any point, then would subtend a lesser angle than ninety degrees, an acute angle. The diameter is the longest chord in the circle and so the angle is larger.*

*Here is what she had drawn:*

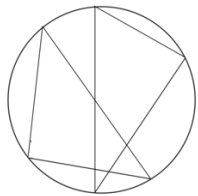

*Figure 16 Circle with two triangles*

*Me: But is that always the case. Let me take this side and this side and draw a triangle.*

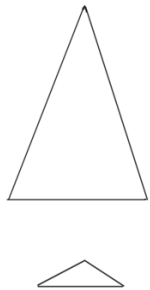

*Figure 17 Two triangles with different bases and angles*

*This side is longer than this side, but which has the larger angle*

*Anya: No, but in this case.*

*Me: Why is this different?*



The apparent problem here with Anya's reasoning is that she used two different triangles and concluded that since a side of one triangle is longer than a side of another triangle, the angle subtended by the shorter side must be smaller. In the examples I drew, you can see that is not necessarily the case. The relative sizes of two triangles doesn't affect the measures of specific angles. However, it could be the case that she meant something else when she stated the claim. That is what I probed into next:

> *Me: Oh, so you are saying something which goes through the diagonal? Goes through the center?*
>
> *(Many in the audience say yes)*
>
> *So, if it goes through the center, then this angle will be less than this one.*
>
> *Anya: No. I'm saying that this (pointing at the chord) will subtend a lesser angle than the diameter.*

I suggested an alternative formulation to Anya's claim which was that if we subtend an angle from a chord which creates a triangle such that one of the sides of the triangle is a diameter, then the angle subtended by the chord on the boundary will be less than the right angle in the same triangle. This is something which can be justified using the fact Anya stated that the diameter is the longest chord along with the claim I had asked them to use. The reason I suggested this claim is that in the picture Anya drew (above), it looks like the two lines intersecting in the middle were diameters.

However, that Anya said no seems to confirm that she was thinking about two separate triangles with two different sets of side lengths. I then decided to challenge the claim she had



made about the relative angle measures by drawing a chord subtending an angle in the minor arc. In that case, what you get is an obtuse angle.

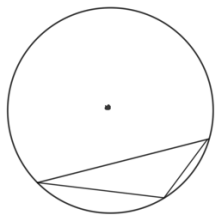

*Figure 18 Circumscribed triangle*

She accepted the problem with the claim, but I am not sure whether the main problem with her earlier reasoning was clear to her. The problem was not that the conclusion was false – you can specify the claim enough so that it is true. The problem was that she was using the claim about relative angle sizes outside of its scope. That claim is about sides and angles of the same triangle and not different triangles. I did attempt to explain that straight after:

> *Me: Notice that just because the side is bigger, doesn't mean that the angle opposite is bigger. That is only in the same triangle.*
>
> *I then changed the original claim on the board to reflect that by adding in the words 'in a particular triangle'.*

An interesting thing to note is that how Anya stated her argument, she seems to have used 'angles opposite greater sides are greater' rather than 'sides opposite greater angles are greater'. The reason I'm saying this is that Anya ends with 'and so the angle is larger'. I did not seem to have spotted that during the session. Listening in to the group discussion, Vivaan was talking about concluding that the diameter is the largest chord given that the angle it subtends is the largest. So, there either seems to be some miscommunication within the group or Anya didn't communicate clearly when she came to present her argument. This would have been a valuable



thing to have spotted since I could have used the opportunity to discuss the distinction between the two claims.

I continued with the group on the left. Vivek spoke for the group. I drew a right triangle on the board since I had talked to the group about their claim while they were discussing it.

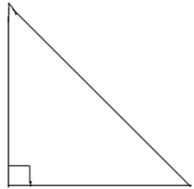

*Figure 19 Right angled triangle*

*Vivek: That angle will be 90 degrees and it will be greater than the other two angles because the other two angles will sum up to 90 degrees, so they will be less than right angle. And in a right angled triangle, the hypotenuse is the greatest side.*

*Pankaj: The side opposite to the right angle.*

*Vivek: So, the greatest angle will be opposite the greatest side.*

There are at least two ways you could interpret what the group said. One is that they are concluding that the hypotenuse is the longest side. However, if you look at what they said, Vivek concludes by saying 'so, the greatest angle will be opposite the greatest side.' This makes it seem like they were concluding that the side opposite the greatest angle is the greatest in the case of a right triangle. This is the interpretation I seem to have had in the moment:

*Me: So, you are saying that in the case of a right triangle, the side opposite the greatest angle is the greatest because you know that the hypotenuse is the*



*greatest side, and the angle opposite the hypotenuse is the greatest angle… all*
*this together gives us that the side opposite the greatest angle in a right-angled*
*triangle is the greatest. This doesn't tell us that for all sides. It doesn't say that if*
*angle C is smaller than angle A, then AB is smaller than BC. It doesn't give us*
*this.*

It would have been interesting to pursue this further later in the course by asking the students how they know that the hypotenuse is the largest side of a right-angled triangle without appealing to the claim that sides opposite to greater angles are greater.

I then moved on to the third group. They claimed that they had done the same thing as the previous group. Since I had talked to them during their group discussion, I knew that they hadn't. Rather, they had justified that the side opposite the right angle (the hypotenuse) is the greatest side of a right-angled triangle. Rather than justifying the claim in the particular situation, they had used the claim to come to a conclusion.

After I pointed the distinction out, I asked the group to state what they had come up with. Arnav responded:

*Arnav: What we were trying to do is, like, we used that, the side opposite, and*
*we were proving that hypotenuse is the greatest side.*

*He then came to the board to justify this:*

*Arnav: We know this angle is 90 degrees and the sum of angles is 180 degrees.*
*So, these two angles have to be smaller than 90 degrees. So this (pointing at the*
*right angle) becomes the largest angle and because this is the largest angle, the*
*side opposite is the hypotenuse, so it is the largest side. [Providing proof]*



Arnav's initial response appears to show that he can see the distinction between the two conclusions and arguments – the difference between justifying the claim in the context of right triangles and using the claim to justify something about right triangles.

At that point, as I mention in my notes, I decided to pick a claim for the students to justify which would require them to write out some of the more basic things about triangles they had initially missed, such as they have three sides, three angles and so on.

I picked something which was used in the proof above: the largest angle of a right triangle is the right angle.

*Anya: The angle sum property of a triangle is 180 degrees. In a triangle any angle cannot be zero degrees. So, even if one angle is as small as 1 degree, the other angle is 89 degrees, so 90 degrees is the largest angle.*

*I wrote:*

*The sum of angles of a triangle is 180 degrees*

*The right angle has measure 90 degrees*

*So, the sum of the two other angles measures 90 degrees*

*Since an angle of a triangle has to be greater than zero, both are less than 90 degrees*

*Hence, the right angle is the largest.*

*I then ask:*

*Me: How do we know there are only two other angles?*



*Multiple students: Triangles have only three angles*

*I added that to the list of claims along with the claims that triangles have three vertices and three sides.*

The point of doing this was to demonstrate the need for stating the most obvious things in order to be rigorous in justification. An alternative way to do this would have been to say that I rejected their proof even though I accepted that the sum of angles of a triangle is 180 degrees since if one angle is 90 degrees, the other one is also 90 degrees. The presumption in that argument is that a triangle has two angles. In order to show that I was wrong, they would have had to explicitly state that a triangle has three angles. Coming from them seeing a need for stating this rather than from me would possibly have made more of an impact.

*Equivalent definitions and more equilateral triangles.*

I then returned to a discussion of equilateral triangles. I asked if we could use equality of angles as the definition of equilateral triangle. There was no response, so I asked:

*Me: What makes two definitions the same?*

*Vivaan: In this case, angles of an equilateral triangle are equal if and only if the sides are equal.*

*Me: Yeah. You want this to imply this and this to imply that (pointing at the two definitions in turn). And then you get that these are equivalent definitions.*

I had discussed equivalent statements with some of the students in this group in previous workshops I had done. It seems like at least Vivaan seems to have remembered something from that.



I decided to discuss equivalent definitions in a slightly different context, that of even numbers. I said that I had sat in on a fifth grade class a few years previously where the teacher had told the students that the definition of even number is something which ends in 0, 2, 4, 6 or 8[2]. I asked whether that is a definition of an even number and took them through why it is equivalent to the usual definition that even numbers are numbers divisible by 2. I do think it is useful to bring in other areas of mathematics so that students see the similarities and can bring their understanding of other areas of mathematics to bear on what they are working on currently. However, this could probably have been more impactful if I had asked students in groups to justify that the two definitions of even numbers are equivalent. This could have also led to a useful discussion amongst students of which definition is better and why.

I then went back to discussing equilateral triangles:

> *Me: You showed that if all sides are equal, then all angles are equal. What about this way? (drawing an arrow from the statement all angles are equal to all sides are equal).*
>
> *The students discussed this in their groups and Vivaan comes to the board.*
>
> *Vivaan: If you take an equilateral triangle like this, this angle is equal to this angle is equal to this angle… therefore, this side and this side are equal because they are equal angles, and this side and this side are equal. And that's why… things equal to the same thing are equal to one another, all sides are equal.*

---

[2] This is something I have seen in multiple schools in India. Students seem to conceptualize even numbers as those which end in 0, 2, 4, 6, or 8. While this is equivalent, it is dependent on working in base 10 and doesn't really give the same understanding of even numbers as defining them as numbers divisible by two or numbers you can break into two equal parts.



> *Me: What you have added in is 'sides opposite equal angles are equal.'*

During the group discussion which preceded this in Vivaan's group, I talked to them about this argument they had come up with. In that discussion, he explicitly mentioned that 'sides opposite equal angles are equal' which he did not do while at the board. Rather than me adding this in in front of the class, I could have written out Vivaan's proof and asked students what was missing. Even though I had done that while talking to the group, I had not done that with the entire class.

I will now describe the first day at Indus before discussing the similarities and differences between the implementation in the two schools.

**Indus.**

In Indus, I started differently. One of the reasons for doing that is that students hadn't mentioned the claims about equilateral triangles the students in Ganga had. I started with their classification of triangles. They had said that triangles are of four types: equilateral, right angled, scalene and isosceles.

*Classification.*

The problem with this listing of types of triangles is that there is no consistent criterion. A few of isosceles, scalene and equilateral have to do with relative side length which right angled has to do with the measure of the largest angle.

I started writing the four down as a tree with triangle being the parent node and all four being direct children. While I was doing that, Tarini said:

> *Tarini: See, there are two types of classification, one is to do with the sides and another is to do with the angles.*



*Me: Yeah. So, let's stick to the sides for now. What are the different types?*

*Multiple students: Equilateral, isosceles and scalene.*

*I drew a tree with triangle as the parent node and the three types of triangles directly under it.*

*Me: Now if we prove something about isosceles triangles – that angles opposite equal sides are equal, then it's not true about equilateral triangles?*

*Multiple students: It is*

*Me: But equilateral triangles are not types of isosceles triangles*

*Multiple students: They are*

*Devika: Equilateral triangles are a type of isosceles triangle but every isosceles triangle is not an equilateral triangle.*

What I was trying to do here was to discuss logical inheritance – if we know something to be true about isosceles triangles and we classify equilateral triangles as isosceles triangles, then that property is inherited by equilateral triangles. However, I phrased it in a misleading manner. It is not that it is impossible for the same thing to be true about the two shapes even if we choose the other classification. Rather, the issue with the first classification is that we would have to prove the same claim separately for equilateral triangles.

*Me: So, you are saying you want this*

*I remove the arrow pointing to equilateral triangles from triangles and add in an arrow from isosceles triangles to equilateral triangles.*



*Me: There is a difference between these right?*

*Multiple students (lead by Tanya and Tarini): You have scalene then isosceles, then equilateral.*

*This was interesting. The students had extended the need for subclassification to scalene triangles.*

*I clarified that to make sure what they meant by drawing the following tree:*

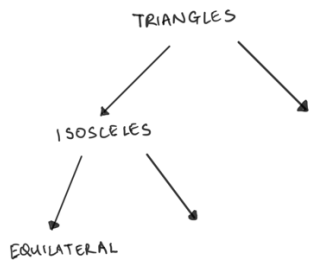

*Figure 20 Triangle classification 1*

*I asked if they wanted scalene to be written under the arrow on the right. Many of the students in the class said that was what they wanted. However:*

*Tarini: No*

*Tanya: Before isosceles, there is scalene*

*Me: So, isosceles triangles are types of scalene triangles*

*Multiple students from the group (including Tarini and Tanya): Yes*

The consequence of what Tarini and Tanya seem to be saying is that scalene triangles and triangles are the same thing assuming their concept of scalene triangle includes the usual notion,



which is all triangles which are not isosceles. If isosceles triangles are now included in scalene triangles, all triangles become scalene triangles. Assuming scalene triangles are types of triangles, the consequence is that scalene triangle is a synonym for triangle. I probed into that:

*Me: So, what is the meaning of scalene triangles?*

*Uday: A random triangle is a scalene triangle.*

*Me: But an isosceles triangle is a scalene triangle*

*Multiple students: Yes*

*Me: So, all triangles are scalene triangles*

*Multiple students: Yes*

*Me: So, what is the difference between the word triangle and scalene triangle?*

*Multiple students: Nothing*

*Tanya: All triangles do not come under isosceles or equilateral. So, there should be another category.*

*Me: That's just a triangle*

*Tanya: But, every triangle is not isosceles or equilateral*

*Me: Every triangle is a scalene triangle, right? So then, basically you have this: every triangle is a scalene triangle.*

*Tanya: But, every scalene triangle is not an isosceles or equilateral.*



*Me: Yeah. But now you have (and I draw this:)*

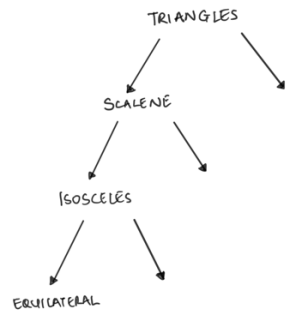

*Figure 21 Triangle classification 2*

*Me: Notice that there is nothing here (pointing at the top left arrow). This is empty. Are you saying that all triangles are scalene triangles?*

*Multiple students: Yes*

*Me: So, basically, these two are the same (scalene and triangle). There's nothing else here.*

*Tanya: If we say isosceles and equilateral are types of triangles, there will be nothing for the triangles which are not isosceles or equilateral.*

*Me: They can come here, right? And I draw the following:*

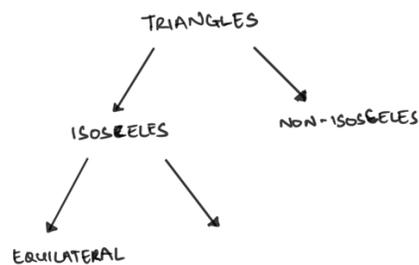

*Figure 22 Triangle classification 3*



*Devika: That is scalene*

*Me: Now you are saying something completely different. This is very confusing. So, what are scalene triangles?*

There is clearly a lot of confusion in the above exchange. There are two seemingly contradictory claims being made by the same students. Tanya and Tarini initially seem to be making a parallel argument to the one made earlier to conclude that isosceles triangles should be types of scalene triangles. The result of this is that triangles and scalene triangles become synonymous. However, this appears to be a problem at least for Tarini, who said, 'there will be nothing for triangles which are not equilateral or isosceles.'

In my notes, I wrote that I didn't know what to make of this. One possibility is that the students did not have a clear understanding of the representational system being used. I didn't think of that in the moment since I had used the same representation in previous workshops where Tarini and Tanya were present. Another possibility is that there were competing values in their mind – logical inheritance on the one hand and having a term to describe a class on the other.

To give more evidence that there were these competing values, as the discussion continued, I drew the following Venn Diagram:

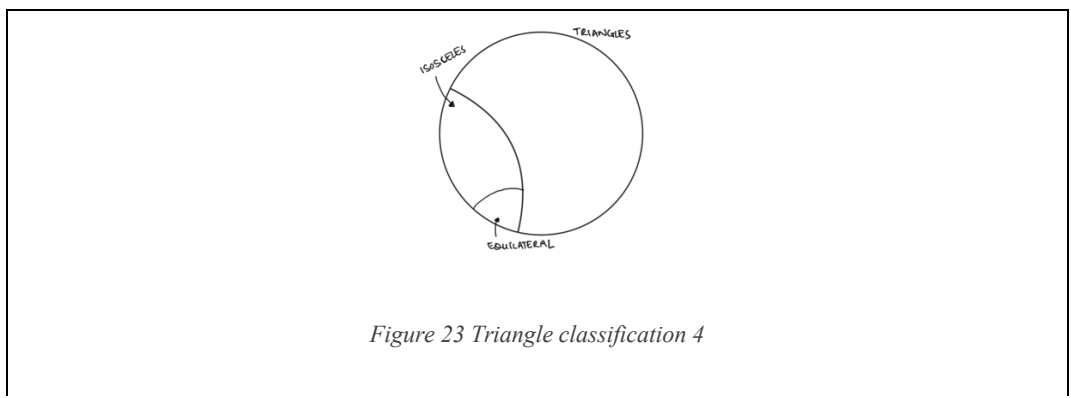

*Figure 23 Triangle classification 4*



*Tanya: We should have a name for the part remaining. So, that is scalene.*

*Me: But you said that is not scalene. You said the whole thing is scalene.*

*Tanya: Yes.*

*Me: So, you're saying even this is scalene (pointing at the isosceles part of the diagram)*

*Tanya: Yes.*

Given that I was now using a different representation, one they had been introduced to in their school mathematics, it makes it less likely that the issue is the representation and more likely that the confusion stems from these two competing values – description and logical inheritance. However, it could be the case that even though they have used Venn diagrams before, they are not comfortable with them.

I left that discussion there and asked the class to classify triangles on the basis of angles in their groups. I asked for what they had come up with and here is what I drew on that basis:

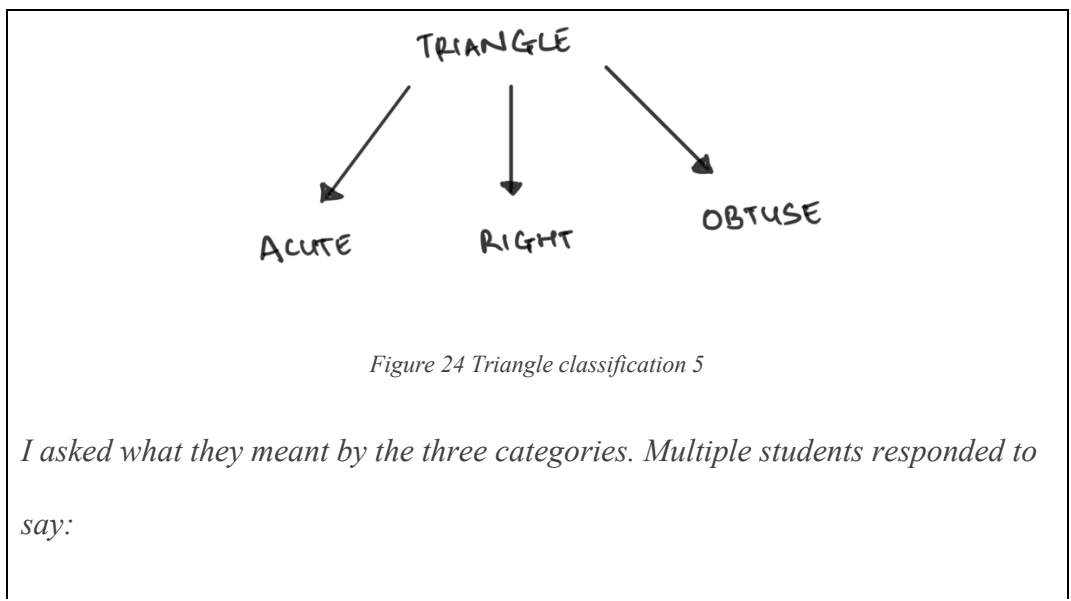

*Figure 24 Triangle classification 5*

*I asked what they meant by the three categories. Multiple students responded to say:*



*Obtuse: One angle greater than 90*

*Acute: All angles less than 90*

*Right: One angle 90*

*Me: You're saying that this is a complete classification, right? By that I mean, you're saying that everything falls into one of these categories? In fact, exactly one of these categories. It is not just complete but exclusive. Everything falls into one of these and all triangles are of one of these types.*

*Multiple students: Yes*

There were significant issues with my initial framing of completeness. I conflated completeness and exclusivity at first. Even when I separated them out, I didn't expound on them sufficiently. An alternative strategy would have been to demonstrate these using classificatory systems outside of mathematics. For instance, consider the classification of food by taste. Suppose for now that there are only two types of tastes we are interested in: spicy and sour. There are foods which are spicy and not sour, foods which are not spicy and sour, foods which are spicy and sour and foods which are neither spicy nor sour. This classificatory system is non-exclusive. However, we could convert this into two exclusive classifications which separate out the dimensions of spice and sourness (I am ignoring the complexity of degrees of spice and sourness for the purpose of this example). In this same example, if we only had the categories sour and spicy without the possibility of having foods which are neither, that would be an example of an incomplete classification. So, if you wish to have a complete and exclusive classification, you will need all four categories.



I decided to discuss how they know that the classification of triangles by angles is complete:

> *Me: How do I know that? How do I know that there can't be something with two right angles? You see that that doesn't fit into any of these, right?*
>
> *Multiple students: What?*
>
> *Me: So, if there was a triangle with two right angles.*
>
> *Uday: If we do that, we will be breaking (and points at the list on the board)*
>
> *Tarini: Angle sum property*
>
> *Me: Yeah. What I'm trying to say is that a triangle with two right angles, you do not know where it is.*
>
> *Udit: But a triangle with two right angles is not possible.*
>
> *Me: But we don't know that yet. If you know it's not possible, then this works.*
>
> *Tanya: It is possible in spherical geometry*
>
> *Tarini: That's what I'm saying. It's not possible in Euclidean geometry (says something else which is not decipherable)*
>
> *Me: So there are places where it is possible, right? But, right now since we are sticking to Euclidean stuff. In Euclidean geometry you are saying it is not possible. In Euclidean geometry you cannot have two right angles.*
>
> *Tarini: Yes.*



*Me: But that doesn't come automatically from the definition. So, you actually need a proof to show that this is complete. Otherwise, it may not be complete. Maybe there could be a triangle with two right angles and that doesn't fit into any of these two categories. We know the proof is very easy. What's the proof?*

*Uday: The sum of all angles is 180*

*Me: Okay. So, can't you have two right angles?*

*Uday: There has to be three angles*

*Me: Yeah. One is zero.*

*Multiple students: It will not be a triangle.*

*Me: Aah. So, triangles must have non zero angles (and I add that to the list of claims. I also add that there is at most one right angle in a triangle)*

*Uday: Also, at most one obtuse angle (I add that in)*

*I explicitly write out the proof of the right angle claim:*

*Sum of angles = 180 degrees*

*If two angles are 90 degrees each,*

*The third is 0 degrees*

*Which contradicts xii (the claim that all angles need to be greater than 0 degrees)*

*Me: This shows that the classification is complete*



As I say in my notes, the point of going through this exercise was to get students to realize the assumptions required to even create the classificatory system. Retrospectively, I could have done this first in the earlier classification to do with sides. In that classification, the assumptions made were minimal – if we make the assumption that triangles have three sides, the only possibilities are that no sides are of equal length, two sides are of equal length and three sides are of equal length.

However, in this case, the assumptions are much deeper. In fact, as Tarini and Tanya said, the classification requires assumptions specific to Euclidean Geometry since this classification is not complete in Spherical Geometry.[3]

At the end of the conversation, I said that 'this shows the classification is complete.' By this, I was referring to the proof that triangles with two right angles cannot exist. That isn't the case. This is only one of the things we would need to show in order to conclude that the classification is complete.

After that discussion, I pivoted to asking why the classifications we had seen were good classifications. I gave the following example:

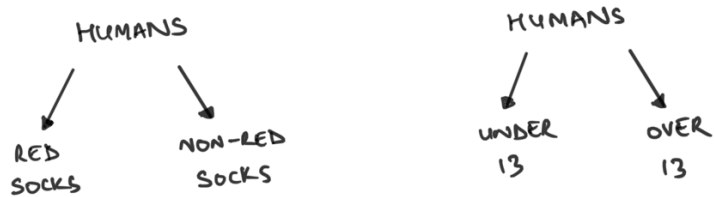

*Figure 25 Classification of humans*

---

[3] In my notes, I say that I talked to them after class about Spherical Geometry. During that conversation, I asked them where they had heard of spherical geometry and they said I had talked discussed it in the workshop I had done a few years previously.



> *Me: We can divide humans into two types: red sock and non-red socks. We are*
>
> *doing biology. The other classification has humans under 13 and over/equal 13.*
>
> *Tarini: These are categories on which we classify?*
>
> *Me: Yeah. Now say we are doing biology and we want to find out about humans.*
>
> *Which of these is a good classification?*
>
> *Multiple students: Second*
>
> *Me: Second. Yeah. Because the second one tells you something, humans*
>
> *development at the age of 13. There is a difference between these two (under 13*
>
> *and over 13). But in this case, it is very unlikely that there is something true*
>
> *about, ok, people who wear red socks are over seven feet tall. Very unlikely that*
>
> *such a relationship exists.*

The goal here was to give students classifications which illustrate the distinction between those which have value and those which do not depending on what we are interested in. If we are interested in biology, it is unlikely that things humans choose to wear will have significant correlations. However, age, especially near adolescence, will have significant correlations. I could have spent more time clarifying the relationship to the domain. For instance, it is plausible that there are psychological correlations to those who choose to wear red socks!

I continued by moving on to the angle-based classification and said:

> *Me: I can add another category of triangles with one angle equal to 47.5*
>
> *degrees.*
>
> *Uday: That will come under acute angle.*



*Me: No, you can have a triangle with one angle 47.5 degrees which is a right triangle. You can even have a triangle with one angle 47.5 degrees which is an obtuse triangle. Right?*

*I add in the category 'One angle equal 47.5 degrees' and Tarini says:*

*Tarini: We need to know all the three angles.*

*Me: Here you don't know all the three angles (pointing at the right-angle category). Here you only know one angle.*

*Tanya: But if we know all the angles, then we can decide.*

*Me: Yeah. But, this is an alternative classification. Let me add in the exact same thing as with 90 degrees. One angle equal to 47.5 degrees, two angles 47.5 degrees, two angles greater than 47.5 degrees, two angles lesser than 47.5 degrees. Everything's covered by this. Why is this a bad classification while this is a good one.*

*Uday: Most probably we won't have a triangle which is 47.5 degrees.*

*Me: Most probably, you won't have a right angle triangle. So, what is it about right angled triangles?*

*Multiple students speak at the same time and it is not decipherable.*

*Me: Yeah. If you have a right triangle, there are many things you can tell me about it, right? What's true about the relationship between the sides?*

*Prerna (along with others): a square plus b square equals c square*



> *Devika: The square of the hypotenuse equals the sum of the squares of the other*
> *two sides.*
>
> *Me: So, you know a lot of things about right angled triangles. Not just that. The*
> *area of a right triangle is that multiplied by that (the sides) divided by two. You*
> *know things about this category (the right angled one). This is useless (the 47.5*
> *degree one). You know nothing about 47.5 degree triangles which you do not*
> *know about 50 degree triangles.*

Triangles which have a single angle equal to 47.5 degrees is not a useful category since it does not have interesting correlations, just as the category of humans who wear red socks in the domain of biology. On the other hand, triangles with one right angle has a large number of interesting properties like pythagoras' theorem.

One potential pedagogical issue here was to do with the board work. If you look at what we ended up with, I used the same parent node for both the classificatory systems we had. This could be confusing since it is unclear whether I am just adding to the existing classification or positing a new one. I am not able to tell from the interaction as to whether that resulted in confusion.

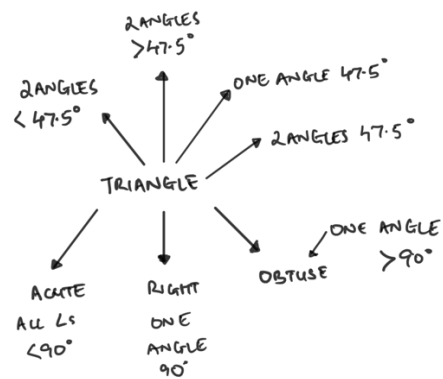

*Figure 26 Classification of triangles 6*



*Equilateral triangles and number of sides.*

I then asked the students to name things they know about acute angled triangles, obtuse angled triangles, equilateral triangles and isosceles triangles, and we discussed them. I added the following on the board:

---

*Equilateral triangles are acute*

*Sides of an equilateral triangle are equal*

*Angles of an equilateral triangle are acute*

*Me: This one. Equilateral triangles are always acute. How would you show that?*

*Uday: Because all angles are 60 degrees*

*Me: All angles are equal*

*Tarini: 3x is 180. So, x is 180 by 3, which is 60*

---

I write out the argument on the board and accept that. The argument used the fact that the angles are equal. So, that is what I asked them to justify next in their groups. I told them to use the stuff on the board as far as possible. After they had discussed in groups, I asked the middle group.

---

*Prerna and Uday: We used i, ix and xvi (these numbers refer to the ordering of the claims on the board).*

*Me: What was your argument?*



*Uday: All sides of an equilateral triangle are equal.*

*Vandana: Angles opposite equal sides are equal.*

*Uday: So, angles are equal*

*Me: So, you don't need i (which was that the sum of the angles is 180).*

*Uday: To say the angle is 60.*

*Me: So that will be in order to find the actual value - each angle is 60 degrees.*

The students interpreted the claim 'angles opposite equal sides are equal' slightly differently to those in Ganga. While in Ganga they had assumed it meant 'given two sides which are equal, the angles opposite them are equal', these students interpreted it as 'given any number of equal sides of a triangle, the angles opposite them are equal'. I didn't press on this point. Rather, I went back to the previous proof and pointed out an assumption we had made there about the number of angles of a triangle. When showing that the angles of an equilateral triangle are equal, the proof had assumed that there are three angles and that was nowhere mentioned on the board.

The students suggested not just adding that but that triangles had three sides, three vertices, is closed and is a polygon. I then asked students why they believed that triangles had three angles. Multiple students responded by saying that is the definition. I accepted that and asked them why a triangle had three sides. I gave this to them as a group discussion task.

All the groups had discussed that two lines make an angle. So, I wrote that out on the board. I then asked for the rest of the argument.



I accepted a very non-rigorous argument for the claim: that angles require two line segments. So, if you have three angles made up of two line segments each, the only way to put them together to form a closed shape without making extra angles is by making a shape with three sides. I didn't ask them to explain why there was no other way to create such a configuration.

I then went back to equilateral triangles and asked students to prove the other direction: that a triangle with equal angles has equal sides. Udit responded:

> *Udit: Sides opposite to equal angles are equal.*
>
> *I wrote that on the board and accepted the argument. I then continued:*
>
> *Me: Now these two are equivalent. What does that mean? We can say that equilateral triangles are triangle with all sides equal. We can now also say that equilateral triangles are triangles with all angles equal. So, that is the case because 'All sides equal if and only if all angles equal, in a triangle. That means that if all sides are equal, then all angles are equal, and if all angles are equal, all sides are equal. We took all sides being equal as the definition of an equilateral triangle, right?*
>
> *Multiple students: Yes*
>
> *Me: But, this is as good of a definition. Because, if you get this, you already get this.*

In Indus, I rushed through equivalent definitions and equilateral triangles unlike in Ganga. I didn't write up proofs in as detailed a manner as I could have done, and this part of the



session ended up with me giving more expositions than I would have liked to. Given the time constraint, I could alternatively have left half the equilateral triangle discussion for the next session.

**Summary of the two set-ups**

The focus of this module is to establish relationships between statements and objects of a theory through constructing arguments which connect them and through definitions. For instance, the definitions of equilateral and isosceles triangles create a classificatory system which relate different types of triangles. The arguments which showed that the two definitions of equilateral triangles were equivalent did something similar.

The two set-up sessions on Triangle Theory Building were quite distinct in terms of content. While equilateral triangles and equivalent definitions were discussed in both the schools, the topic was the main focus in Ganga but a more rushed afterthought in Indus, and consisted of a lot more exposition from me.

Classification was not discussed in Ganga, but took up most of the time in Indus. The downside of having students suggest the initial claims is that you then have much less control over the direction of the sessions. Even though the instructor can choose amongst the given options, that option set is constrained by the students. So, if an instructor wants to touch on certain topics during a session, it makes sense to give the students a list of claims to begin with.

Another difference between the two schools is the open question I asked in Ganga. I decided not to do that in Indus. The focus of Triangle Theory Building is rigor. The value to the open question is that it brings in some creativity into the session. However, once again, it makes the session much harder to facilitate since the instructor has to be able to evaluate arguments as they are presented.



I will touch on this again at the end of the next chapter, but one interesting thing to note is that students in both the schools didn't appeal to authority or to empirical arguments in the excerpts presented. As has been shown many times in the literature (Harel & Sowder, 2007; Martinez & Pedemonte, 2014; Schoenfeld, 2013), students usually don't see value in deductive proofs. Since these excerpts represent the first session of the module, this is where you would expect this to be most visible. The reason for these students seeming to value deductive proofs could just be that many of the students have attended previous sessions with me and know what I am looking for – they did what the course facilitator expected them to do. It would be valuable to explore how these same students interact with mathematics in their regular classes in order to judge whether they do actually value deductive arguments.

The next chapter explores a few more episodes from the implementation of the module. It also discusses students feedback and ends with general observations and reflections from the course. These general observations and reflections will make reference to this chapter as well.



## Chapter 7 – Triangle theory building: Other episodes and feedback

This chapter continues from the previous chapter and is focused on the triangle theory building module. The first part is concerned with the implementation of the rest of the module after the set-up. In it, I present entire chunks of the module by which are coded by the mathematical content codes. I have picked chunks of the module which I judged to be of interest and which are somewhat representative of the rest of the module. I decided on this by focusing on parts I mentioned in my teacher notes along with those chunks which covered a large number of codes based on the coding I did of students' mathematical actions (such as conjecturing, proving, defining, etc.).

The next section of this chapter deals with students' reflections on the module. There are two different types of student reflections I have access to. The first is in the form of written feedback at the end of a particular session. In the case of triangle theory building, I have feedback on the first and third sessions in both the schools. The second is student reflections on the course as a whole which they gave on the last day of the course. My focus in this chapter is on parts of that specifically related to triangle theory building.

The final section of the chapter consists of general observations related to the research questions. As I mentioned, the main focus of my analysis is local to specific parts of the course. However, there are a few general themes worth exploring. This section will involve discussion of data from this chapter as well as from the previous one.

Here is a table of contents for this chapter:







**Other episodes**

In this section, I will be focusing on three episodes, two from Ganga and one from Indus:

1. Congruence by RHS, SAS and SSS in Ganga

2. Altitudes meeting at a point in Indus

3. Straight Lines in Ganga

**Congruence in Ganga.**

The motivation for discussing congruence came from the justification of the isosceles triangle theorem, that the angles opposite equal sides of an isosceles triangle are equal. The



students used RHS (also called Hypotenuse Leg congruence) in their justification, and also said that they could use SSS instead by dropping a median instead of an altitude.

What I wanted to achieve here, as I wrote in my notes, was to get students to see the equivalence between different ways of checking for congruence, at least where they apply. The equivalence between RHS and SSS is a good place to start since it is relatively easy to justify.

So, I asked:

*Me: How do you know that RHS gives you: triangles are congruent? What do you mean by congruence?*

*Sana: The triangles are equal*

*Gauri: They coincide*

*Sana: When kept one on top of*

*Imran: Overlap*

*Vivaan: If we rotate or move or reflect the triangle, then they coincide. Using only these three operations.*

*Me: So what you are saying is that: SSS shows congruence and you are telling me RHS shows congruence. So, in some sense, at least when we are restricted to right angled triangles, there is this relationship: RHS is true if and only if SSS is true. If the right angle, side and hypotenuse are true, the sides are true.*

*Sana: No, only two sides.*

*Me: No, you're telling me that if this was true, if you had RHS. The right angle, hypotenuse and side were true, you also get the sides is equal to the other one.*



> *And therefore you get that if this is the case, then this is the case, and if this is*
>
> *the case, then this is the case. Is that what you are saying?*
>
> *Sana: Yes.*

In my notes, I mention that I could have spent more time on the conceptual aspect of congruence and connected it to the more formal aspects better. The students were thinking of congruence as placing a shape over another and seeing whether they overlap. I instantly took it to the equivalence of different ways of checking for congruence in triangles. Instead of that, I could have discussed the consequences of the overlap notion of congruence for the different properties of triangles.

I was also not very careful about using the word 'true' in this exchange. When I said 'RHS is true', I meant given two triangles, there is a correspondence between the right angle, hypotenuse and side in those two triangles. I could have made this more explicit. However, the more egregious use of the word true was when I said, '*If the right angle, side and hypotenuse are true, the sides are true*'. Here, the word 'equal' would have been much more appropriate. This is a good example of me not modeling the clarity and precision in language I was expecting from the students.

I continued:

> *Me: In one case, that seems very straight forward, right? So, if you have,*
>
> *assuming we are restricting ourselves to a right angled triangle…*
>
> *Vivaan: We can prove RHS without SSS also.*
>
> *Me: What do you mean by prove RHS?*



*Vivaan: We can show RHS always shows congruence.*

*Me: But, what do you mean by congruence? Overlapping?*

*Vivaan: No*

*Me: Then, what do you mean by congruence?*

*Vivaan: That all the dimensions are equal.*

As I said in my notes, my intention with this session was to get students to justify that the methods of testing for congruence are equivalent to each other. One way of doing that is to demonstrate their equivalence by showing if and only if relations between them. The other way, which Vivaan seems to be suggesting is to show that all of them lead to corresponding aspects of triangles being equal. This would show that all of them give us congruence. However, it would not show equivalence between these different ways of getting at congruence in their domain of applicability. As an analogy to show why, we can think about similarity of triangles. Here are three possible ways in which triangles can be similar:

1. All the corresponding angles are equal

2. Two of the corresponding angles are equal

3. All of the sides are equal

The first two are equivalent since, if two corresponding angles are equal, the angle sum property gives us that the third is also. The other direction is obvious. However, while the third does give similarity, and implies the first and second, the third is not equivalent to them since angles being equal does not give us that sides are equal.



I did not make this intention to show equivalence completely explicit during the session and I did not adequately address why I decided not to follow the direction Vivaan's seemed to be suggesting. Rather, I said:

*Me: Yeah. So, that's what we are trying to do, right? So, because now you are saying that if we have SSS, all sides are equal, in a right angled triangle you get right angle. You get RHS, right?*

*Anya: Yes.*

*Me: That's very clear, right? Say we have SSS on a right angled triangle, and we know there's a right angle. So, of course, if we know SSS, we get RHS because ignoring that this is equal to this (one of the sides). Everyone see that?*

*Multiple students: Yes*

*Me: So, SSS implies RHS. Okay. So, why does RHS imply SSS? So, spend a minute thinking about that.*

In this exchange, I justified that SSS in a right angled triangle gives us RHS. It was an easy justification and it seems like I was just trying to get it out of the way and focus on the slightly harder direction. However, it could have been useful to ask students to justify the claim. Given that it is relatively straightforward, it would be something which most of the students would have been able to do.

After they came back from the group task, the discussion was as follows:



*Imran: So, by Pythagoras Theorem, we could get the third side as well. Any then by, and we know the equation is going to work for both the triangles. If two of the sides are the same, then the third side should also be the same.*

*Me: Okay, so if this is c, this is a, then this is square root of $c^2 - a^2$ and this is square root of $c^2 - a^2$. So, there you go.*

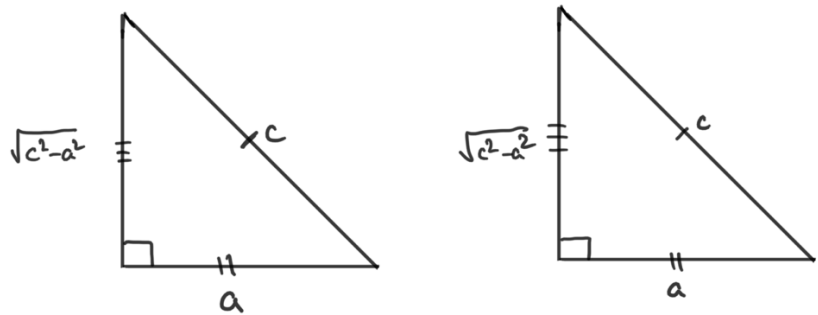

*Figure 27 Congruence of right angled triangles*

*Anya: But if we have just the triangles and we don't know the measures?*

*Me: You know that this is equal to this, right? We don't know the measures. We know that this is equal to this. Just call it c.*

It seems like Imran understood the goal of the exercise. However, given what Anya said, it is unclear whether she did. In the group discussion, she did not say much apart from agreeing with Vivaan's argument which used Pythagoras' theorem. One possible interpretation of what she said is that she thought of the question as asking, 'how do we know that two given triangles are congruent?' rather than, 'given two triangles whose are the same by RHS, how do we know that they are the same by SSS?'. It might be useful to make this distinction explicit in future



sessions of this type. That ended the discussion on RHS and SSS. I then gave them a group task to figure out why SSS and SAS are equivalent.

### *Summary of Congruence.*

One of the reasons this short excerpt is useful is that there seem to be many ways in which I could have done better. For instance, I could have given students some relatively easy claims to justify on their own. I could also have been clearer with my instructions as well as more precise in my language. Finally, I could have probed more into what Anya said to understand where she was coming from rather than responding as I did.

The intuitive concept of congruence is to do with being able to put one object on top of the other. Now the question which arises is how we can formalize that idea in the case of triangles. If we know that the corresponding sides and angles of two triangles are equal, that results in congruence.

Students have learned many rules on how to test for congruence. The idea behind these rules is that they predict congruence. If we know that any of these rules hold, then the other aspects of congruence also hold. In fact, all the other rules of congruence also hold.

One way of justifying the concept is to show the equivalence of these rules. If we can show that SSS iff SAS iff ASA, and so on for the same triangles, then we have justified the concept. That is what we were trying to do in this session. I started with a pair of rules which are relatively easy to justify from each other – RHS and SSS.

### Altitudes meeting at a point in Indus.

I was initially intending to discuss area of triangles and the sum of angles of a triangle. As I wrote in my notes, I picked the meeting of altitudes as a claim to discuss since students,



specifically Tarini and Tanya, had expressed interest in discussing it after the previous session. I wrote down the claim on the board and said:

> *Me: So, that's kind of surprising, right? Because, it wouldn't be surprising if that two, let this be an altitude and this be an altitude, that two altitudes meet at a point. That's not that surprising if that were the case. The surprising thing is that the third passes through the same point. So, why?*
>
> *Meghna: That's a theorem.*
>
> *Me: Why is that a theorem? Did you understand the question? So, altitudes of a triangle are these perpendiculars we drop from the top to the bottom, right? That two of them meet at a point is not surprising. Two lines meet at a point is not surprising. But, that the third passes through the same point is what is surprising. So, why do they all meet at the same point?*

Unlike the case of Congruence in Ganga, I was much clearer here in setting up what we were trying to achieve, and what the claim was. In other words, I clearly articulated what we should be surprised about and hence need justification for. One part of this which could have been made clearer what the phrasing of the question. Rather than just saying 'why do they all meet at a point?', it might have helped to have added that I was looking for justification for the claim that they do.

One assumption I made above is that all students believed the claim. Given that some of the students had mentioned the claim to me, I seem to have assumed that all the students would have accepted the truth of the claim. However, that was not the case:

> *Uday: They were not altitudes. They were perpendicular bisectors.*



*Me: No, this is not a perpendicular bisector.*

*Uday: No. You have drawn it in a way that they are perpendicular bisectors.*

*Me: So you think altitudes don't meet at a point?*

*Uday: Not altitudes.*

*Me: But, everyone else is saying that altitudes always meet at a point. Do you think altitudes always meet at a point? (to the whole class)*

*Multiple students: Yes*

*Uday: Sorry I got mixed up.*

*Me: Okay. So, why do you believe that altitudes meet at a point?*

There are two interesting aspects to this exchange. The first was Uday's response where he said, 'you have drawn it in a way that they are perpendicular bisectors.' He seems to be saying that perpendicular bisectors necessarily meet at a point. Altitudes do not necessarily meet at a point. However, when these altitudes are also perpendicular bisectors, they do. The accusation he seemed to be making was that the picture I had drawn was misleading. Whether it applied to this situation or not, this insight that pictures of examples can mislead us about the entire class of objects is a valuable one. I could have picked up on this and made it explicit to the rest of the class.

The other thing of note in the above exchange is the way I dealt with Uday's concern. A potential path I could have taken would have been to get Uday to try drawing some examples on a piece of paper and convincing himself of the claim through that process. Rather, what I did was



to ask the rest of the class whether they believed the claim and potentially put pressure on Uday to accept it, whether or not he actually did.

---

*Uday's objection was not the last objection to the truth of the claim.*

*Devika: If there is an obtuse angle triangle, then the altitudes won't meet.*

*Tarini: They will.*

*Devika: Altitudes is right angle.*

*Me: That's a good point. This altitude is easy (drawing an altitude on an obtuse angle triangle from the obtuse angle). Where is this altitude? (pointing at another vertex)*

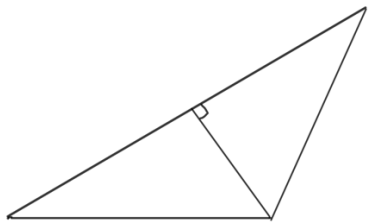

*Figure 28 Obtuse angled triangle*

*Devika: It will be outside.*

*Me: So, you are saying you have to extend it, that's your altitude.*

*Devika: But they are not meeting. Or are they?*

---



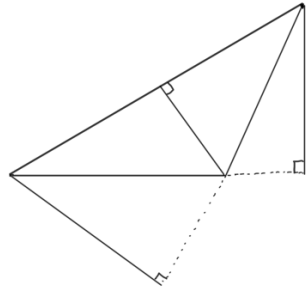

*Figure 29 Obtuse angled triangle with altitudes*

*Tarini: So, if you extend it, then it will meet.*

*Me: So, if you extend them, then they meet at a point.*

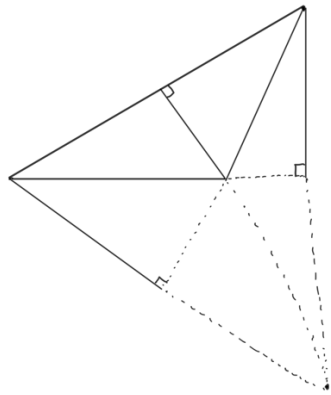

*Figure 30 Obtuse angled triangle with altitudes meeting*

Devika's objection was similar to Uday's in the sense that if she turned out to be right, then my initial figure of an acute angled triangle was misleading. However, it isn't clear to me from the exchange whether she believed that altitudes did not meet for obtuse angled triangles or whether she was questioning the very existence of altitudes in obtuse angled triangles.

By clearly working through an example and not polling the class, I dealt with this situation quite differently to how I dealt with Uday's objection.



I gave the task of justifying the claim as a group task. Once they returned, Varun came to the board and drew the following.

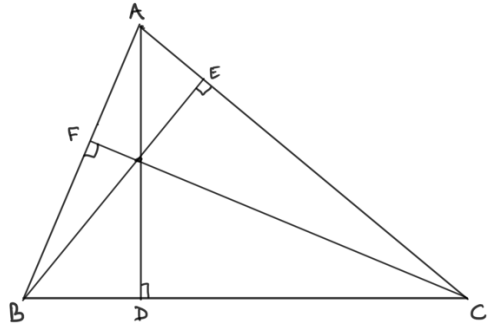

*Figure 31 Altitudes meeting*

*Varun: First, you join these (E and F)*

*Me: No. How, did you create this, first.*

*Tarini: We just created it*

*Me: No. It's not exactly that*

*Varun: Okay. So, If I consider this as O (the point where they meet). C, O and F are not colinear.*

*Me: Not necessarily, right? What he's done is he's drawn two altitudes. The point where they intersect, he called that O. Then he connected CO and dropped a perpendicular from O to F.*

*Tarini: And they are not colinear.*

*Me: We don't know. So, that's what he wants to show. That they are colinear.*



*Varun: So, we join EF and we get a cyclic quadrilateral here.*

*Me: Okay. So, why a cyclic quadrilateral?*

*Varun: So, if this is x, then this is also x.*

*Me: The angle?*

*Varun: Yes the angle.*

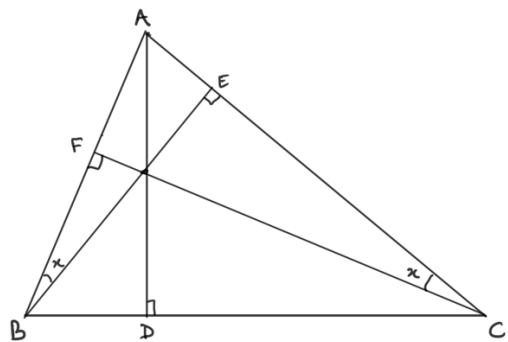

*Figure 32 Altitudes meeting - equal angles*

Varun's proof used the concept of a cyclic quadrilateral. As I wrote in my notes, the proof I had in my mind involved an elaborate construction using parallelograms and perpendicular bisectors. My intention was to share this proof with students if they could not come up with a proof and then work through the premises required for that proof. However, Varun's proof is not something I was familiar with. Even since then, I have not seen a proof using cyclic quadrilaterals in the way he set it up. The proofs I have seen involving cyclic quadrilaterals use that the quadrilateral formed by joining D and E in the figure above is cyclic.

Once Uday assumed that the quadrilateral was cyclic, he was easily able to show that the COF was a straight line. I wrote down some of the things he used in the proof:



1. YOU CAN DRAW A CIRCLE SUCH THAT
   ALL VERTICES OF A CYCLIC QUAD.
   ARE ON THE CIRCLE
2. ANGLES SUBTENDED BY THE SAME
   CHORD ON ITS CIRCLE ARE EQUAL
3. SUM OF ANGLES a) 180 (△)
                   b) 360 (□)
4. LINEAR PAIR IMPLIES STRAIGHT LINE
5. CYCLIC QUAD. IS ONE WHERE SUM
   OF OPP. ANGLES IS 180°

*Figure 33 Cyclic quadrilateral claims*

I then asked Uday:

*Me: How did you come to the conclusion that this is a cyclic quadrilateral?*

*Uday: Because we have proved that when two angles are subtended.*

*Me: You know that if it was a cyclic quadrilateral, then those angles are equal, right? But, how do you know it is a cyclic quadrilateral?*

*Uday: We are using the same thing.*

*Me: No. No. You used the cyclic quadrilateral to show me that this equals this, right?*

*Uday: Yes.*

*Me: Now you can't use this to show me that it is a cyclic quadrilateral.*

It doesn't seem like Uday had a justification for why the quadrilateral is cyclic. He had earlier used that the quadrilateral was cyclic in order to justify that two angles were equal and was now using the two angles being equal in order to justify that the quadrilateral is cyclic.



The class spent a little more time on trying to justify that the quadrilateral was cyclic till I decided to move them on to discussing the other premises that Uday had used. However, I did come back to discussing altitudes later in the class and gave the students the proof I had initially intended to share.

### *Summary of altitudes meeting at a point.*

I said in my teacher notes that this episode did not help all that much in achieving the goals of the module. The goal of the module is not that students learn Euclidean Geometry. Rather, it is that students learn to build theories. Hence, concentrating on claims which I am quite certain students would be able to prove would have probably been a better choice. Alternatively, I could have given students the initial proof right at the outset and then asked them to dig into the premises required for that proof.

What ended up happening is that we got stuck on the initial proof of the claim. It did not result in interesting discussion and we did not really touch upon any of the important learning outcomes associated with theory building.

### Straight lines in Ganga.

We had used the term straight line many times in the module till that point without defining it. So, I asked:

*Me: What is a straight line?*

*Gauri: A set of non-collinear points.*

*Me: A straight line is a…*

*Gauri: Set of collinear points.*



> *Me: Set of collinear points (while writing that out). What does that mean?*
>
> *(pointing at the word collinear)*
>
> *Vivaan: Collinear lies on a line. So, we will use straight line in the definition.*

Vivaan here seems to see the problem with defining a straight line as a set of collinear points – that we will use the notion of straight line to define collinearity. Given the ubiquity of this definition in this course in both the schools, it seems like it is something they are taught. However, this is not how straight lines are addressed in the textbook they use so if they are learning it in their mathematics classes, it must be something their teachers are saying.

> *Me: Points on a straight line (while writing that out). So, notice the problem with this. You define straight line as collinear and collinear as straight line. You see that being a problem, right? Circular. So, this is an important thing. Remember I drew those diagrams? This statement about opposite angles being equal. Let's call that statement A. That relies on many things but amongst them it relies on statement F. We used statement F in the proof. Here we needed many things. But here we used the definition of triangle and definition of perpendicular, and for both of these, we used straight lines, and for this we used collinear for which we used straight lines.*
>
> *So, this is where problems arise in this sort of diagram. It is okay if you use straight line in things which are in different parts of this graph, in different branches. You can use that as many times. But, in the same branch, you have a branch where you use straight lines definition and you use it again, that's a*



*problem. Because that means there is circularity. Basically, it means you are*

*doing this. (I draw an arrow upward)*

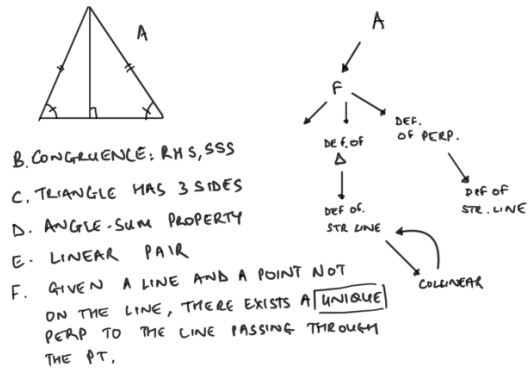

*Figure 34 Argument in diagram form*

The is potentially an issue here with how I talked about the diagram. When I said that 'this is where problems arise in this sort of diagram,' I did not mean problems in the diagram but problems in the structure represented by the diagram. This distinction between a representation and the entity it is representing is something I want students to appreciate. While making that distinction every time is not possible, and referring to the representation as shorthand for the entity is inevitable, it would be useful to occasionally make the distinction explicitly.

Imran then suggested an alternative definition:

*Imran: Can we use shortest path?*

*Me: So, you want to use shortest path.*

*Imran: So, all the points between them are collinear points.*

*Me: So, you are using shortest path. Now, what does this rely on?*

*Vivaan: What is a path?*



*Me: One is 'what is a path?' What else does it rely on?*

*Vivaan: What is shortest?*

*Me: So, what is length? What do you mean by length? What is the distance between two points?*

*Imran: Maybe the points between them.*

*Me: So, when we were doing discrete geometry, it was very easy to figure out shortest path. Right, because there's only these many points. You just count the number of things.*

*Vivaan: In Euclidean Geometry, then all lines will have equal length. See, there is infinite number of points between two and there is also infinite between these two.*

*Me: Yeah. So, shortest path requires distance. And distance, what you were trying to say (looking towards Imran), was length of the straight line. If you pick two points, then distance is the length of the straight line.*

The problem with defining straight lines as shortest paths in Euclidean Geometry is what I said at the end – that distance is defined as the length of a straight line segment. Unlike in the discrete geometry module, there is no obvious metric we can use and, unlike in coordinate geometry, we cannot impose an externally defined metric. Taking the discrete geometry strategy of using the number of points has the problem which Vivaan mentioned – that lines in Euclidean Geometry, no matter how small they are, don't have a finite number of points.



However, it is unclear as to whether Imran was 'trying to say' that distance between two points is the length of the straight line between them. I seemed to be putting words in Imran's mouth when I said that. An alternative strategy there would have been to elicit that from Imran by pushing him, or others in the class, to define distance and then pointing out the circularity in the definition like I did above with collinearity.

After this, Gauri suggested that we define straight line as collinear but we do not define collinearity:

> *Gauri: What if we consider collinear points as the axiom. We don't have to further define it.*
>
> *Me: So, you don't want to further define collinear. Sure, you can do that. But, you need something to tell me what collinear points are, right? Because, why are these points not collinear? (I draw three points on the board) So, tell me why are these points not collinear?*
>
> *Gauri: Because they are not joined together by a straight line.*
>
> *Me: That's one way to think about it, but can you think about it without using straight lines?*

As I said in my notes, I was avoiding using the term axiom as far as possible during the session even though I did end up using it soon after. The reason for this is that I wanted students to understand the concept before putting a label on it. Gauri used the term axiom here to seemingly mean undefined entity since she wants it to be something you do not define further. Rather than explicitly saying that the term axiom is not synonymous with undefined entity in mathematics, I decided to rephrase what she said without using the term.



This ended the interactive part of this session. For the rest of the time I spent on this topic, I did two things. I first gave the analogy of defining a t-shirt. I said that whatever you give as your definition, I can keep asking you what you mean by the words you have used. If you are dealing with real world objects, at some point, you will have to point at the object. However, in mathematics, you cannot do that since the lines we draw on the board are just representations of lines rather than lines. So, we have to treat certain objects as undefined and put some rules on them called axioms/postulates.

Then, I gave the example of the axiom, 'given two points, there is a unique straight line which they lie on.' This doesn't tell us exactly what a straight line is, but it does constrain the space of things the term straight line could be referring to. This was the last few minutes of the last session on triangle theory building. As I said in my notes, I decided against discussing other axioms of Euclidean Geometry including the parallel postulate due to time constraints.

### *Summary of straight lines.*

This was the first time in the course that we arrived at undefined entities and axioms. This came right at the end of the course and hence we did not spend a lot of time on it.

One of the mathematical insights was that you cannot keep defining objects in terms of other objects. Eventually, that will result in circularity as happened when students defined straight lines as collinear points and collinear points as points on a straight line. However, we cannot just leave the object as undefined. We need to constrain its interactions with other undefined objects by using axioms.

We have to decide what object to leave as undefined. In the session, we chose straight line. We didn't get very far in discussing the axioms - I didn't want to introduce the parallel postulate since that would have taken a significant amount of time to work through.



**Student reflections on the module**

Students wrote reflections at the end of the first and last days of the module. They were asked to answer the following two questions:

1. What did you learn?

2. What did you not understand?

Apart from this, students were asked for general reflections on the entire course at the end. In the reflections, they were asked the following four questions (the template for the survey is available in the appendices):

1. What was your favorite part of the course? Why?

2. What was your least favorite part of the course? Why?

3. What did you find valuable in this course? Why?

4. What are one to three specific things about the course which can be improved?

The goal of this subsection is to give a sense of what they said related to this module. As mentioned in the methods section, I coded sentences in students' daily reflection responses on the basis of whether the responses were related to specific mathematical content or to the types of general theory building abilities I discuss in Chapter 3. Some sentences could have both codes. The sentences without codes turned out to be things like working in a team or to do with the instruction. So, I added in two additional codes which were 'non-mathematical learning' and 'reflections on instruction.'

The same set of codes were applied to the end of course reflections along with an additional set of codes to do with whether a particular sentence applied to a particular module. Here, I will be mainly concerned with those instances related to triangle theory building.



**End of Course Survey.**

Starting with the end of course survey, there seems to have been a significant difference between the sentiments between the two schools. Almost all of the students interpreted the first two questions in the survey, to do with their favorite and least favorite parts of the course, to be referring to the three modules. The responses to the last two questions barely mentioned the modules and so I will not be discussing them here.

In Indus, 12 of the 16 students specifically mentioned triangle theory building amongst their favorite parts of the course and no student mentioned it amongst their least favorite. In Ganga, on the other hand, only five out of 15 students mentioned triangle theory building amongst their favorite parts of the course while five mentioned it amongst their least favorite parts of the course. To try to understand what is going on here, it is useful to see what the students gave as their reasons.

The negative reactions to this session included:

- '…it was a repetition and too easy for my grade.' – Gauri, Ganga

- '…we kept going endlessly…' – Anuj, Ganga

- '…(it) gave minimal opportunity for finding alternatives and doing something different/new.' – Anjali, Ganga

On the other hand, the positive reactions included:

- '… in a classroom we only use these theorems blindly…' – Vivek, Ganga

- '…finding proofs, proving them wrong, again finding another proof.' – Anya, Ganga

- '…I liked the way we figured out the interdependency between properties.' – Tushar, Indus



- '…I had some knowledge on the triangles and I can answer to the questions asked.' – Aryan, Indus

- 'During the triangle session we were asked to define everything… and I enjoyed that.' – Meghna, Indus

- '…it intrigued us to prove things we normally would not question.' – Uday, Indus

It is interesting that some of the reasons students gave for the module being one of their favorites is very similar to the reason others gave for the module being amongst their least favorites. For instance, 'finding proofs, proving them wrong, again finding another proof' is similar to 'we kept going endlessly.' They appear to be opposite reactions to the same stimulus.

However, it isn't clear whether the difference between the two schools is best explained by the students, the way I conducted the sessions, the time of day, or even a single influential student who swayed the rest.

Two other opposite reactions worth comparing are 'it intrigued us to prove things we normally would not question' and '(it was) too easy for my grade.' It seems like Gauri was talking about conclusions and saying that she already knew the conclusions, while Uday was talking about the justification process.

However, the two reactions I want to highlight are Anjali's and Aryan's. Aryan mentioned during his interview that he was lost during large parts of the discrete geometry module which was completely new to him while he was able to contribute to this module since he had some knowledge about triangles. On the other hand, Anjali says she did not enjoy the module as much since it did not give an opportunity for her to do something new.



**End of Day Reflections.**

Most students did mention specific understanding they gained about the mathematical content – about triangles, congruence, and so on. I will not be focusing on that here as that is not a central research question. My focus here will be on the more general statements students made about what they learnt from the two sessions in this module for which they wrote reflections.

On the first day in Ganga, every students said that they learned about the relationships between statements. For instance:

- Vivek said, 'We learned relationships between statements or definitions of properties of triangle. It is that there is an interrelationship between them. So, if we have 4 statements, 1st statement can be proven by other 3 of them. 2nd statement can be proven by the 4th and other permutations and combinations.'

- Deepa said, 'I got to know how each property of triangle are interconnected.'

- Aditya said, 'I learnt or improved to build and find relationships, use them to prove something else, this cycle continues to go on'

In Indus, nine out of 16 students mentioned relationships between statements as part of their learning. For instance:

- Aditi said, 'Today was an interesting class on theory building of triangles. It made me realize how all properties are inter-related. Proving various theorems was fun and engaging.'

- Vandana said, 'Today we/I learnt how different theorems about triangles are inter-related. They are one and all the same and how different congruencies are related to each other.'



- Serena said, 'Found new techniques of proving. And also how to prove from one point to other point in different ways.'

In Indus, additionally, four students mentioned classification amongst what they had learned.

On the last day in Ganga I only got 8 reflections. This is because some of the students had to leave for some school event. Out of those 8, four mentioned definitions. For instance:

- Anya said, '(we learned about) Definitions & Axioms & difference between them.'

- Sana said, 'One thing can be defined with other things that we know, Not all things can be defined'

In Indus, one of the main themes mentioned was not blindly accepting what textbooks say and not taking things for granted – eight out of 15 students mentioned these. For instance:

- Tarini and Tanya (who wrote exactly the same thing) said, 'Not to take anything given in the textbooks for granted'

- Vandana said, 'I started looking at mathematics from different perspective. Now I do not restrict myself to textbook related theorems.'

- Arun said, 'today I learnt to think and find about the reasons behind certain things which we take for granted in our day to day world.'

Another theme was thinking in different ways – ten students mentioned this. For instance:

- Vandana said, 'It helped me to think about things from out of the box.'

- Udit said, 'I learnt different ways to prove the proofs and also think in different ways.'

- Sanya said, 'I learnt how to work in teams, discuss & explore ourself to mathematics in different way which we had never thought of. The knowledge which we got from



this particular session was something we knew but it made our mind think a different

perspective.'

Many of the students in their reflections used words which are closely aligned to the

learning outcomes of the module and the course. This does not mean that all the students actually

achieved the learning outcomes to any extent. However, it does at least points to them being

sensitized towards the learning outcomes.

Another interesting thing to note here across these reflections is that many of the students

who say they learnt things which are important learning outcomes of the module did not speak

up much during class discussions. For instance, in the episodes I chose to focus on, I did not find

a single instance where Sanya said something during the class discussion. While she did engage

in the group discussions, it is worth thinking about ways in which to encourage such students to

participate in class discussions – for their own learning but more importantly for collective

learning. This does not have to be oral participation – it could involve written work which then

gets used during the class discussion.

**General observations and reflections**

   **Definitions and classification.**

One of the questions of theory building addressed in this module was 'what should be

part of our inventory?'. For instance, the discussion on '47.5 degree triangles' in Indus. The

argument for not including 47.5 degree angles in Indus was that there are not theorems about

them and that the classification becomes complicated. That isn't strictly true. I did not bring up

that you could rewrite the Pythagoras Theorem using the law of cosines to be about 47.5 degree

triangles. The reason for that is they would not have seen the law of cosines before. Had they

seen the law of cosines before, a better argument would have been that the expression would be



unnecessarily complicated. Another important consideration is that Euclidean Geometry is a model for compass and straight-edge constructions. Making arbitrary angles is not possible using this technique.

Just before this discussion in Indus was a discussion about classification of triangles by the equality of their sides. The first impulse of a number of students was to separate equilateral and isosceles triangles. Once they were convinced that equilateral triangles should be types of isosceles triangles, a few of the students generalized the same type of reasoning to say that isosceles triangles are types of scalene triangles. Since doing this results in scalene triangle being synonymous with triangle, this is not useful. One possibility move here is to keep scalene and isosceles triangles as distinct. Another possibility is to drop the concept of scalene triangles from our inventory.

There was an opportunity to discuss choosing between definitions at the beginning of the module in Ganga which I missed. Anya justified that the right angle is the largest angle of a right triangle by appealing to the assumption that a triangle cannot have a zero degree angle. Suppose we allowed for degenerate triangles to be triangles where the degenerate case here is when one angle is 180 degrees and the other two are zero degrees. The claim that the right angle is the largest angle of a right angled triangle would still be true. However, the argument would have to be different since we cannot rule out the possibility of zero degree angles in a triangle. This was an opportunity to see a relationship between definitions and arguments for a conclusion.

**Flawed reasoning and proof.**

There were a few instances of what looks like flawed reasoning. This seemed to be especially evident when I left students with 'the open question' in Ganga. Anya's initial proof of the angle subtended by the diagonal being a right angle involved using two different triangles



and the claim that angles opposite greater sides are greater. In this case, it isn't completely clear as to whether the reasoning was bad or her and her group's understanding of the assumed claim was lacking. If the underlying claim was true for different triangles, then the reasoning was sound. Hence, it may be the case that the reasoning was valid.

However, the case of the third group to present on the open question is more interesting. They claim they were proving was the converse of the claim the second group had presented. However, they said that it was the same claim. This is more clearly a matter of bad reasoning where the students seemingly concluded P => Q from an argument for Q => P.

An advantage of asking an open question as opposed to asking students to justify a certain claim is that with the open question, there are many more ways for students to make mistakes and hence a chance to engage in interesting discussion.

Another two instances of bad reasoning worth bringing up are to do with cyclic reasoning. One example was of straight lines in Ganga where students defined straight lines as collinear points and collinear points as points on a straight line. This was something which was immediately noticed by the students. It is also inevitable when reaching the foundations of a theory given that students do not have a deep understanding of axioms and undefined entities.

The more intricate example of cyclic argumentation was when Uday used the notion of cyclic quadrilateral to justify altitudes meeting at a point. He seemed to use the existence of a cyclic quadrilateral to justify that the sum of opposite angles was 180 degrees. He then used the sum of opposite angles being 180 degrees to show that the quadrilateral was cyclic.

**Proof schemes.**

What I did not see in these students was any significant attempt to give proofs based on an authoritative or empirical proof scheme. In fact, there was an instance when Devika suggested



that the claim I was making for altitudes meeting was not valid for all triangles, only for acute ones. She was saying that the claim I was making was false, at least for a large class of the objects I was making the claim about. While that may not have been true, the valuable part is that she was looking for counter-examples and cases where the claim may not even make sense.

Given that the literature in mathematics education is full of examples where students use empirical and authoritative proof schemes (Schoenfeld, 2013; Harel & Sowder, 2007; Martinez & Pedemonte, 2014), this is surprising. There are a few possible reasons for this. One possible reason is selection bias, given that this course was something students chose to be a part of. Had they not seen value in either being a part of previous courses I taught or seen value in what their friends or teachers told them about the course, it is unlikely that they would have chosen to be part of it. There is further room for selection bias in those who chose to speak during class discussions.

Another possible reason is that many of these students have been through previous sessions I have done and sessions by others including Mohanan, who was present for the first session in Ganga. Most of the students in Indus were present for a workshop I did in 2015 which I discussed in my Master's Thesis. In Ganga, some of them were present for that workshop, but many of the others were part of a course I was involved in designing when they were in the sixth grade. It is hard to say what the impact of these interventions was, but given that students chose to be a part of this course, they seem to have seen value in those previous interactions. It is plausible that they saw value in the deductive proof scheme through those interactions.

A third possible reason is that the Eucldiean Geometry course they did was extremely good. This seems unlikely since other things they said indicate otherwise. For instance, defining straight lines as collinear points and collinear points as straight lines. Also, some of the students



in Indus thought of equilateral triangles as distinct from isosceles triangles. It is unlikely that a good Euclidean Geometry course would have taught them that. Also, if they had a good Euclidean Geometry course, then erecting the structure I wanted them to about triangles should have been very straightforward.

**Productive struggle.**

There were a few instances in this chapter where I mentioned that I could have gotten students to think more on their own rather than helping them. There is always a hard balance to strike between giving students the chance to come up with things on their own and helping them when they really need it. The goal is for students to engage in productive struggle (Heibert & Grouws, 2007). This balance is heavily context dependent, based on the students' background and even their particular mental state on that day. It is hard to say much generally on what an instructor should look out for in their decision to intervene.

However, it is important that the struggle which students endure is targeted towards the goals of the lesson/course. An example of a situation in the implementation of this module when that was not the case was during the session on altitudes in Indus. A large part of the session was spent on coming up with and understanding a particular proof for the claim. It was almost as if the goals of the course shifted from theory building to understanding Euclidean Geometry. Two ways to have avoided this in that situation would have been to have chosen something simpler to inquire into or I could have presented a proof to the claim myself. In either of those cases, the focus could have shifted from the particular proof to the relationship between statements of the theory, a return to the goals of the session.

In this instance, students struggled with something which did not directly contribute to the goals and such situations ought to be avoided.



**Interpreting student responses.**

There were a number of times that I interpreted what students said when there appeared to be ambiguity. I said things like, 'you mean that…' or 'what you mean is…'. I would have much less concern doing this if this was a group of students I had been working with for months on theory building as long as I gave them the opportunity to push back against my interpretation. However, in this situation, two possible downsides of doing this are:

1. Students may accept my interpretation even though they meant something different due to that interpretation coming from an authority figure.

2. Students at this stage may not have the ability to clearly articulate their interpretation and hence could struggle to push back.

This does not mean that the instructor never interprets students' responses. Asking students to clarify everything they say without help would require a huge amount of time. They also may not be able to. A compromise position would be to ask them initially to clarify what they said. If they are unable to, then, if possible, offer multiple possible interpretations of their statements for them to choose from. However, always give them the opportunity to suggest an alternative in case they did not mean any of the options given.

**Following students' direction.**

This module is the most open ended of the three in terms of the directions it could take. There are a large number of claims an instructor could choose to explore. This needs to be done in negotiation with the students.

There is a compromise to be made between achieving the intended learning outcomes of the module efficiently and following the direction students take the class in. It is useful here to make a distinction between the beginning of the module and later on in the module. At the



beginning of a module of this type, it is important to cover a few different types of claims so that students see a range of different ways in which they have to think. In the implementation of this module, I did not discuss classification in Ganga since students did not bring it up. Since classification is an important intended learning outcome of this course, it is unfortunate that students in Ganga did not get to spend time on it.

So, while listing out the initial set of claims, it is important that they would cover a large number of learning outcomes if taken up for inquiry. This can be done in a few ways. One would be for the instructor to give the initial set of claims. A second possibility is to add in some additional claims when the students give theirs and a third option would be to get each group to write out a large number of claims and then pick out a few from the written set of claims to discuss.

While it is important to maintain some flexibility in the beginning in terms of the particular claim being inquired into, it is crucial for the instructor to be able to relate claims to intended learning outcomes as they go about the session.

As the module continues, there can be more and more flexibility on what claims to pursue and even what learning outcomes to aim at based on students' interests. However, even in this case, it may be best to avoid situations like the case with the altitudes in Indus where the claim picked was relatively hard to prove. Also, the more flexibility given to the direction of the module, the more comfortable the instructor has to be with Euclidean Geometry. Especially when a teacher is teaching the course for the first time, the best strategy is likely to be them giving the students the initial list of claims – claims they are comfortable with. After a few iterations, they might be able to take a more open-ended approach.



## Chapter 8 – Discrete geometry: Changing worlds

This and the next chapter are focused on the third course module – discrete geometry. A sketch of a lesson plan for this module is available in Appendix A. Like with the previous results chapters, this chapter describes the implementation of the module in the two schools. The intention of doing this is to explore how students engaged with the course and the role played by the facilitator.

After describing the mathematical content and intended learning outcomes of the modules, the rest of the chapter deals with the set-up of the module. In this case, by the set-up, I mean the first question in the module in the two schools. I present the entire set-up in both the schools in a sequential manner. The next chapter will include a description of other episodes in the course as well as students feedback and general reflections and observations on the implementation of the module.

Here is a table of contents for this chapter:







## Description of the module

This module starts with a question: In a world with exactly six points, can every straight line be bisected? This is an ill-formed question since we don't know what a world with six points looks like, what a straight line in that world is or even what we would mean by bisection in such a world. The idea in this module is to borrow ideas from what students are familiar, namely Euclidean Geometry, and make them meaningful in the context of the question.

Students have to first interpret the question and then answer it. Different interpretations will result in different conclusions. Once they have created some worlds which are meaningful, the next part of the module is to explore those worlds and explore concepts like circles and triangles in those worlds. The learning outcomes associated with this module include an understanding of:

1. definitions

2. the relationship between assumptions and conclusions

3. conjecturing, and generalizing claims and theories

## The set-up

The set-up of the Discrete Geometry session is much easier to demarcate clearly that in the case of triangle theory building. I'm using set-up here to refer to the part of the Discrete



Geometry module where students were answering the initial question I posed. I will be dividing this section into four parts:

1. Posing the Question

2. Clarifying the Concept of Straight Line

3. Transferring the Concept of Straight Line to the World Implied by the Question

4. Engaging with the Concept of Bisection and Answering the Question

In each of these sections, I will look at how these activities occurred in the two schools, with a special focus on how they were different and similar to each other.

**Posing the question.**

The goal of the initial interaction was to clarify the question – give meanings to the words so that we can answer the question can be addressed. In initial posing of the question, there are at least two directions which can be taken. One possibility is that the instructor poses the question without suggesting a strategy initially. Only when students have struggled sufficiently does the instructor point out potential entry points. The other possibility is that the instructor starts with outlining the strategy to be taken. When this session was done in Ganga, I took the prior strategy while when the session was done in Indus, I was closer to the latter strategy. Here are the excerpts from the video, first Ganga and then Indus.

Ganga Excerpt:

*Me: In a world with exactly 6 points, can every straight line be bisected? So, in Euclidean Geometry, what is your answer to this question? Not a world with 6 points, but in Euclidean geometry, can every straight line be bisected?*

*Multiple students: Yes*



*Me: Yes. So, what do we mean by bisected in Euclidean geometry?*

*Sana: Divided into two equal parts*

*Me: Two equal parts.. in some sense of the word equal. What about in a world with 6 points?*

*(No response)*

*Me: Do you know what a world with 6 points looks like? Exactly 6 points. No more, no less. What does it look like? What would such a world be?*

*Pankaj: Because in every two points, there will be a line segment. And that can be bisected.*

*Me: There are six points, right? So, what do you mean that between any two points there will be a line segment*

*Pankaj: From one point to another*

*Me: But the line must contain some points? This has how many points? A line. A Euclidean line has how many points?*

*Multiple students: two*

*Me: I mean how many points does it have on the line?*

*Multiple students: infinite.*

*Me: However long or short you make it, it will have infinite points. So, this (the diagram I drew) is more than 6 points.*



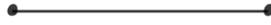

*Figure 35 Euclidean line*

*So, what would a world with 6 points be? Exactly 6 points. Any thoughts?*

*I asked them to discuss this in their groups.*

Indus Excerpt:

*Me: Yes or No? In Euclidean Geometry, we know that every straight line can be bisected.*

*Once again, I left out the 'line segment' bit but this time, I stopped added it in myself.*

*Me: In Euclidean Geometry you agree with me that every straight line can be bisected, right? That is something you learnt pretty early on in Euclidean Geometry – you have a straight line, you can bisect it. So, now if the world had exactly six points, would it be the same?*

*No Response*

*Me: Euclidean geometry – how many points (pointing to the line I had drawn). In Euclidean Geometry, how many points are there?*

*Multiple students: Infinite*

*Me: Yeah. Even on the line, there will be an infinite number of points. Right? Now, there aren't an infinite number of points. Only six points.*



> *Uday: Will there be straight lines between the points or just points?*
>
> *Me: There are only six points in this world. There's no between.*
>
> *Uday: then there can't be a straight line.*
>
> *Me: What do you mean by straight line?*

The tact I took in the two different sessions were quite different. In Ganga, I posed the question and didn't comment on it or help student with it as much initially. In fact, the instructions for their first group discussion was: what does a world with six points look like? In contrast, in Indus, I instantly pushed the students to clarify the concept of a straight line. This distinction is made clear by my summary of what had happened so far two minutes after this initial interaction in Indus. This summary was the result of a few students arriving late:

> *Me: The problem with the question is: I don't understand the question. Which is why I asked it.*
>
> *Uday: The actual problem with the question is there won't be any line segment.*
>
> *Me: We don't know. We don't know what line segment means. We don't know what a world with six points means. What does it mean? So, that's what we are trying to figure out. So, I thought, let's start with: what do you mean by straight line in Euclidean Geometry.*

In the case of the Triangle Theory Building module, the students already had experience with Euclidean Geometry. This module is much more of a departure from their school mathematics and hence it is harder for the instructor to strike the right balance between struggle and frustration.



### Clarifying the concept of straight line.

In both the schools, the first way the claim was attacked was by clarifying the concept of a straight line. In Indus, as mentioned above, this was because I explicitly asked them to do that. In Ganga, things were slightly different. Pankaj was the first student to try to clarify the concept of a straight line by saying that, '*...in every two points, there will be a line segment.*' Later, during the first group discussion, Aditya called me over to clarify the question – I had written 'line' in on the board and he said that even in Euclidean Geometry you cannot bisect lines, only line segments. I changed the question, but this led to a discussion on what a line segment is.

In Ganga, the first attempt to define a line segment came from Arnav who just said the phrase 'infinite points.' I responded:

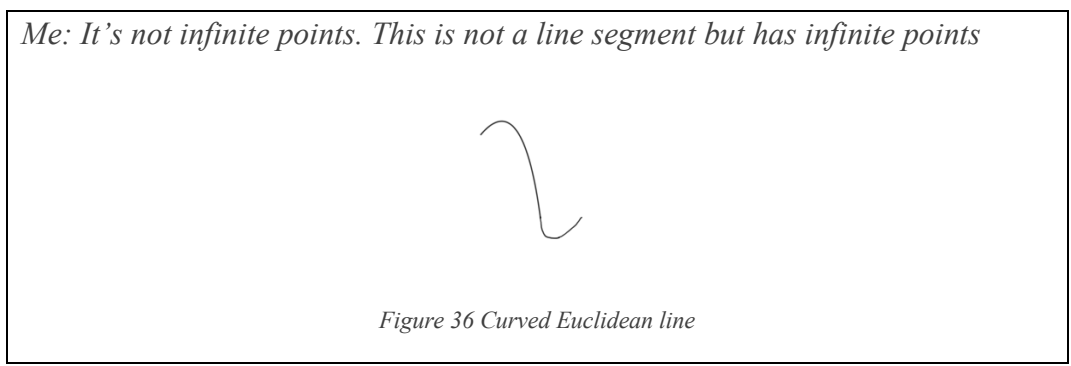

*Me: It's not infinite points. This is not a line segment but has infinite points*

*Figure 36 Curved Euclidean line*

The phrasing of my response here is not great. What I said can be easily interpreted as 'line segments do not have infinite points' rather than what I intended to say that 'the definition of line segment is not infinite points.' Modeling good reasoning and communication of that reasoning is a crucial component of helping students learn. While making mistakes in these things is inevitable, a way to at least notice the mistakes and correct them would be to write out the reasoning on the board. This is especially helpful in the first session of a module.

Immediately after that, Vivaan said that a straight line is the shortest path between two points. This is the most helpful conceptualization of straight lines to transition to a world with six



points. In contrast, in Indus, it took a while to get to this. Here is the initial exchange on straight lines:

> *Me: What do you mean by straight line?*
>
> *Uday: According to Euclidean geometry, there can't be a straight line.*
>
> *Me: So, what's a straight line?*
>
> *Uday: A line is a combination of points.*
>
> *Me: So, what makes it a straight line?*
>
> *Arun: Collinear points*
>
> *Me: What does collinear mean?*
>
> *Prerna: They lie in one line.*
>
> *Arun: They form one single straight line*
>
> *Me: But what's a straight line? So, collinear is something which lies on a straight line. A straight line is something with collinear points. There's a little bit of circularity there? So, what makes a line segment?*

Uday starts by suggesting that there cannot be straight lines. Assuming his next response reflects his reasoning, he seems to be saying that straight lines need have an infinite number of points and hence there cannot be straight lines when there are only six points.

Arun's response is something which I have heard repeatedly, including in the Triangle Theory Building session in this very course. He conceptualizes lines as a set of collinear points



and then says collinear points are points which lie on a line. One possible explanation for students doing this is a lack of understanding of the basic, undefined entities of a theory. At some point in the potentially infinite regress of asking what things mean, you have to stop and treat some objects as undefined and state certain rules about them. It may be the case that Arun, and others, do not have access to this tool as yet. It isn't even clear that those who define straight lines in terms of distance have this tool. In Euclidean Geometry, this would be circular since the only way of defining distance is in terms of straight lines and congruence. However, in geometries when we are able to define a notion of distance externally, these are not necessarily the same. For instance, if you have an externally defined metric imposed on a space, you can then define straight lines in terms of that metric. That is the tactic, without using the language of metrics, which conceptualizing straight lines as shortest paths allows us to use in a world with six points.

After this discussion in Indus, I proposed conceptualizing straight lines as shortest paths unlike in Ganga, where students did. In both the schools, I then gave a similar exposition on straight lines as 'shortest paths.' Here is a snippet of that from Ganga:

> *Me: So, now imagine I was an ant and I wanted to get from here (point at a spot on the board) to, and I can only walk on walls. The surface I'm walking on is this (pointing at the board) all the way till here (point at a chair in the class).*
>
> *What's the straight line between this point (point on board) and this point (point on chair)?*
>
> *(No response)*



*Straight line is shortest path. Let's make it more straightforward. Let's say I want to go from the top of this to the side of this (object is outside the shot)*

*Vivek: then the shortest distance will be directly between the two (and makes a gesture which looks like the direct path – through the air)*

*Me: But you can't go there. For the ant, the shortest path is this (I gesture but you can't see the object) Now suppose you took a flight between continents.*

*(I draw a sphere on the board)*

*I want to take the shortest path to save fuel and so on.*

*I want to take the shortest path. I'm going from India to Mo is in Argentina – on the other side. What is the shortest path?*

*Multiple students: through the Earth*

*Me: Through the Earth, right? And that is how they went. Right?*

*What will the shortest path be which they can actually take? It will be on the surface and will be something like this (and I draw something approximating a great circle). So, that is the straight line path along the surface of the Earth. Just as this is the straight line path from this point to this point along this surface. So, now if you think of straight line as shortest path, then this might make a little more sense.*

What I was trying to do here was to take the concept of straight line we had discussed earlier and explore its consequences in different contexts. The goal of doing this was to move



students away from their visual intuition of what a straight line is and follow the consequences of the definition rather than their concept images.

However, this could have been done in a more structured manner, going slowly through multiple examples of different contexts. I have laid out an alternative exposition in the appendices as part of the lesson plan on this module.

**Concept of straight line.**

In both the schools, students were unable to come up with a way to conceptualize worlds with six points. This was even after the straight line discussion. Hence, I stepped in and helped with an analogy. I used the idea of teleportation. I did this in two slightly different ways in the two schools. In Ganga, I did this by discussing travel from a spaceship millions of light years away with a teleportation device to a point on the Earth. Suppose you were on somewhere on the spaceship and you wanted to get to a spot on the Earth. The shortest path would involve walking to the teleportation device on the spaceship, taking it to Earth, and then travelling on the Earth to get to wherever you wanted to go. I then gave the students a group task to try to use this analogy to understand shortest paths in six-point worlds.

In Indus, I directly introduced them to teleportation in the situation on six-point worlds. I said that suppose three of the points represented Mars, Delhi and Pune. Suppose we had access to teleportation devices in each of these three places. However, the teleportation devices did not allow you to go from Delhi to Pune, but they allowed for instantaneous travel along the other two routes. Now, what is the shortest path from Delhi to Pune? It is the path from Delhi to Mars and then Mars to Pune.

After the introduction to the teleportation analogy in Ganga, I drew the following, said that the connections represent teleportation being allowed between the two points, and asked



about various shortest paths. While doing so, I said that for now we can assume that each individual teleportation takes the same amount of effort or, in other words, is of the 'same length'. I asked for the shortest path between A and C, and E and F.

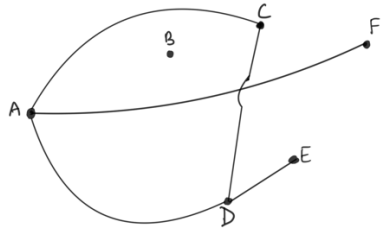

*Figure 37 Ganga six point world*

In Indus, I drew the following simpler diagram and took them through a few examples of shortest paths.

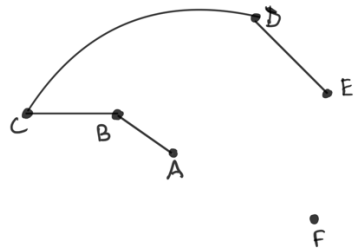

*Figure 38 Indus six point world*

The advantage of the more complicated diagram used in Ganga is that there are more often more than one path connecting two points and hence the shortest path is not the only path. In the simpler example, paths and shortest paths mean the same thing. Given that this is something which comes up again and again in this module, it would have been a good idea to spend time on it before jumping into more complicated questions.

After this, in both the groups, I said that we should concentrate on a simple world which looks like this, where each point apart from two end points have exactly two neighbors, and the end points have one neighbor.



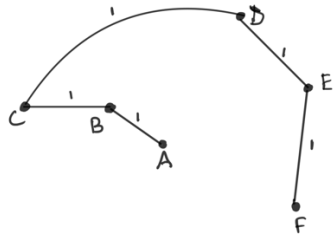

*Figure 39 Simple world*

I asked them to answer the question I had started with, about bisection, in this world.

Before going into that, it is useful to think about the representation I used. In both the schools, I used graph-type edges and vertices to indicate that two points were neighbors. I made sure to put the points on various parts of the board and not make the connections exactly as long as each other. The goal of doing this was to get students to differentiate physical proximity in the representation from proximity in the space. To make this clearer, it could have been a useful exercise to introduce alternative representations. For instance, putting two points right next to each other like you would with pixels on a screen to indicate adjacency. Another possible representation is a table like you use to represent graphs. Having multiple representations could have allowed students an easier escape from their intuitions about geometry.

**Engaging with the bisection question.**

The two approaches for dealing with this question I took were different enough that it is worth spending time on each of them separately before comparing them.

***Ganga.***

The discussion straight after returning from their groups was initially mostly a dialogue between Vivaan and me. There were a lot of different things which came up during the dialogue, some which were addressed and others which were not.



*Dialogue with Vivaan.*

> *Me: You guys had something (pointing to the group on the right)*
>
> *Vivaan: Bisection mean that if a shortest path can be broken into two parts by any other path, that means the path is bisected.*
>
> *Arnav: By a shortest path or by any other path?*
>
> *Vivaan: Any other path*
>
> *Arnav: Not shortest path*
>
> *Vivaan: Not necessarily. Shortest path can be divided into two by any other path.*

This is an interesting exchange between Arnav and Vivaan. It points to Arnav thinking carefully about the definition. He is interrogating a very specific aspect of the definition. It isn't clear why he thinks this is an important distinction – did he have a specific example of a situation in mind where this distinction would make a difference? I have not found anything in the preceding group discussion regarding this. Also, this was not something taken up for discussion in the class since Vivaan continued with the following:

> *Vivaan: The other path doesn't contain the same two points. Otherwise, every line can be bisected.*
>
> *(Vivaan gets up to come to the board)*
>
> *… what I'm asking is divided into two parts by any other path.*



*(He draws the following)*

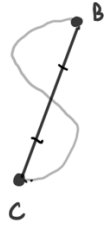

*Figure 40 Vivaan's drawing*

Vivaan drawing this seems to point to him not having completely understood at least one of the setup or the representational system. He seems to be thinking of the edges indicating a connection as somehow divisible entities rather than just a representation indicating a particular relationship. This points to the need to spend even more time on the representational system and on developing various alternative representational systems before continuing with the question.

I dealt with this in the following manner:

*Me: So, there's two paths between these two points. And somehow they meet?*

*Vivaan: Yeah*

*Me: But, what does that meeting mean?*

The reason I asked this was for the students, and Vivaan specifically, to realize that what he was doing here was introducing a new point, and hence the world no longer had six points. However, this did not happen since Vivaan took the discussion in another direction by asking if there can be more than one connection between two neighboring points. I went along with that even though that doesn't deal with the issue at hand directly even though it is tangentially



related. Alternatively, I could have opened up the question I asked to the rest of the class and stopped the dialogue I was having with Vivaan.

*Me: We had said, at least for now, all the paths would have the same length, whatever you want to call it. So, this path and this path (referring to the two paths between the same two points) have the same length. At least, that's what we had agreed on. And even this path and this path and this path are the same length (referring to other paths between neighbors). Whether they are the same path or not is a question you need to answer. You can choose.*

*Vivaan: Are all paths between points the same. Like, how do you define shortest.*

*Me: So, that's something you can decide, right? So, say there were two paths between A and B. So, think about it like two teleportation machines, one right next to the other. They have the same length – we agreed on that. So, you can count these as two separate paths or the same path.*

*Vivaan: If we count them as the same path, then what do we mean by shortest. There is only one path – that should be the shortest path. Then all paths will be the shortest path.*

Vivaan seems to be saying that if we do not allow more than one path between two points, every path in this world will be a shortest path. That is reasonable. However, adding in another path (of the same 'length') between neighbors wouldn't change that. So, in retrospect, I am uncertain what he saw as the connection between this and the question about two paths between neighbors. In the moment, I did not push him to make that connection. Rather, I picked



up on his concern that 'every path is a shortest path' by returning to the world I had initially created.

> *Me: So, on this one, definitely, shortest path does not mean… all paths are not shortest. (Referring to the original world with six points I had drawn)*
>
> *Vivaan: I don't know*
>
> *Me: No, on this one. ACD is longer than AD.*
>
> *Vivaan: Yeah*
>
> *Me: If you assume they all have the same length. So, in this one, path does not mean shortest path. In this one (referring to the newer diagram), you are right. Then every path is a shortest path. But even if you add this in (the extra path between two neighbors), both of these are now the shortest paths.*
>
> *Vivaan: So, then that definition is meaningless in this case.*
>
> *Me: What?*
>
> *Vivaan: Shortest path divided into two paths by any other path.*

There are two interpretations I have been able to come up with for the term 'meaningless'. The first is that the definition makes the term a synonym for another term. An example of a definition that could be thought of as this type is the definition of straight line in the world we were investigating. Since all paths are shortest paths, all of them are straight lines.



Hence, straight line is just a synonym for path[4] and doesn't serve an alternative purpose. Another possible interpretation of the term 'meaningless' is that, given a particular definition, there are no examples of the object in world we are dealing with. An example of this would be if you defined bisection as Vivaan did but did not count a single point as a path. In that case, in the particular world we were dealing with, no straight line can be bisected.

A reason for believing that he was raising the latter is that he explicitly mentioned the definition of bisection when I asked what he thought was meaningless. We do not have any obvious synonyms for bisection so it is unlikely that he was referring to the synonym objection.

However, straight after this excerpt he says:

*Vivaan: All paths are shortest, so any path can be divided by another.*

This could indicate that he was not referring to the definition of bisection, but to the definition of straight line as shortest path and its impact on the definition of bisection. He had raised the same issue about straight lines in the previous excerpt presented when he was concerned about all paths being the shortest path. I could have investigated this further with Vivaan in order to figure out what his concern was, but I instead decided to open up the discussion to the rest of the class.

There were a few different things which came up in this dialogue with Vivaan. While the dialogue lasted for a while, I seem to have failed to understand a lot of what he was saying in the moment. Alternatively, I could have concentrated on one or two particular aspect of what he was saying and tried to get him to clarify his ideas. This could have led to a more productive exchange which would have been of benefit to the entire class.

---

[4] The truth of this depends on your definition of path. If you allow for paths which go back and forth (eg. A to B and then back to A on the same connection), then paths and shortest paths could mean two different things.



*Returning to the Group.*

Me: Let's pick one path. Let's pick this path. Path BCD. Can this be bisected?

Karan: Yes. C is where it can be bisected.

Me: So, you are saying it can be broken up here. What is it broken up into?

Karan: BC and CD.

Vivek: BC equal to CD

Vivaan: What about BC. Can BC be bisected?

Arnav: No. There exists no point between B and C.

Me: So, if it has to be bisected by another point or another path, then, this cannot be bisected.

Vivaan: Then, not every line can be bisected.

Me: So, I guess that by that definition of bisection, not every line can be bisected.

This is the first exchange where students have the tools at their disposal to be able to directly address the question we had started with. They have a definition of bisection and of straight line, and they have a conception of what the worlds look like. What I decided to do then was to pivot from the original question and ask a related question.

Me: ...what sorts of things can be bisected?

Arnav: A path which passes through another point.



> *Me: Okay. So, let's try that. Let's say the path that passes from A to D. So, A, B, C, D. Can this be bisected?*
>
> *Multiple students: No*
>
> *Vivaan: A path which passes through an odd number of points.*
>
> *Me: So, A,B,C,D,E. Can this be bisected?*
>
> *Multiple students: Yes*
>
> *Arnav: C is the bisector.*

In this exchange, I decided to give a counter-example to Arnav's claim that paths which pass through other points can be bisected. I'm assuming here, and seemed to assume in the moment, that 'other point' means a point apart from the two end points. Rather than giving the counter-example myself, I could alternatively have opened up the claim to be evaluated by the entire class. This was a claim that students could have easily come up with counter-examples to.

I then suggested that students come up with an alternative definition of bisection which allows for paths with even points to be bisected.

> *Me: What would that notion of bisection be?*
>
> *Vivaan: Whether it can be divided into two separate paths which do not share a point in common.*
>
> *Me: So, if you have ABCD, that could be broken up into AB and CD.*

As I wrote in my notes, this alternative definition is something students have come up with on their own when I had asked this question before. Since it didn't happen automatically, I



pushed students to come up with it. The reason for doing that is that one of the goals of this module is for students to understand the relationship between definitions and conclusions. Seeing that requires multiple definitions of the same term at hand.

*Worlds where all straight lines can be bisected.*

I then gave the students a group task to come up with worlds with six points where every straight line can be bisected by each of the two definitions. In the case of 'bisection at a point', this is only possible in a world with no paths. This is something all the groups came up with and we discussed after. However, the case of 'bisection between points' is more interesting.

Arnav came up with something, not directly answering the question, but in order to clarify something:

*Arnav: So, what we can do is – if we have (this – he draws A-B-C-D-E). What we can do is break this up here and he writes:*

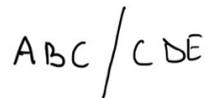

*Figure 41 Arnav's drawing*

*Me: So, you are saying that if you add this path, this is what you guys were talking about (in your group), then you break it up like this.*

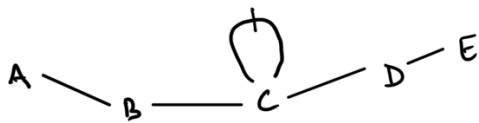

*Figure 42 Reporting group work*



His claim seems to be that if you add in a path from a point to itself, then paths with an odd number of points can be bisected by this notion of bisection. I didn't take this up for discussion. Rather, I chose to say that for now we won't allow this. However, there is something worth thinking about in Arnav's idea. Arnav's idea is challenging the particular claim made earlier about odd and even number of points. This is an example of a path with an odd number of distinct points which can be bisected. However, there are at least three possible ways to deal with this objection without losing the original claim completely. The first is to say that the example is not a straight line – it is not a shortest path. Our definitions of bisection involved shortest paths. Another way is to say that the path Arnav gave as the example may have five distinct points, but has six points on it since you travel to the point C twice. A third way would be to change the original claim to be about connections between points rather than points themselves. In that case, paths with an odd number of connections can be bisected by this definition but not paths with an even number of connections.

It would have been interesting to bring up each of these, either encouraging students to do so or doing so myself, rather than just rejecting the premise that points can have connections to themselves.

Vivaan came to the board next to represent his group's discussion. He said:

> *Vivaan: Let any two points. In this world. Let's join these two points. If two points are joined, we call a set and if we can break up that world into several sets such that no two sets share a point in common, then we can bisect any line.*
>
> *Now, I'll prove that no other world can be bisected. For example, if you have two set, there is a point in common, then it cannot be bisected by this definition because… for example.*



*Sabareesh: By a set you mean any two points.*

*(Sabareesh, as mentioned in the methods section, is an undergraduate student who attended a few sessions in the schools)*

*Vivaan: Yeah. This is a set and this is a set. They share a point which cannot be bisected. Now if we add another point, it can be bisected into this and this, but there is a subset of those four points where it can't be.*

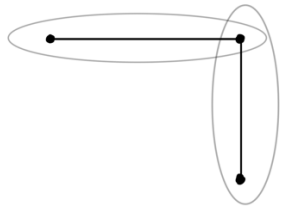

*Figure 43 Vivaan's 'sets'*

Vivaan makes an interesting move here by giving the object he is interested in a name – he calls any pair of points with a connection between them 'a set'. This allows him to use that notion without elaborating on what he is describing each time. The other valuable thing Vivaan does is to attempt to show that all such worlds allow for every line to be bisected and only such worlds allow for every line to be bisected. Even though it turns out that he was wrong about the latter, as we will see in the next excerpt, the important part is that he sees the need to justify both these claims.

*Sabareesh: Consider this world. Then, what are your sets?*



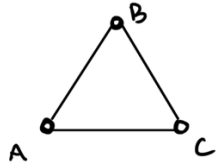

*Figure 44 Sabareesh's world*

*Vivaan (pointing at each of AB, BC and AC): This is a set. This is also a set. And this is also a set. And now, none of this can be bisected.*

*Me: So, you are not saying none of this can be bisected, right? You are saying that there are lines which cannot be bisected.*

*Vivaan: There are lines which cannot be bisected*

*Sabareesh: Which lines cannot be bisected?*

*Vivaan: This one. (apparently pointing to a line with three points, say A-B-C)*

*Sabareesh: But, this is the shortest path (apparently pointing to the direct connection A-C)*

*Vivaan: But, it didn't say shortest path. It said can any line be bisected, not any shortest path.*

*Me: But what did we mean by that?*

As I wrote in my notes, I hadn't thought of this example which Sabareesh brought up. However, the value of this example is that it highlights the need to carefully interrogate the connections between definitions and conclusions. Even though not every path in the world which



Sabareesh drew can be bisected, if we restrict ourselves to straight lines, that is no longer the case. This is because all straight lines consist of exactly two points.

There potentially was a little bit of confusion in the exchange above due to the use of the word line. I was using the words straight line and path. However, the term line is often used synonymously with the word straight line and that is the sense in which Sabareesh seems to be using it.

This exchange was between Vivaan, Sabareesh and me. I decided to go through it once again for everyone else. I drew the following:

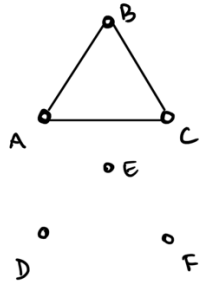

*Figure 45 World where all straight lines are bisectable*

*Me: What are the straight lines?*

*Multiple students: AB, BC and CA*

*Me: Is ABC a straight line?*

*Arnav: No*

*Me: Why not?*

*Arnav: It is not the shortest path.*



> *Me: It's not the shortest path between A and C. So, ABC is not a straight line.*
>
> *So, these are the straight lines. Can every straight line be bisected?*
>
> *Vivaan: Yes. Every straight line can certainly be bisected.*
>
> *Me: Yeah. But not every path. Because this is a path which can't be bisected.*

This ended the discussion on the bisection question. This strategy I used at the end of taking something which was the result of a dialogue and going through it once again with the rest of the class, is something I could have used in the earlier dialogue with Vivaan. In classes like these, it is inevitable that there will be conversations where only one or a few students are involved. Doing this is way of allowing for that while not leaving the rest of the class behind.

**Indus.**

*Examples of lines which can and cannot be bisected.*

I asked the students at Indus to have a group discussion to decide on whether all lines in the world can be bisected. When they came back, I asked who believed that every line can be bisected and who believed every line cannot be bisected. The students who responded said that every line cannot be bisected.

> *Me: Give me an example of a line which cannot be bisected?*
>
> *Uday: DE*
>
> *Me: D-E. The line containing D and E cannot be bisected. Why not?*
>
> *Uday: Because there is no point to bisect them.*
>
> *Me: Oh. So, what do you mean by bisection?*



*Uday and others from his group: Divided into two equal parts.*

*Me: But, here is my division into two equal parts. DE. I just.. two points (while making an action conveying chopping into two)*

*Tarini: But, there's nothing in between. There is no other point.*

*Me: Oh, so you want to bisect at a point? You have to have a point to bisect at?*

*Tarini: Yes*

Uday's example illustrates that he was using the notion of bisection at a point. However, his group's definition did not convey that. The definition just says 'breaking into two.' From Tarini's response to my concern and Uday's reasoning about there being no point to bisect DE, it appears that they had grasped the discrete nature of the world we were working in – that between two points there need not be other points.

In the above extract, I asked the students for the definition of bisection. They didn't explicitly state the need for it themselves. A possible strategy for getting them to do that, in this situation where there seemed to be agreement, could have been for me to take the other stand – that all lines can be bisected. I could have just asserted that without giving an argument. In that case, there is a chance that they would have noticed the need for clear definitions.

Straight after the previous excerpt, Tanya said:

*Tanya: You said the distance is not similar. So, you cannot bisect that also.*

*Me: EF and ED we said is one.*

*Tarini: You said the time taken, the effort.*



*Tanya: The distance is different.*

*Tarini: You didn't say it.*

*Me: So, the question is: what do we mean by distance?*

*Tanya: We want to know whether the distance between DE and EF is the same.*

*Me: It depends on what you mean by distance. So, by distance, you usually mean what? So, in Euclidean,, the length of what? So, is this the distance?*

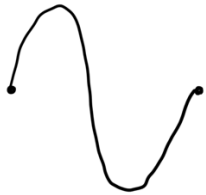

*Figure 46 Length of curved line in Euclidean geometry*

*The length of that?*

*Multiple students: No*

*Me: The length of the shortest path is the distance. What did we mean by length? What was the length of DE? One. We are just taking it as one.*

*Tarini: So, it is not the effort that is one?*

*Me: So, now the question is: when you say length, what do you mean in this situation? Have you seen a map of a metro, like an underground thing, like Pune metro is starting. They look something like this:*



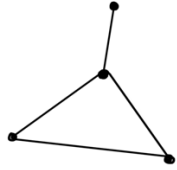

*Figure 47 Metro map*

*They don't represent any sort of Euclidean distance. Otherwise, the map would look really weird. It would look something like this:*

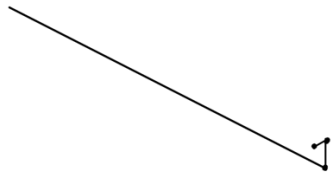

*Figure 48 Euclidean metro map*

*In this map (the first one) this connection just means that this is the next station. The connection on subway lines on a maps mean that it is a neighbor, next station. So, when I ask you 'when should I get off?', you tell me 'get off after five stops.' You won't tell me 'get off after 800 kilometers.' So, what's the useful notion of distance here? The useful notion of distance here is that we say 'we are counting the number of stops.' On your teleportation. And we say we spend about the same amount of effort on each stop. Maybe later we can start changing that and say that this path is more troublesome to get through than this path and so on.*

This excerpt is not directly about bisection but is related to the question of bisection via the notion of distance. I had previously stated that for now we can assume that the effort required to teleport between neighbors was the same and let us call that one. I had also explicitly used the



word length when saying that. However, it seems like Tarini and Tanya did not see the association between length and 'effort' I had tried to make. Given that Tanya seemed to understand the discreteness of the space given the last excerpt, it is unlikely that the issue was to do with discreteness. It could have been to do with my communication of the concept of distance we could use in this world. Examples like the metro-map example and other similar examples could have been brought up when I introduced these ideas.

*Types of lines which can and cannot be bisected.*

The next thing I asked students to figure out was which types of lines can and cannot be bisected.

---

*Me: So, you showed me one line which can be bisected and one line which can't.*

*Tanya: AB*

*Me: So, you are saying AB can or cannot?*

*Tanya: Cannot*

*Me: So, AB cannot. DE cannot. So, what sorts of lines cannot?*

*Tanya: Lines which have points in between*

*Serena: Odd number of points*

*Tarini: Lines which have odd number of points can be bisected.*

*Me: Ah. So, you are saying lines with odd number of points can be bisected because there is something in the middle.*



When Tanya said 'lines which have points in between' can be bisected, that could have been taken up for discussion as it was in Ganga. However, in this case, Serena and Tarini immediately corrected her. Even so, I could have taken that claim and asked for counter-examples if students thought it was false. Even though rigorous justification is not the main goal of this module, this was an opportunity to bring it in naturally. This could also have been a way to get other students involved in the discussion.

I ended that part of the discussion by asking students whether straight lines with an even number of points cannot be bisected. They said yes and I didn't push them for a justification. Rather, I introduced the other concept of bisection.

*Me: Here is an alternative notion of bisection. You have A, B, C, D. They are neighbors of each other. And, we break it up into AB and CD. This is another possible definition.*

*We broke it up into two equal lines, but just now we are not doing it at a point.*

*Multiple students: That is not bisection*

*Uday: Break it up into two equal parts and no part should be left out*

*Me: That's what you want (pointing at the original definition of bisection). You want no part to be left out. Say you did this sort of bisection. Now, what sorts of lines can be bisected? Let's call this Bisection-B. Now, what sorts of lines can be bisected and what sorts of lines cannot?*

*Tanya and Devika: Odd number of points can be bisected.*

*Me: Odd number? This has an even number.*



> *Tanya: All lines can be bisected.*
>
> *Me: So, how would this be bisected?*
>
> *Devika: By the point.*
>
> *Me: So, if you put these together. I'm saying if it was just this type of bisection.*
>
> *Uday: Then even ones*
>
> *Me: Then odd ones cannot and even ones can*

It seems like the students who responded saw the new definition not as an alternative, but as an addition – bisection is something which falls under either of these two definitions. This idea of having multiple definitions and following the consequences of each one is at the heart of this module. In this group, unlike in Ganga, I gave the second definition a name – Bisection B. In retrospect, I would have given it a more descriptive name – Connection-Bisection for instance. I didn't spend any more time on this question. Rather, I gave an exposition recapping what we had done so far. Then, I asked students to find worlds where all lines can be bisected, similar to what I had done in Ganga.

*Worlds where all lines can be bisected.*

> *Me: So, now that we've created this world, now let us try and explore different worlds of this type. You can have worlds like this:*



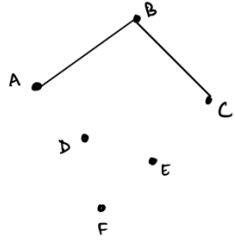

*Figure 49 A six-point world*

*D, E, F are isolated and A, B, C are connected. So, can you create a world, which by this definition (the original definition of bisection), every straight line can be bisected?*

*Worlds where every line can be A-bisected (and I added in the letter A to where the original definition of bisection was), and B-bisected.*

This was given as a group task. These instructions are quite ambiguous. It isn't clear whether I am asking for a single world where all lines can be both A and B bisected or whether I am asking for worlds where all lines can be A-bisected and worlds where all lines can be B-bisected. While I meant the latter, from what I wrote, the former is probably the more likely interpretation. However, from listening in to the group discussions, it seems like students interpreted it in the latter way, the way I wanted them to. After returning to the class discussion, I asked:

*Me: Can you create a world where all lines can be A-bisected?*

*Multiple students: No (or nodding their heads to indicate no)*

*Me: You guys said yes (pointing to Tarini's group)*

*Tarini: We said no*



*Me: You had one*

*Tanya: The only points one?*

*Me: Yeah*

*Multiple students in the group: Only points*

*Me: Okay. So, look at this world. Can every line be A-bisected? In this world can all lines be A-bisected?*

*Uday: There are no lines.*

*Tanya: That is the point. There are only points.*

*Me: You could think of a single point as a line with the same end points, potentially. Otherwise, if no lines exist, then of course A-bisection. On all lines which do exist, it works – there are no lines.*

In the group discussion, Tarini and Tanya had come up with the example of a world with no connections as a world where all lines can be A-bisected. I am uncertain as to why they didn't share it when I asked till I directed the question to them. I am not sure whether students had seen the idea of vacuous truth before. If they had not, this is unlikely to have been a sufficient introduction to the idea. I continued with the question on B-bisection:

*Me: So, what about B-bisection*

*Tanya: Same world*

*Me: Same world?*



*Tanya: Yes. Same world*

*Me: Okay*

*Uday: A world with two points also*

*Me: Two points in what sense? Two connected points?*

*Uday: Two connected points.*

*Tanya: We cannot bisect that.*

*Me: Why? You can break A-B into A and B. But, there is a question here. Because, bisection does what? It breaks a line into lines. The question is: do we want to treat these as lines?*

*Tanya: We can*

*Me: We can. So, let's call this Line A – includes single points. Line B does not. And, Line A kind of works with our idea of shortest path. What is the shortest path from this point to itself? Stay there. Don't move. So, it is reasonable. So, with this notion, Line A, does this work?*

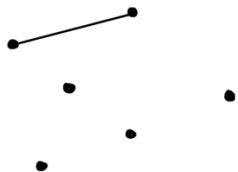

*Figure 50 World where all lines can be bisected*

*In this world, can every straight line be bisected?*



> *Uday: That one cannot. There is no point in between.*
>
> *Me: Sorry. Can every line be B-bisected? Let me rephrase that: can every Straight Line-A be B-Bisected?*
>
> *Multiple students: Yes.*
>
> *Me: So, this line can be B-bisected (pointing at the two connected points). What about this is now a straight Line-A (pointing at a point), right? Can C on its own be B-bisected? What are we breaking it into?*
>
> *Multiple students: Same point*
>
> *Me: But, B-bisection, we have to break it at a connection.*
>
> *Multiple students: Then no.*

Here, I continued with giving different names when I change a definition – Line A and Line B. This did make the question easier to phrase. However, once again, descriptive names might be better since they would probably be easier to remember. The issue here is that we are trying to create a world where all lines can be B-bisected. However, by one notion of straight line, you cannot break up a two-point line into two straight lines and by the other notion you introduce new lines which you cannot B-bisect. This is interesting as we cannot change the definition of straight line in such a way that we can bisect all lines. This did not come up in Ganga. In Ganga, we accepted that we could break up a line into points. That is where I went next:

> *Me: By line-B, you can B-bisect it, but you don't get a line out if it. You get points. Everyone sees that? If you choose Line-B which doesn't allow single*



*point, you can 'bisect' this, but you get two points. Not two lines. Let's allow it. Just for now. Assume we allow for a line to be broken up into two points. Then, this can be B-bisected. And all lines can be B-bisected. Assuming no point in a line.*

*Tarini: Yeah.*

*Me: Right? There is only one line and that can be B-bisected. Okay. So now, what sorts of worlds can be B-bisected? What if we connect this (two other points)? Can every line be B-Bisected?*

*Multiple students: Yes*

*Me: There are these types of worlds where two points are connected. Can we connect three things together?*

*Multiple students: No*

*Me: Is there no example of a world where three points can be connected together? In some way and everything can be B-bisected.*

*Multiple students: No*

*I then asked Sabareesh to come to the board to give the same example he gave in Ganga.*

*Me: In this world, can every straight line be B-Bisected?*



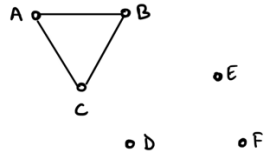

*Figure 51 World where all straight lines can be bisected*

*Tarini: No. ABC.*

*Me: What is a straight line?*

*Tushar: AC then*

*Tanya: Shortest path.*

*Tushar: Yeah. AC is the shortest path*

*Me: So, Sabareesh's world seems to work. The shortest path is AC. That (ABC) can't be bisected but that is not a straight line. So, can every straight line B be B-bisected?*

*Multiple students: Yes.*

*Tanya: Even if you join DEF.*

In the above extract, Tushar spotted that all lines in Sabareesh's world can be B-Bisected. Tushar is somebody who hadn't spoken much till this point in any of the sessions.

The strategy I used throughout this session with Indus of giving names to the different types of lines or bisection we were dealing with, is something worth exploring further. I left the bisection part there and moved on to discussing other shapes in these worlds.



**Summary of the initial question**

It is valuable here to reconceptualize the initial question in terms of the questions of theory building and the concepts of theory building outlines in Chapter 3. The question was phrased in such a manner that it seemed like it was meaningful, but it was not. It was not something students could respond to within frameworks they were familiar with.

By asking them questions like, 'what is bisection in Euclidean geometry?' or 'is this true in Euclidean geometry?', the goal was that they take concepts from Euclidean geometry and look to extend them to a context where they currently do not exist. Taking the concept of straight line as an example, it is not defined in Euclidean geometry. However, we do conceptualize it in a few different ways. Two conceptions of straight lines are 'paths in the same direction' and 'shortest paths.' On a flat, gradient, plane both of these have the same consequences. However, in a world with six points, they do not. In fact, in a world with six points, it is not clear as to what you could mean by direction or angle. However, shortest path is something which may be conceivable in such a world. Hence, it is a good candidate for extending.

There were also often multiple possible directions which could be taken for definitions. For instance, bisection had two definitions with two different sets of consequences. In these sessions, we did not choose between them. Rather, we observed the consequences of both the definitions and saw the relationship between definitions, even those which are similar, and their consequences.

While we were not explicit about axioms, they resided in the background. The concept of neighbor, of distance between neighbors, of points, and so on, were constrained by axioms which were implicit in the question and in the analogy which I used of instantaneous travel. If this



module were to continue for a longer period of time, the focus could move from exploration of these different worlds to making these assumptions explicit, and formalizing the theory.

In the next chapter, I will discuss other episodes in the course as well as student feedback. The chapter will end with general observations and reflections from the entire module during which I will make references to the data in this chapter as well.



## Chapter 9 – Discrete Geometry: Other episodes and feedback

This chapter is a continuation of the previous one. It focuses on the discrete geometry module and consists of three parts. The first part of this deals with the implementation of the rest of the module after the set-up discussed in the previous chapter. In it, I present entire chunks of the module. When I say entire chunks of the module, I decide on this based on my coding of mathematical content which I describe in the methods chapter. I have picked chunks of the module which I judged to be of interest and which are somewhat representative of the rest of the module. I decided on this by focusing on parts I mentioned in my teacher notes along with those chunks which covered a large number of codes based on the coding I did of students' mathematical actions (such as conjecturing, proving, defining, etc.).

The second part of this chapter deals with students' reflections on the module. There are two different types of student reflections I have access to. The first is in the form of written feedback at the end of a particular session. In the case of discrete geometry, I have feedback on the first and last sessions in both the schools. The second is student reflections on the course as a whole which they gave on the last day of the course. My focus in this chapter is on parts of that specifically related to discrete geometry.

The final part of the chapter consists of general observations related to the research questions. As I mentioned, the main focus of my analysis is local to specific parts of the course. However, there are a few general themes worth exploring. This includes discussion of data presented in this chapter as well as the previous one.

Here is a table of contents for this chapter:







## Other episodes

I will be discussing two episodes from Indus and one from Ganga. They are:

1. Initial exploration of circles (Ganga)

2. Initial exploration of triangles (Indus)

3. Straight lines and line segments (Indus)

### Initial exploration of circles – Ganga.

We started discussing circles at the end of the first day after finishing the discussion on

bisection. I asked:



*Me: What is the definition of a circle? Set of points…*

*Arnav: Equidistant from a given point.*

*Me: Set of points equidistant from a given point which we call the center. Let's say this is my center. (pointing at A in the diagram below)*

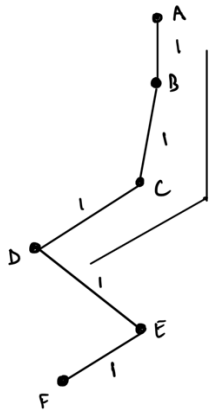

*Figure 52 Circles in a simple world*

*As we said, let's just call this distance one. All of these are distance one.*

*Vivaan: Then AB is a circle.*

*Me: Is A part of the circle?*

*Vivaan: No*

*Me: What do you mean by a circle? When you say a circle with center O, is O part of the circle or no?*

*Two students (unidentified): Yes*



> *Me: So, what is a circle? A circle is the set of points at some distance from the center. Is the center, let's say this is one (drawing a representation of a Euclidean circle with radius 1), is the center at distance one from the center?*

What I'm attempting to do here is to borrow the definition of circle from Euclidean Geometry and see how it applies in the worlds we created through answering the bisection question. In the exchange above, before I could ask a question about what circles look like in this world, Vivaan said that AB is a circle. Rather than asking why he thought that, I seem to have interpreted what he said as that the points A and B lie on the circle.

However, once I interpreted it that way, it is interesting to see that at least two of the students judged the center of a circle to be part of the circle. Given the usual definition of a circle in Euclidean Geometry, the one we discussed in this extract, the center is not part of the circle since it is not a radius distance from the center. The discussion continued with Vivaan saying:

> *Vivaan: It depends on context. If we want the center to be a part of the circle.*
>
> *Me: Of course. Usually in Euclidean geometry, is it? What is a circle? Repeat the definition. Set of points…*
>
> *Multiple students: given distance from the center.*
>
> *Me: So now, what's a circle with radius one? The set of points at a distance one from the center. Is the center distance one from the center? In Euclidean Geometry? Yes? No?*
>
> *Karan: Yes.*
>
> *Me: What is the path of length one from the center to the center?*



*Karan: Center to the center or center to any given point on the circle?*

*Me: Center to any given point on the circle has distance one.*

*Karan: Yes.*

*Me: We are asking if the center is a part of the circle?*

*Karan: No.*

*Me: Why not?*

*Karan: Because in the definition of a circle, circle is a set of points which are equidistant from the center. We don't say including the point from which you…*

*Vivaan: The center is not equidistant from the center.*

*Me: Of course you could say that a circle is the set of points equidistant from the center and the center. Add in the center to the circle. If you want to.*

It seems like in this exchange either Karan did not understand what I asked or he changed his mind through the questioning. One interesting thing I did not pick up on during the session was that Karan said that a circle is 'a set of points equidistant from the center' rather than 'the set of points equidistant from the center.' By the first definition any subset of a circle would also be a circle while with the latter that is not the case. It isn't clear whether Karan actually meant to say 'a' rather than 'the' but it could have been a good opportunity to discuss the distinction.

At the end of the excerpt, I said that you could define a circle such that the center is a part of the circle by stipulating that. I did not take that up for discussion later, but it could have been a useful arena to discuss choices of definitions. The problem with defining a circle in such a



manner is that for various theorems, you would have to add in the caveat 'except for if the point in question is the center.' Also, at least as far as I'm aware, there aren't interesting theorems about all points on a circle if the points include the center. We then went back to discussing circles of radius 1 in the world we were talking about.

*Me: What is the circle with radius one with center A? What is the definition of a circle of radius one?*

*Vivek: A set of points equidistant from A with length of one.*

*Me: So, what is the set of points equidistant from A with length one.*

*Multiple students: B*

*Me: B. Just B. So, the circle with radius one and center A is just the point B. What about the circle of radius four with center F. Center F and radius four?*

*Karan: B*

*Me: B. Because it is the only thing at distance four away from F.*

It seems like many of the students were able to deduce the consequences of the definition of a circle in the world we were working in. It could have been interesting to point out here that in this world, there are two circles with exactly the same points but two different centers. This is not true in Euclidean geometry where the points on the circle determine the center. However, I did not do that this since Vivaan interrupted with:

*Vivaan: What if we connect F directly to B. Then it will be distance one.*

*Me: Then it will be distance one.*



*Vivaan: Since all paths are equal.*

*Me: Yeah. Now if we do this (I connect F to B on the diagram), is there something of radius four? Is there a circle of radius four? With center F.*

*Vivaan: No.*

*Karan: Yes. If we take from F to B. Then B to C. Then C to D. then D to E. Then that will be four length.*

*Me: But the shortest path form E to F is one.*

*Vivaan: So, radius should be shortest path?*

*Me: You said distance. So, what do you mean by distance?*

*Vivaan: Distance is shortest path. Distance between two points is the shortest path.*

*Me: Yeah. Distance between two things is the shortest path, right? So, if I ask you: what is the shortest path from here to Bombay? You won't go like via Moscow. That's also a path, right? But, it's not the shortest path.*

Vivaan suggested creating an alternative world here where we have a direct connection from F to B. Given that connection, the set containing only B is no longer a circle of radius 4 with center F since the distance from F to B is now one.

When I asked whether there was still some circle of radius four with center F, Karan responded that there is one. What he suggested was a path from F to E of length four which would result in E being distance four away from F. The problem with this is that there is a



shorter path from F to E, namely the direct path. When I pointed that out, Vivaan asked whether radius should always be the shortest path. I responded by saying that radius is to do with distance and it would be strange to use distance to refer to anything other than shortest path, giving the example of a path from Pune to Delhi via Moscow. While this is a good example to say why distance should be the length of the shortest path, it does not say why we should consider radius to be distance. A good way to show that would have been to return to a Euclidean circle and show a curved path of length one from the center and ask whether the point at the end of that path is part of the circle.

That was the end of the conversation for the day. We continued with circles in the next session.

### *Summary of the initial exploration of circles.*

In this conversation, the goal was to borrow the definition of circle from Euclidean geometry and apply it in the worlds we were concerned with. I chose a specific world to explore first – one where each point had exactly two neighbors apart from the two end points.

The idea of a closed shape came up. It isn't clear what that exactly means here, but the 'circles' in this world do not fit any intuitive idea of what a closed shape is. So, we have a choice. Either we keep the requirement for a closed shape and conclude there are no circles or we reject that part of the definition, which results in the existence of a number of circles.

In Euclidean geometry, that the circle is closed is redundant in the definition. If you define a circle as a set of points at a particular distance from a center, a consequence is that it is closed. I did not bring this up in the discussion, but it is an interesting point to keep in mind for the future.



The discussion above on whether the center is a part of the circle is worth thinking about. A circle in Euclidean geometry is completely determines in terms of location and size by specifying a center and a radius. Both of these are parameters we can vary to create new circles. Thinking about the center as a part of a circle is similar to thinking of the radius as part of a circle. The difference between an object and what is specified in order to define an object may not be apparent to the students. I have not been able to find any literature on this.

**Initial exploration of triangles – Indus.**

After a short discussion on circles in Indus, we moved on to triangles. I asked:

*Me: Are there triangles in this world?*

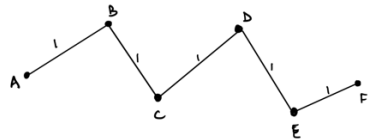

*Figure 53 Triangles in a simple world*

*In this world. Six points.*

*Multiple students: No.*

*Me: In this exact world.*

*Multiple students: No.*

*Me: What is a triangle?*

It is interesting that students jumped to the conclusion that there is no triangle in this world even though they just went through an exercise where they showed that depending on your definition of circle and bisection, there may or may not be circles or bisection in the world.



So, I asked them to define triangle to which Tarini responded:

*Tarini: It is a closed figure. With three sides.*

*Tanya: Three vertices.*

*Me: Three vertices and straight lines connecting those vertices.*

*(I wrote down a definition)*

*Is this good? A closed figure with three vertices and every pair of vertices is connected by straight lines. Okay. So are there any triangles?*

*Multiple students: No.*

*Uday: Unless we connect any other points.*

Given the definition of a triangle I wrote down, three collinear points would count as a triangle even in Euclidean Geometry unless the notion of closed figure requires there to be a non-empty inside. So, by this definition, there are triangles in this world, namely any three collinear points since each pair of them is connected to the others by a straight line. The only ambiguity is the meaning of 'closed.' The students did not come up with any triangles, so I gave a suggestion:

*Me: Let me suggest one. Let me take the vertices A, B, C. The straight lines AB, BC and ABC. So, this is three vertices each connected by three straight lines. Are these three straight lines?*

*Multiple students: Yes.*



*Me: So, all these three are vertices. All these three are straight lines. So, is this a triangle?*

*Multiple students: No.*

*Tarini: No. It's ABC, not AC.*

*Devika: Oh. I didn't see that.*

*Tarini: No, it's not a triangle. Because ABC… You said AB is a line segment (and then says something which I was unable to parse)*

*Me: AB is a line segment. BC is a line segment.*

*Tarini: BC is a line segment.*

*Me: ABC is also a line segment.*

*Tarini: But, that's not a triangle.*

*Me: Why?*

*Tarini: There is no AC.*

*Me: There is AC. ABC is AC. It's the shortest path.*

*Uday: It only has one angle.*

*Me: I don't know what angle means here, right? What is an angle here?*

*Tanya: It is not a closed figure.*

*Me: What do you mean by closed figure? Why is this not a closed figure?*



There are three objections to ABC being a triangle. The first is the objection that ABC is not a line segment from Tarini. The second is Uday's objection about triangles and the third is Tanya's objection about the shape not being closed. There is an interesting difference here between Uday's objection and the other two objections. Uday's objection is unrelated to the definition we had agreed to. The definition says nothing about angles and hence the number of angles or whether there are angles at all doesn't impact whether the shape is a triangle. In fact, it is unclear what you would even mean by angle in this world. On the other hand, the other two objections are rooted in the definition. If ABC turned out not to be a straight line, then the shape in question would not be a triangle. The same would be true if the shape was not closed.

I continued:

---

*Me: Why is this not a closed figure. It contains three things*

*Devika: They should be adjacent.*

*Me: They are. This is adjacent to this. This is adjacent to this.*

*Uday: (Says something which is not audible) So, AC is not connected.*

*Me: It is. We agreed that ABC was a straight line, right?*

*Prerna: AC is not connected.*

*Me: But, ABC is a straight line because it is the shortest path from A to C. Straight line means shortest path.*

*Multiple students: Oh*



*Me: A and B is a straight line. B and C is a straight line, shortest path. ABC is a straight line, shortest path. So, all these three are straight line, these are vertices. I'm not sure what closed figure means over here. But, maybe it is a closed figure, maybe it's not. That is the only point of contention.*

*Me: But, notice that this problem exists in Euclidean Geometry as well. Let me draw this here. Let this be a straight line. So, in Euclidean Geometry, is this a triangle?*

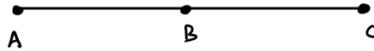

*Figure 54 Degenerate Euclidean triangle*

*Multiple students: No.*

*Me: Why? Because it has three straight lines. It has three straight lines?*

*Uday: But, there are two vertices.*

*Me: There are three vertices (pointing at them in the diagram). There are three sides – AB, BC and AC.*

*Tanya: There shouldn't be a vertex in between.*

*Me: So, you don't want collinearity. You don't want this to be allowed.*

*Uday: Two lines should not overlap.*

*Me: So, you don't want them to be collinear.*



The analogy to Euclidean Geometry is useful here since the definition we had allows for three collinear points to be vertices of a triangle. We have to explicitly add in something which disallows that if we want to do so.

This seems to show that the students' objections were not to things like whether something is a straight line or whether the shape is closed. Rather, it seems to be to the definition itself. At least the students who participated in the discussion seem to want a definition which disallows such degenerate triangles.

The session ended by separating out the two types of triangles. I called triangles which include the degenerate case Alpha-Triangles and I called triangles which have the stipulation of non-collinearity Beta-Triangles. We concluded the session by saying that there exist alpha triangles in this world, but not beta triangles.

### *Summarizing the initial exploration of triangles.*

Exploring triangles was similar to the exploration of circles. Once again, what we were trying to do was to extend a definition from Euclidean geometry to these discrete worlds. One difference was that the question of degeneracy came up here (without mentioning the term). If we allow for three collinear points to be the vertices of a triangle, then what looks like a straight line is a triangle.

Taking this further, what looks like a triangle could be a quadrilateral or a pentagon or any polygon as long as we specify points on that shape as vertices. This was not a discussion which happened in the session but could be something worth discussing here and in the Triangle Theory Building module in future editions of the course.



**Straight lines and line segments – Indus.**

In this part of the module, the focus was on lines as distinct from line segments in the worlds we have created. The goal is to figure out whether there is a coherent concept of lines which is somewhat similar to the concept in Euclidean Geometry.

*Me: In Euclidean Geometry, you have these two concepts: line and line segment. You have these two things. What is the difference in Euclidean Geometry between a line and a line segment?*

*Uday: Line is infinite and line segment is finite*

*Tarini: It is part of a line.*

*Me: Yeah. This is infinite, as you said. And line segment is part of a line.*

*Tanya: It has infinite points.*

*Me: But, it is not infinite in length. It is between two points. In this world, let us look at the bead world. No let's look at this. What would you like lines to be? This is an example of a line segment. Just this.*

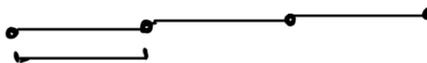

*Figure 55 Lines and line segments 1*

*You pick these two. It is a line segment, right?*



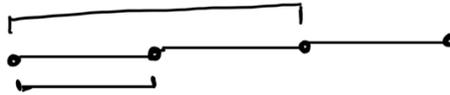

*Figure 56 Lines and line segments 2*

*What would you want lines to be? Lines are bigger than line segments. They are kind of like the biggest things which line segments can make.*

*Devika: Line in the six point world?*

*Me: Let's assume it is a four point world. Let's leave it as a four point world for now. Just to make things easier. So, one way of characterizing line is they are infinite.*

*Tarini: But in this world, we can't have them.*

*Me: Then we can't have them.*

There are at least two ways in which we can borrow the concept of a line from Euclidean Geometry. One is to borrow the idea of an infinite length straight object as suggested above. Given that notion, as Tarini said, we cannot have straight lines in this world. However, there is another way to extend the concept of line from Euclidean Geometry which I discussed straight after this:

*Me: But another way of characterizing lines is as an extension – the largest extension of line segments.*

*Tanya: Then here we put the first to the fourth points.*



*Me: So, notice that these are two ways of characterizing lines. At least conceptually we can understand these two, right? Think of it as something extended. The other way is to think of it as something infinite. If we take this one, that is meaningless in our world. Because there are only a finite number of points. So, if we want the concept – now this is a choice we have to make – if we want the concept of line, then this might be usable, which in Euclidean geometry will result in an infinite line, but here it doesn't. What's an example of a line here?*

*Tanya: The left most point to the right most point.*

*Me: Okay. So, what's the characteristic of a line in relation to line segments?*

*Tarini: Combination of line segments make a line.*

*Me: Yeah. Can you say it the other way around? Lines from line segments.*

*Tanya: Parts of lines.*

*Me: So, if I pick two points on a line in Euclidean Geometry, I get a line segment. So, that should happen here also. If I pick two points on a line, do I get a line segment? Is the line segment on the same line?*

*Tanya: Yes.*

*Me: So, if we take the entire thing to be a line, then it's on. So, you can't take just this to be a line, right (pointing at just the first three points)? It is not the largest extension.*



By conceptualizing a line as the largest possible extension of a line segment in some sense, we do have lines in this world. In fact, we have exactly one line which contains all the four points. By this conceptualization, all lines in this world, and in finite point worlds generally, are line segments, but all line segments are not lines.

In the other excerpt, I was quite unclear when I asked, 'Can you say it the other way around? Lines from line segments.' What I meant to say there was: what can you say about points on a line and the line segments between them?

I continued by making the situation more complex by adding in another point as in the figure below:

*Me: And, what about if we had this world?*

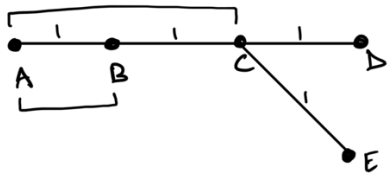

*Figure 57 Straight lines and line segments 3*

*Tanya: There are two lines.*

*Me: Two lines, right? Are there two lines or more than two lines?*

*Tanya: Two lines.*

*Me: So, what are the two lines?*

*Tarini: AD and AE.*

*Me: ABCE and ABCD. Those are two lines. What about… (I gesture at DCE)*



*Tarini: That's not a…*

*Uday: The value of CE is greater than the value of CD.*

*Tanya: They are all equal.*

*Uday: It is not given there.*

*Tanya: but we assume that.*

*Uday: The value's not given.*

*Me: So, let's assume they are all equal for now. Just for now.*

We had discussed giving different values for lengths of different connections before this. I didn't explicitly say that I wanted them to be equal for this example. However, even if they were of different lengths, it is hard to figure out what Uday is pointing out. Even if the length of CE and CD were different, that wouldn't change whether DCE was a straight line or not.

*Tarini: DCE, the value of the distance is less. So, it is not a line.*

*Me: Okay. So, if we say this. The largest possible extension of a line segment. By largest we mean in length, right? If we take this line segment AB, what is the largest possible extension of this?*

*Tanya: AD*

*Me: Why can't all the points be on the line? What about a line stops all points from being on the line? You don't want all points to be on the line, right? Why?*

*Tanya: I didn't get that.*



*Me: So, here, before we added the E, the line consisted of all points. I've added point E. Now, you are saying ABCDE all together do not form a line.*

*Tanya: Yes.*

*Me: Why?*

*Tanya: Because if we take ABCD CE, they would all be overlapping.*

*Tarini: It wouldn't be shortest path.*

*Tanya: I'm saying ABCD CE.*

*Me: So you are saying a line is not just some points, but a list of points. So we said ABCD is a line. But, if you add E, then notice that this whole thing is not a line segment.*

*Multiple students: Yes.*

*Me: Even this BCDE is not a line segment. Because you can't go from D to E. In fact DE is not a line segment.*

*Tanya: We can put DCE.*

*Me: So, this (writing ABCDCE) not every subset is a line segment. But every subset is a path. Earlier, not everything was a path because DE was not a path. Here, now, CDC is a path. It is a reasonable path. But, is it a line segment? What are the end points of CDC?*

*Multiple students: C.*



*Me: Is CDC the shortest path between C and C?*

*Multiple students: No.*

*Me: No, right? Stay put. Don't move. So, this is not a line segment. So, it violates that the subsequences should be line segments.*

I used the terms 'subset' and 'subsequence' interchangeably in the exchange above. I should have stuck to one of them, ideally subsequence, and clarified the concept with the students. The question I asked midway through the exchange was whether all the points constitute a line. Tanya's response is interesting. She seems to have understood lines as a type of path or a sequence of connected points. A path which connects all the points is ABCDCE and that is not a straight line. While that isn't sufficient to prove that there is no straight line which connects all the points, it does contain the intuition which the proof would be based on.

We then returned to discussing DCE:

*Me: DCE. Why is that not a line?*

*Tanya: DCE is a line.*

*Uday: It is not a line.*

*Tarini: DCE is a line. We are not saying the line – the big line. Line segment.*

*Me: It's a line segment. But is it a line?*

*Tarini: No. It's not.*

*Me: What about it?*



*Tarini: We have the general perception that line segments add up to make a line. But here if you take a line, we know that the line is three distance. But, if we take DCE, then it's only two distance. So, it's not the largest…*

*Me: Possible extension.*

*Tarini: Possible extension.*

*Me: But, is it not? What about this line segment – DCE. That's a line segment.*

*Tarini: Segment yes. It is a line segment.*

*Me: DCE is a line segment. It is the shortest path from D to E. What's the largest possible extension of DCE which is a line.*

*Uday: DCE.*

*Me: Then DCE is a line.*

*Tanya: But, only if we consider D and E as the end points.*

*Me: You said that this is the largest possible extension of a line segment which is also a line segment.*

*Tanya: (Something inaudible)*

*Me: But, it is the largest possible extension of this line segment.*

*Multiple students: Yes.*

The ideas here are a little counter intuitive since, for example, the line segment CD would be part of two lines if we take the definition of line we came up with seriously – ABCD and



DCE. This is not true in Euclidean geometry where a given line segment is part of exactly one line. While DCE is not the 'longest straight extension' of CD, it is the longest straight extension of DCE.

This seems to be the exact issue which Tarini brought up when she said, '… we know it is three distance. But, if we take DCE, then it's only two distance' and which Tanya brought up when she said, 'But, only if we consider D and E as the end points.' I don't think I really dealt with this. I could have worked through a few more examples of a few more similar worlds.

After this, I had a short discussion on straight lines in 'necklace' worlds like this one:

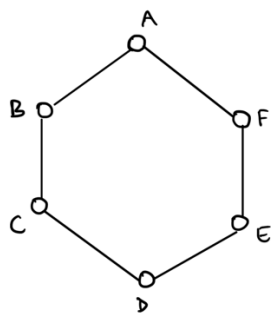

*Figure 58 Necklace world*

That ended the day and the module in Indus.

***Summary of straight lines and segments.***

In this section, we returned back to a concept we started the module with. The relation between line segments and lines in Euclidean geometry is something students are familiar with. Given a line segment, it is part of exactly one line and if we cut a line at any two points, we get a line segment. We already had a concept of line segment by this point.

The question to ask is whether we can extend the concept of line from Euclidean geometry such that all line segments are parts of lines uniquely, given the concept of line segment we already had. For general discrete worlds of the type we are dealing with, that doesn't



seem to be possible. We have to lose at least some of the properties of a line. The direction we took during the class was to allow for line segments to be a part of more than one line by saying that a line is the longest extension of a line segment which is also a line segment.

**Students' reflections on the module**

Students wrote reflections on the first and last days of the module. They were asked to answer the following two questions:

3. What did you learn?

4. What did you not understand?

Apart from this, students were asked for general reflections on the entire course at the end. In the reflections, they were asked the following four questions (the template for the survey is available in the appendices):

1. What was your favorite part of the course? Why?

2. What was your least favorite part of the course? Why?

3. What did you find valuable in this course? Why?

4. What are one to three specific things about the course which can be improved?

The goal of this subsection is to give a sense of what they said related to this module. As mentioned in the methods section, I coded sentences in students' daily reflection responses on the basis of whether the responses were related to specific mathematical content or to the types of general theory building abilities I discuss in Chapter 3. Some sentences could have both codes. The sentences without codes turned out to be things like working in a team or to do with the instruction. So, I added in two additional codes which were 'non-mathematical learning' and 'reflections on instruction.'



The same set of codes were applied to the end of course reflections along with an additional set of codes to do with whether a particular sentence applied to a particular module. Here, I will be mainly concerned with those instances related to discrete geometry.

**End of course survey.**

Like in the case of triangle theory building, there was a significant difference between how the students in the two schools received this module. As mentioned in the Triangle theory building chapter, it is not clear as to what the explanation of the differences between the two schools could be. However, I will discuss this further in the discussion chapter.

In the case of discrete geometry, students in Ganga had positive reactions and with only one student listing it amongst their least favorite parts of the course while eight out of 15 students listed it amongst their favorite parts. In Indus, seven out of 16 students listed it amongst their least favorite parts of the course while only two students listed discrete geometry amongst their favorite parts, though there were a few who said they enjoyed everything.

The positive reactions included:

-   'I also liked the discrete geometry session because it was something which was very different from what we usually learn in geometry. It made me think outside the box. When we used to define what lines, points, etc. are I used to feel kinda irritated because we never got it right but it was really fun.' – Devika, Indus

-   'My favourite part of the course was 'discrete geometry'. Because it was very different in terms of Euclidean geometry and challenging for the mind in terms of 'drawing circles'.' – Gauri, Ganga



- 'My favourite part of the course was discrete geometry. Especially the questions relating to worlds with only a given number of points. It encouraged us to find alternatives and loopholes and also use our creativity.' – Anjali, Ganga

The negative reactions included:

- 'The least favorite part of mine in the course was defining circles in discrete geometry. I did not like it as I was not able to understand the different types of circles in different worlds.' – Arun, Indus

- 'The least favorite part of the course was the discrete geometry as I cannot help in my group as for giving answers and I least knew about discrete geometry.' – Aryan, Indus

- Sometimes I was unable to find interest in discrete geometry sessions. that was so because we were digging a lot into stuff which made me feel that it was complicated.' – Arvind, Indus

The negative reactions to discrete geometry are all to do with not understanding aspects of it or finding it too hard. Those who had a positive reaction seem to have had that reaction for a similar reason – that it was challenging.

Aryan's reaction to the previous module was that he enjoyed it because he was able to contribute to discussions due to his prior knowledge of triangles. A module like this one is much further from what students are used to than the triangle theory building module. Hence, it isn't very surprising that it was harder for students to engage with. This module also builds up upon itself more than the previous one – you start by creating objects and then use those objects later on while constantly adding more objects to your repertoire.



Amongst the positive reactions, Devika's is worth highlighting since she mentions that she was irritated at times during the module but yet found it fun. There were three other reactions to this module which explicitly said that they found it hard at first but then enjoyed it.

**End of day reflections.**

It is interesting to compare the general reactions to the course to the end of day reactions on the first and last days of the discrete geometry module implementation in the two schools. Aryan and Varun were the only two students in Indus who answered the second question about what they did not understand. Both of them mentioned the discussion on triangles.

In response to the first question in Indus on the first day, most of the students discussed details about discrete geometry that they learnt such as 'different ways of bisection and 'a new perspective on triangles.' A few other students did mention very general things like 'logical reasoning' and 'improving our thinking.' The only theme seemed to be about interpretation. Three students explicitly mentioned that they learnt about interpretations of questions/statements. For instance, Aditi mentioned that, 'We learned how to answer questions with different interpretations and come to conclusions considering the same.'

In Ganga, six of 12 students mentioned interpretation or something similar. For instance:

- Aditya said, 'I learnt to interpret the question before answering and understand what does it mean.'

- Sana said, '=> Interpret questions => Form definition => Answer changes with change in definition => looking for possibilities'

- Anjali said, '- Understanding questions - Finding loopholes in vague questions - Proving by contradiction - Thinking of possibilities and trying to prove them (looking for alternative possibilities) - Changing answers after change in condition'



On the last day in Ganga, students mentioned a large number of things and there is not really a theme which runs across the responses. However, Vivek's response summarizes what other students said. He says, 'Discrete geometry is a very vast subject with many contradiction. If you make 1 contradiction, you derive another question. So it is extremely exciting. In this session there was 1 question which ate up the brain of all. But to questions there are conclusion. Thus we came up to one. The question which were derived were of similar kind but had different solution which was an interesting part. In today's session was a continuation but with more tougher level which made us think a lot and that is the objective one should have to think. The solution of every student was different but had some connection. The different worlds with different points having various properties of figure this is what we had to guess and was an interesting part.'

On the last day in Indus, all students barring one talked about things they learnt about discrete geometry. For instance, 'about the types of circles', 'prove and contradict the possibility of circles,' and so on. Aryan was the only student who did not discuss details about discrete geometry. He said, 'I learnt today that we don't need bookish knowledge but we have to think out of the box.'

**General observations and reflections on the module**

**Working with definitions in this module.**

In the triangle theory building module, students were using definitions and maybe clarifying definitions they were already familiar with. In the definition guessing module, students were attempting to uncover definitions in somebody else's mind. In this module, students had to create definitions in order to imbue the initial question with meaning as well as to import other objects from Euclidean Geometry.



There were five major examples of defining objects in the excerpts above – line segments, bisection, circles, triangles and lines. One strategy which I used throughout was to allow for multiple definitions for each term. For instance, bisection was defined in two ways – 'bisection at a point' and 'bisection at a connection.' For triangles, we had the type which included triangles whose vertices are collinear and the type which didn't. This allows students to see the relationship between the definitions and their consequences. Given one definition, there may not be triangles in a given world while, given another definition, there may be.

It might be the case that this is the part which students found the hardest since they are used to being given definitions rather than creating them and playing with them.

**Productive Struggle.**

From the feedback at the end of the course, it seems evident that students found the discrete geometry module challenging. Those who had a positive impression of it said they found it challenging and even irritating at times. Those with a negative impression seem to have found it too hard. This is almost the opposite of what happened with triangle theory building where those who had a negative reaction said it was 'boring' or 'easy'.

This is not very surprising since this session was the furthest away from anything they had learnt in their mathematics classes. The journey from struggle and even frustration to learning is exactly what productive struggle means (Heibert & Grouws, 2007). It seems like some of the students did go through that journey and get value from the module. However, it also seems to be the case that for some students the struggle did not result in learning, at least not as much.



One thing which differentiates this module from the triangle theory building module is that it builds up on itself a lot more. If a student did not understand something in the first session, that will impact their understanding later on.

Hence, for a module like this, there does need to be more scaffolding available to those who need it so that they are not left behind. There also needs to be some way of evaluating what students have understood, especially at the end of the first session.

**Following students' direction.**

Early on in this module in Ganga, I had a relatively long dialogue with Vivaan which touched on a lot of different things. However, the question I had initially asked was not really addressed and neither were many of the things Vivaan brought up. The reason for this is that Vivaan brought up many different things one after the other and I followed along with him.

I do think there is room for one-to-one dialogues in a class since others can learn from listening to them. However, this dialogue went on for a long time and I don't think it was very fruitful. It also strayed the general conversation from what we were discussing. So, while I think it is important to go down paths individual students are interested in, that needs to be done keeping in mind the learning outcomes being aimed at and the learning of the other students.

**Representations in this module.**

The initial representation of the 6-point world is something I gave the students in both the schools. I have so far been unable to get students to construct a representation themselves including when I have done this session with other groups. This is not surprising since the initial question is vague and most people do not have experience dealing with such questions in mathematics.



We ended up using that same representation for all the worlds we created in this module – that of dots representing points and edges representing connections. I could have introduced or gotten students to create alternative representations. For example, representing the worlds as an adjacency matrix. This may have helped initially when students were struggling with the idea that there weren't necessarily points 'between' two points.

**Justification in this module.**

I took a conscious decision at the beginning of this module to not push too hard for rigorous justification. The reason is that the triangle theory building module was already doing that and I wanted to focus attention on a slightly different set of learning outcomes. That isn't to say that justification was not a part of this module. In fact, justification was a part of everything but just not at the level of rigor as it was in the triangle theory building module.

For example, we never actually wrote out explicit axioms of the world and instead relied on diagrams and intuitions. If this module were to continue for a few more sessions, this is the direction I would have taken it in. I would have insisted on more rigor and explicitness.

There is another strategy which can be used in a module with similar content. However, this strategy would change the learning outcomes. Instead of starting with a question and some pictures, you could start with a set of constraints (axioms) and get students to create the representations and worlds from the given constraints.



### Chapter 10 – Tests and feedback

In this chapter, I will present the results of the student feedback and tests. I have already discussed the end of day reflections and the post course survey in relation to particular modules. Here, I will focus on those aspects of the feedback which are more general in the post-course survey.

With regards to the tests, I will be presenting the changes in student scores on the same question types from the beginning to the end of the course. I am not looking to arrive at any conclusions from the test due to a small sample and that this test has not been previously used

**Student Feedback**

I had asked students to provide feedback at the end of the course. I asked them the following questions:

1. What was your favorite part of the course? Why?

2. What was your least favorite part of the course? Why?

3. What did you find valuable in this course? Why?

4. What are one to three specific things about the course which can be improved?

The first two questions largely focused on the modules and I have already discussed them in the previous chapters. In this chapter, I will focus in the last two questions. Since each question had just thirty responses, my method here was to look through all the responses and look for themes which emerged.

**What was valuable.**

There were four themes I found across student responses to the third question on the survey. Students mentioned that the course made them think, the value of intellectual skepticism, understanding mathematics, and the value of teamwork.



### Making us think.

All of the comments students made mentioned the course making them think. For example, Imran from Ganga said, "The most valuable lesson I find from this course was to never stop thinking and a question may have many answers." Adding to that, Bharat mentioned thinking about a problem from different points of view. He said, "This course increases your ability to conquer a problem by thinking about it from different points of view. It also developed the thinking capacity for tackling a problem." Zoya from Indus brought up curiosity. She said, "I like how they taught us to be curious & try to comprehend 'what, how & why' for everything." Tarini from Indus mentioned curiosity as well and brought up thinking out of the box. She said, "It helped me think out of the box. It will help me to be curious."

### Intellectual Skepticism

Many of the students mentioned the idea of being skeptical. Tanya from Indus said, "What I found valuable is that we learnt how we shouldn't stick with textbook knowledge & we should prove or justify the things we know to understand it deeper." Zoya from Indus said, "Earlier my mind used to believe everything that came in front of me. But now I ask questions like why. I try to explore mathematics try finding the meaning of words which we usually take for granted." Prerna from Indus said, "The course taught us not to assume the given facts to be true. We should cross check and go to the root of things to understand them better. It made us ponder a lot and question things."

### Understanding Mathematics

Students mentioned understanding of mathematics in two ways. One way was in terms of understanding mathematical concepts. The other was about understanding mathematical practice. Examples of the former include Vandana from Indus who said, "Because of this workshop I was



able to understand the real depth of mathematical concepts and its real meaning," and Sania from Ganga who said, "The proof and characteristics of a triangle."

The latter was more commonly mentioned. Tanya from Indus said, "We understood how mathematicians build up chains of derivations & conclusions from basic knowledge or axioms." Pankaj from Ganga said, "I actually came to know that no one proof/theorem is enough - we need more and more." Prerna said, "I try to explore mathematics try finding the meaning of words which we usually take for granted."

### *Teamwork*

The ability to work in teams was not one of the learning outcomes the course aimed at. However, many of the students mentioned it as a valuable learning for them. Pankaj from Indus said, "I learnt to function well in groups." Aditi from Indus said, "the group discussions, logical reasoning and debates were all very engaging and I learnt a lot from these." Aryan brought up the various points of view in the groups. He said, "When we (my group) together solves a problem we look at the thing from different angles which is really valuable." Sanya from Indus mentioned that they were previously not used to working in group. She said, "As we were all not used to group discussions earlier regarding the topics Madhav had given, we were excited and amongst ourselves were arguing for small reasons. So, it was a determined course where we all would speak & not textbook knowledge but share different views and opinions about everything."

### Changes to the course.

Five students left this question blank or said that there were no improvements. A few students mentioned changing the length of the course and of parts of the course. Some students



mentioned adding more topics and others mentioned adding assessments during the course. There were a few other individual suggestions I have included at the end.

### *Timings constrained.*

Multiple students, all from Ganga, mentioned that they wanted the course or parts of the course to be lengthened and some mentioned that they wanted the sessions to be longer. Deepa said, "Some points were left undiscussed because of the time limit." Sania said, "I think that the duration of the course must be increased a little so that more time can be spent for discussing." Anjali said, "More time given to a particular topic." Karan said, "Increased time (upto 3 hours). Increased frequency (upto 4-5 times a week)." However, in the other direction, one student, Anuj, said, "Some of the discussions can have and endpoint and not go on for hours."

Thinking carefully and rigorously for long stretches of time can be tiring. However, very short sessions could result in lack of depth. It isn't clear what the ideal length would be for a course of this type, and it is probably the case that stamina for a course of this type will vary across students.

### *More topics.*

Quite a few students mentioned that the course include more topics. Gauri from Ganga said, "The topics we learnt, maybe some more can be added." Aditya from Ganga said, "I think you should also include parts like altitude, perpendicular bisector, angle bisector and median concurrency." Aditi from Indus said, "I feel if some algebra was discussed too, then it would be more knowledgeable."

The focus of the course in terms of content was geometry. However, there are theory building activities which can be done in areas like number theory. In future iterations of the course, these could be valuable additions. In the case of Gauri and Aditya, this reaction could be



interpreted as wanting more variety of topics to think about. However, Aditi, who wanted to be 'more knowledgeable', may have missed the purpose of the course. The purpose was not that they learn geometry or mathematical content. It was that they learn to do mathematics, specifically construct theories.

### More assessment.

A few students mentioned the need for more assessment. Aditya from Ganga said, "I also think that you should increase the number of questions asked in both assessments (*talking about the pre and post-test*)." Pankaj from Ganga said, "I think that after every lecture there should be a worksheet to brainwash the things and your worksheets are one which need thinking." Aditi from Indus said, "And tests could've been more better, as in more questions should have been added."

The students who mentioned tests seemed to have found the questions enjoyable. These were low stakes tests – in fact, I did not even give them scores on the test. It would be interesting to see what their reaction to the tests would have been if these tests determined a grade on the course.

### Other suggestions.

Aditya from Ganga was the only student who mentioned something about me. He said, "I think he should be a bit louder because sometimes he speaks in a low tone." Imran wanted, "more activities using arts and crafts."

Two students from Indus mentioned concerns related to their groups. Tushar said, "When a student does not understand something that is being discussed he is completely left out of the discussion and it makes it hard to understand things related to it later. Something has to be done



about this." Udit said, "There should be neutral distribution of team like one should not have all intelligent or one not participating amongst them."

There are interesting concerns here related to group distribution. In a group, if a few students disengage because they find something too difficult, that could hurt the learning of others. On the other hand, if there are students who are really excited by what they are learning, they could inspire others in their group to try harder to understand what is going on. During the workshop, I did not pay attention to concerns about group dynamics. I allowed for students to create their own groups.

**Pre and post-tests**

The pre and post-tests used in this course each consisted of two questions on theory building. Each of the questions had multiple parts. The first question was to do with definitions. Like the podgon module, it gave students a term and a list of objects, some of which were examples of that terms and others which were not. Given this information and a set of definitions, students had to figure out whether those definitions were possible definitions of that term and had to pick which of the objects served as a counter-example if they were not. There were six definitions given. In terms of grading, they got a point for a correct answer as long as it was accompanied by the right counter-example. They got no points if they stated a wrong counter-example.

The second question was to do with reasoning within a system. Students were given a list of assumptions and a list of conclusions. They had to say whether those conclusions were true, false, or not determined, given the assumptions. The statements in the pre and post tests were the same form with changes to the language used. There were seven assumptions given and seven conclusions. You can find the pre and post-tests in Appendix B.



On the first question, out of six, the average score on the pre-test in Ganga was 3.3 while the average score on the post-test was 4.7. Of the eleven students who took both the pre and post-tests, eight saw an increase in their scores and nobody saw a decrease. In Indus, the average score in the pre-test was 4.7 and the average score in the post-test was 5.6. Six students saw an increase in their scores, one saw a decrease by a single point, and the rest of the nine stayed the same.

The second question had the scores reversed in terms of the two schools. At this point, it is unclear what to make of this difference. The average score in Ganga went up from 3 to 3.6 while the average score in Indus went up from 2.8 to 3.1, both out of seven. In Ganga, seven students saw an increase in their score, one saw a decrease and the remaining three were on the same score. In Indus, six students saw an increase in their scores, two saw a decrease (by a single point) and the remaining eight had the same score. Only three students across the two schools saw a change of more than a single point and no student saw a change of over two points.

There is not much we can draw from this data since the instrument itself had not been validated and the sample size is relatively small. However, there is a clear increase in the scores from the pre-test to the post-test, especially in the first question. The questions in both the tests are relatively similar. However, it could be the case that the seeming minor difference in the tests resulted in the change in the scores. Even so, given the change, it might be useful to explore the use of this instrument as an assessment tool in future research.



## Chapter 11 – Discussion

The goal of the last few chapters of results was to give the reader a sense of what a course on theory building looks like inside the classroom – how students engaged with the course and the role played by the facilitator. General conclusions about implementing theory building courses, or even this theory building course, are not warranted given the narrowness of the implementation of the course – it was implemented once in two schools, part of the same group of schools, with a single facilitator who also designed the course. Also, many of the students have had experience working with me or my associates in the past. It is impossible to separate their relationship with me and preconceptions of me from their engagement with the course materials.

However, there is some learning from the implementation we get from looking carefully at specific instances, and there is some plausible speculation we can make, extrapolating from this experience. This chapter is divided into three parts. The first part involves putting together various general observations and reflections from the various results chapters related to the concepts of theory building. The second part is on other learnings from the course implementation, with a focus on the considerations we ought to keep in mind while iterating on this course or creating other theory building courses. The third part is on connections between this research and other areas of research and education – in this section, I will discuss mathematical practice, transdisciplinarity, equity, developing mindsets, and various curricula from around the world.



**General observations and reflections on concepts of theory building**

In this section, I will pull together various aspects of the general observations and reflections from the various results chapters related to the concepts of theory building discussed in chapter 3.

**Proof.**

Different modules in the course varied on the dimension of rigor in reasoning expected from the students. The triangle theory building module was the most rigorous since a large percentage of the time was spent on proving. The podgon module was not very rigorous due to it being an icebreaker. However, it can easily be made into a more rigorous exercise as I discussed in the chapter. The discrete geometry module focused on the creative aspects of theory building. Hence, rigor in reasoning was not given as much importance as in the triangle theory building module.

***Proof schemes.***

Harel and Sowder (1997) introduced the concept of proof schemes. Using these, they try to capture the various ways people justify claims and what types of justification they find acceptable. At the broadest level, these proof schemes are of three types: deductive, empirical, and authoritative. I have unpacked these further in the literature review.

I mentioned in the triangle theory building chapter that students seem to be using a deductive proof scheme. I was unable to find a single instance in that module when that was obviously not the case. As I mentioned in the chapter, this is surprising given the literature is rife with students using empirical and authoritative proof schemes. This could be a matter of selection bias given that students chose to join this course, the education system in India being



substantially different from the contexts discussed in the literature, or it could be that students learnt from previous workshops I had done at the schools.

On the other hand, in the podgon module, there were some examples which look like empirical justification. For instance, from a single example of a non-closed shape which was not a podgon, a student concluded that podgons must be closed.

However, the reason this is more excusable is that they were not attempting to prove a claim within an existing theory in this module. Rather, they were trying to guess the definition in the heads of the instructor(s). Given that, it is impossible to use a deductive proof to come to conclusions in this module. Even so, this at least indicates poor communication of reasoning. I will discuss this further in the next subsection.

### *Flawed reasoning and communication.*

There were a few examples of what appeared to be flawed reasoning in the episodes I presented in the previous few chapters. Four types of flawed reasoning I highlighted were:

1. Applying a theorem outside of its scope

2. Deducing P => Q from Q => P

3. Cyclic reasoning

4. Coming to general conclusions from examples

An example of the first was when a student applied a theorem about two sides of a single triangle to sides of two different triangles. It is clearly not the case that a theorem depending on the relative lengths of the sides of the same triangle would apply to the relative lengths of the sides of two different triangles unless some other constraints were met. An example of the second is when a group said that the claim they had come up with was the same as a claim the other group had made. However, the claim made by the other group was the converse of their



claim. An example of the third type came up when discussing straight lines in both the triangle theory building and discrete geometry modules. Students defined straight lines in terms of collinear points and collinear points as points on a straight line. When they spotted this, they immediately realized that there was something wrong. However, they did not know how to fix the problem. Examples of the fourth type were mentioned in the previous subsection and were largely limited to the Podgon module.

All these cases could be explained by how students reasoned, communicated their reasoning, or both. For instance, it could be that students wanted to say something weaker when generalizing from definitions in the podgon module but did not have the language to communicate that. Looking back at the course from the point of view of the instructor, my guess is that it is a mixture of both, with each contributing to the other. In fact, it was my intention to run a course over the summer with the same group of students on clarity and precision in language. I had designed the outline of the course but was unable to implement it due to COVID 19. I will discuss the relationship between language and mathematics later in this chapter.

**Conjecturing.**

Conjecturing did not play a significant role in the triangle theory building module. However, it was an important part of the discrete geometry module and the podgon module.

While the claims in the podgon module were not truth-claims, the thinking involved did bear resemblance to conjecturing and judging whether conjectures are plausible. There were two parts to the activity in the module. One was students guessing definitions and the other was students looking for examples to test the definitions they had guessed. The guessing of definitions is very similar to conjecturing – they had to guess a definition which would fit all the examples they had access to.



Examples to test the conjectures were of two types. The first were counter-examples which allowed the students to rejects a definition or aspects of a definition. The other type of example, which turned out to match the definition, gave empirical support to the definition it was testing. This empirical support makes the definition more plausible.

The conjectures in the discrete geometry module were interesting for another reason – they were made at various levels of abstraction. There were conjectures at the level of particular worlds. For instance, a world with six points arranged in a specific manner. However, there were also conjectures at the level of all worlds with a particular shape. For instance, about all worlds such that each point has exactly two neighbors. Finally, there were conjectures about all worlds of the type we were dealing with. For example, about the impossibility of all lines being bisect-able by a one definition of bisection.

**Definitions.**

Defining was involved in all three modules. The podgon module was about defining, equivalent definitions were a part of triangle theory building, and extending definitions played an important role in discrete geometry.

The idea of multiple definitions for a given word was an important part of all the modules. The podgon module had two different 'definers' with different meanings for the same word. In the triangle theory building module, the idea of equivalent definitions came up, most notably when discussing equilateral triangles. Students showed that equilateral triangles are equiangular and equiangular triangles are equilateral. This resulted in being able to define equilateral triangles either in terms of angles or in terms of sides.

The discrete geometry module involved multiple definitions for the same word, each with different consequences. It started with the definition of bisection and included the definition of



triangles and circles. Each of these had two definitions. However, the goal here was not to choose between them but to pursue the consequences of both definitions independently.

An important aspect of the discrete geometry module was to extend definitions from Euclidean geometry. This was how the worlds themselves were created. We borrowed the concept of straight lines. To answer the initial question, the concept of bisection was also borrowed. Later in the module, triangles and circles were extended from Euclidean geometry to the worlds we had created. The strategy for extending definitions was to look for a concept from Euclidean geometry which best extends. For example, straight lines in Euclidean geometry have various properties associated with them which do not even make sense in the worlds we had created. We did not have an obvious concept of direction. Hence, ideas like angles are probably not useful in these worlds. So, the direction/angle property of straight lines in Euclidean geometry is not something which can be used in the worlds we had created. However, the property of straight lines related to distance is something which can be extended and is what we chose. Another example of this was not stipulating that a circle needs to be closed since if we made circles closed shapes, they would be exceedingly rare.

**Classification.**

Classifying and defining are closely related. There were two ways in which classification came up in the course. One was the discussion of classification of triangles in the triangle theory building module. Another was the classification which happened as a results of proving certain theorems in the discrete geometry module.

While classification through theorems was not explicitly discussed with students, they did engage in classification during the discrete geometry module. For example, they divided the worlds into those with even and odd points. This classification was not done apriori as was done



when classifying triangles in the triangle theory building module. Rather, it was a classification which emerged as a result of interesting statements which turned out to be true.

In the triangle theory building module in Indus, there was a discussion on classification of triangles. In chapter 3, I discussed the classification of squares and rectangle. The answer to the question of whether a square is a rectangle rests on similar arguments to the answer to whether an equilateral triangle is an isosceles triangle. It depends on the concept of logical inheritance. The other type of triangle classification was to do with angles. In the discussion, two things were in focus. The first was that the classification depended on theorems – that there cannot be a triangle which is both obtuse and right at the same time depends on the claim that the sum of the angles of a triangle is two right angles. The second was on why a classification hinging on right angles is better than one which hinges on another arbitrary angle.

**Undefined entities and axiomatic systems.**

While axioms were not explicitly discussed in the discrete geometry module, they were present implicitly. However, in the triangle theory building module, there was some discussion on axioms and undefined entities towards the end of the module.

The way I introduced axioms and undefined entities in the triangle theory building module was by getting students to engage in cyclic reasoning. To define straight lines, they used collinear points. To define collinear, they used straight lines. There must be a way to escape this loop. The way to do that is to stipulate that straight lines and points are undefined entities and set up rules or axioms to govern their relationship with each other.

In the discrete geometry module, axioms and undefined entities were not explicitly discussed. However, they were implicit in the discussion. For instance, the stipulation that the world has six points or the neighbor relation in a particular world. If this module is run for



longer, these can be made explicit, and students can spell out the foundations of these worlds more rigorously than was done in the course implementation.

**Summary.**

The concepts of theory building form a framework for the development of theory building courses alongside the questions of theory building, which engage with different aspects of these concepts. As I mentioned in the introductory chapter, there are various aspects of these concepts which I have not addressed in this dissertation. For instance, the distinction between explanatory and non-explanatory proof. Future work in these areas needs to address these to get a more well-rounded understanding of theory building in education. In the next section, I will shift focus from learnings directly related to the concepts of theory building to other learnings from the implementation of the course.

**Other learnings from course implementation**

This section focuses on other learnings from the implementation of the course, with a focus on pedagogical choices. I will first discuss the difference between the two main modules of the course, with a focus on student feedback. I will then discuss the potential role of representations in a course of this type. Following that, is a discussion on Productive Struggle and its role in the course after which I discuss the need for the facilitator to be a learner. I will then talk about the considerations the facilitator ought to keep in mind when following students' direction, interpreting student responses, and using terminology.

**Implementation of triangle theory building vs discrete geometry.**

Triangle theory building and discrete geometry, the way they were implemented, have very different focuses even though they share many of the same learning outcomes. The focus of



discrete geometry was an exploration of new worlds while triangle theory building focused on rigorous reasoning within a world that students were already familiar with.

It isn't very surprising that different students had different reactions to the two modules. The common reason given by students for preferring discrete geometry is that it was fun to explore while triangle theory building was repetitive and hence boring. Those who preferred triangle theory building had a very different response. There were a few students who indicated that it gave them a better understanding of something they already knew and others who said that they were able to contribute to discussions in the module while they found discrete geometry too hard.

There was an asymmetry between the two schools. The students of Ganga tended to prefer discrete geometry while the students of Indus preferred triangle theory building. It could be the case that this was a matter of peer influence – one student had a preference and discussed it which in turn influenced how others saw the modules. While that may be true in some cases, the reasons students gave for their preferences are quite compelling – triangle theory building is repetitive and discrete geometry is new and hence requires effort to wrap one's mind around even though it may be more exciting.

Assuming these preferences represent something real, there could be many reasons for them including student background, how I conducted the particular sessions, and how the grouping of students was done, amongst many other reasons. I do not have any way at this point to make a case for one or more of these reasons. However, this difference between the two reactions is valuable to keep in mind, not just for future research but also for educators wanting to teach a course such as this one.



**Representations.**

Different elements of mathematical practice had different levels of focus in this course. For instance, the students in Ganga did not discuss classification and students in Indus spent a short amount of time on it. On the other hand, a large amount of time was spent on extending definitions in the discrete geometry module.

Students dealt with various types of representations at various points in the course. For instance, they worked with trees and Venn diagrams when discussing classification, they used graphs to represent worlds in discrete geometry, and they used tree diagrams to represent assumption digging. Representations are obviously valuable when constructing knowledge, whether in mathematics or outside of it.

However, representation as a learning outcome is something I largely ignored in this course. It is also an area I did not unpack when discussing the learning outcomes associated with theory building. Students didn't really create their own representations. Even in the case of discrete geometry, I gave them the representations. I also was not explicit when talking about them. I did not even make a distinction between a representation and the idea it represents. This is an important area for research to improve the course and to improve our understanding of mathematical practice.

Concepts like translating between representations are understandings which could form an important part of a different version of this course.

**Productive struggle.**

As I mentioned earlier in the thesis, productive struggle has its roots in Buddhism, Stoicism, and in many other intellectual traditions. The idea is that learning happens when students struggle with something. Students struggle when solving the same problem repeatedly



on tests and homework. That is not the type of struggle which is of value here. We may also struggle if we are working on something which is beyond our capacity to even understand – for example, asking a five year old to learn General Relativity. That is also not an example of productive struggle.

Heibert & Grouws (2007) use struggle to mean that students expend effort to make sense of mathematics, to figure something out that is not immediately apparent. Placing a student in a situation filled with perplexity, confusion, or doubt, as Dewey (1929) said, is crucial to this endeavor. If we wish to move away from an education system where quick answers are given paramount importance, we need to create such situation so that students learn to resolve these incongruities.

The theory building activities are designed to create productive struggle. It is unclear how far they succeeded in the implementation. In fact, it seems like some of the students found discrete geometry too hard. One aspect of challenging mathematical tasks (Sullivan, 2016) is that tasks should be accessible at multiple levels. While reworking the materials for this course and for other theory building courses, this ought to be kept in mind.

**A learning facilitator.**

There are many aspects of the implicit pedagogy used in the implementation of this course which are similar to those of Problem Based Learning (PBL) which makes use of ill-structured problems as a stimulus for learning (Barrows, 2000). In PBL, the facilitator is an expert learner, able to model good strategies for learning and thinking, rather than providing expertise in specific content (Hmelo-Silver & Barrows, 2006).

Similar to PBL, the role of a teacher in a course of this type is to eventually become irrelevant once students are able and ready to figure things out on their own. Till that point,



teaching a course of this type requires us to learn on the job in various ways. There is pedagogical learning we pick up from interacting with students, getting feedback from them, and reflecting on that.

Then, there is mathematical learning. Many of the claims students came up with during the course were not necessarily things I had spent time thinking about before they came up. More often than that, there were arguments students presented which were not clear or were flawed. A facilitator needs to be able to react to that. That doesn't mean that they need to be correct. A learning facilitator needs to be willing to work with students to figure out things together. Without a facilitator who is willing to learn from the students and be wrong, it is hard to see how a course like this one could achieve its goals.

### Following students' direction.

When facilitating a course of this type, one of the most important considerations is the level of involvement of the facilitator in deciding the direction of the course. Students ought to be the ones constructing knowledge and as a general principle it is valuable to follow their direction. However, as I pointed out in the podgon and triangle theory building modules, this can result in missing out on certain learning outcomes. For instance, in the triangle theory building module, classification was not pursued in Ganga since it was not on the initial list of claims students gave me.  Another downside of following the direction students take is that it requires the facilitator to react in real time potentially to claims and arguments they may not have seen before.

One solution to this is to start with a facilitator-controlled discussion and slowly start giving students control. Another possibility is to extend the length of the course. For instance, if



the course was semester long, as Fawcett's (1938) and Healy's (1993) were, it would not matter if certain learning outcomes were not addressed in the first few days.

**Interpreting student responses.**

While interpreting student responses is inevitable, there were a few times in the course where it appears as if I over-interpreted. At times, especially in the triangle theory building module, I used phrases like 'you mean …' and went on to say something which had a very tenuous relationship to what they had said.

Doing this once the students are comfortable with theory building would be acceptable. By that point, they would be able to push back against interpretations which did not match what they had said. However, early on in their learning journey, they may accept an interpretation even though that is not what they had meant.

Probing into their statements is probably a better strategy early on – asking questions like, "what do you mean by 'x'?". If students are unable to answer, rather than giving a single interpretation, the facilitator can give multiple interpretations and allow students to choose between them.

**Use of terminology.**

During the implementation of the course, I had decided to avoid mathematical terminology as much as possible. I barely used terms like 'axioms', 'undefined entities', 'degenerate', 'proof by contradiction', and so on. Avoiding using these terms was based on the goal of getting students to actually understand the concepts behind these terms and the assumption that words can often become a substitute for understanding. It isn't necessary that even seemingly 'correct' usage of terms is a sign of understanding. Without a word, a student is forced to think carefully about the underlying concept.



However, on the other hand, words can function as useful thinking tools. Imagine not having a term for what we call a triangle, and constantly thinking about the definition when wanting to use the concept. Both of these considerations need to be balanced when teaching a course of this type and these are important to keep in mind in future iterations of this course.

**Connections with research and education**

In this section, I focus on the relation between theory building and other areas of research and education. I discuss mathematical practice, transdisciplinarity, the Indian national curriculum, the Common Core state standards, mindsets, and equity.

### Mathematical Practice.

As I mentioned in the introduction, mathematical practice in the philosophy of mathematics is a relatively new area of research. It is not just concerned with what epistemologists are usually inquire into – things like how we can know the nature of mathematical objects, etc. It is also concerned with how mathematicians engage in their craft. It is partially descriptive using examples of how mathematicians in the past have come to conclusions. As I have mentioned earlier, research in this area is important to inform theory building activities.

Weber, Mejia-Ramos & Dawkins (2020) raise concerns about the relation between mathematical practice and education research. A few of their concerns are to do with research into mathematical practice itself while others are to do with the relation between mathematical practice and education. Three concerns they raise regarding inquiry into mathematical practice are what they call *the mathematical community identification problem, the heterogeneity problem* and *the accuracy problem*. The first is to do with identifying which mathematicians we are referring to when we say 'mathematicians' engage in some practice. The second is that there



is a difference between how different mathematicians approach their practice, and the third is that mathematicians may not be accurately reporting their practices.

Some of their other concerns are to do with education research. They include *the advanced content problem* and *the resources problem*. The advanced content problem is to with many education researchers not having the background to engage with mathematics which mathematicians engage with. This is a concern whose solution requires mathematicians and philosophers to work closely with educators so that the valuable learning outcomes associated with mathematical practice. There are also a number of education researchers who do have backgrounds in mathematics, some of whom are professional mathematicians as well. Contributions by them would be valuable.

The resources problem is to do with the appropriateness of certain mathematical practices to students since they may lack some of the resources which professional mathematicians have. I would add to this that even if students have the resources to learn a particular practice, it is not obvious as to whether that practice is valuable for students to learn. Of course, those mathematical practices which are not appropriate for students ought not to be included in a classroom.

While there are many important concerns about mathematical practice in mathematics education, we need to be focused on the educational goals. If we constrain ourselves to general education, as opposed to specialized education aimed at developing mathematicians, mathematical practices only serve a purpose inasmuch as they give students certain tools which they can use in their public, professional and personal lives regardless of what career they decide to pursue. I will discuss transdisciplinarity in the next section which discusses the relationship between mathematical practices and inquiry practices in other disciplines.



**Language and mathematics.**

Various curricula mention 'learning to communicate mathematically' as one of the goals of mathematics education. One aspect of that is to be able to communicate mathematical ideas to an audience at the right level of sophistication. Another aspect is to communicate clearly and rigorously – using that both to understand mathematics better and refining those abilities while engaging in mathematics. A large part of the course was focused on clarifying words and sentences – understanding the relationships between the meanings of words and the consequences of the propositions encoded in sentences containing those words. Both for students to be better prepared for a course like this one and for students to be able to transfer their learning from this course, there needs to be a partnership developed with language teachers and educators.

**Transdisciplinarity and theory building.**

I use the term transdisciplinarity in a way which distinguishes it from inter and multi-disciplinarity. I use interdisciplinarity to refer to the intersection between two disciplines – for instance, physical chemistry. Multidisciplinarity refers to the practice of using multiple disciplines to solve a problem. For instance, to understand climate change, we need tools from disciplines like statistics, physics, sociology, geology, and so on. Transdisciplinarity refers to a concern for concepts and tools of knowledge constructions which are abstractions of particular disciplines. For instance, the concept of theory, while not exactly the same in different disciplines, has some level of commonality. That common concept is a transdisciplinary concept.

While this thesis is situated within mathematics and mathematics education, my eventual goal with this material is to create a theory building course which cuts across disciplinary boundaries, including scientific, mathematical and ethical theories. Scientific and mathematical



theories are different in many ways but also share much in common. For instance, definitions and reasoning are involved in building both, but scientific theories allow for different modes of reasoning and different types of definitions. Justification also works differently for differently for different types of theories. Mathematical theories are not justified apart from on the grounds of internal consistency, while scientific theories are justified based on observation. However, conclusions from the assumptions of the theory are justified in both by showing the steps of reasoning from assumptions to conclusions.

There are also tools which cut across these different types of theory building. For instance, assumption digging is something which is useful everywhere – asking why we believe something and then asking why we believe the assumptions made in the justification is a major part of any discipline and of activities such as reading a newspaper.

Integration is one of the main goals of any type of knowledge construction – without seeing similarities between things and abstracting those away, knowledge would just be a set of disconnected facts. The idea of taking a theory, changing some aspects, and figuring out the consequences, can also be an important tool. This is especially evident in ethical theory building when we may want to consider the consequences of exchanging one principle with another, possible a principle somebody else holds and we don't.

A single course which involves scientific, mathematical and ethical theory building would allow students to be able to see the connections and differences between these endeavors. It is not just incumbent on students to see these connections, but also on the instructor to make these connections explicit as far as possible. To do this rigorously, this would require creating a framework for such a course with help from experts in various disciplines including philosophers of science, scientists, ethicists, and so on.



I have previously done sessions or courses where I move between different disciplines. In what follows, I will demonstrate some examples of transdisciplinarity in a classroom, which I have done in the past with students.

### Classification in a transdisciplinary manner.

Classification as a tool lends itself naturally to different disciplines. I have often done the following activity: I start with classifying shapes and then have a discussion on classification of biological organisms. Many of the considerations while doing these two classifications is similar.

I usually start with asking whether a square is a rectangle (something similar was done with classification of triangles in the triangle theory building module). The choice between the two classificatory systems implicit in this question is very similar to the choice between two classificatory systems of biological organisms: one where humans are animals and the other where they are not. Things like purpose of classification, logical inheritance, and simplicity come up when discussing these choices.

### An example of scientific theory building.

In 'The Evolution of Physics' (Einstein & Infeld, 2007), Einstein gives the analogy of constructing a scientific theory being similar to figuring out the mechanism of a clock on a wall. Imagine all you could see was the front of the clock. You might be able to move the dials, but you cannot open the back. You have to guess what the mechanism is based on observing the behavior of the dials of the clock.

An activity I have created on scientific theory building is what I call 'Theoretical Anatomy.' The idea is to imagine yourself as a human thousands of years ago. You know nothing about the internal workings of the human body and for legal reasons in your community,



you are not allowed to cut open a human body. So, the task is to figure out how the human body

works through external observations and through non-invasive experiments.

The way this works is that you postulate some claims about the working of the human

body, deduce consequences of those claims, and then compare the consequences to observations.

There are similarities and differences between this activity and theory building in mathematics.

The similarity is in the structure and in the fact that you are using reasoning to deduce

consequences of assumptions. The differences are that you are allowed to use types of reasoning

which would not be valid in mathematics, such as abduction or probabilistic reasoning, and that

you at the end compare the consequences to reality.

### *An example of ethical reasoning.*

Reasoning about ethics is much more complicated that mathematical reasoning.

However, there are relationships between the structures of arguments. The idea of 'ethical

principle' is not completely different to the concept of 'axiom'.

Let us take a situation often used in ethics courses – the idea of euthanasia. Suppose you

were transported back in history to a place where a widow was being burnt alive on her

husband's funeral pyre. You are too far away to save her, but you have a sniper rifle, and you are

a good shot. Would you kill her?

There are two competing ethical principles here. One is 'killing people is morally wrong'

and the other is 'reducing suffering is morally right.' The consequence of the first principle is

that you should not kill her but the consequence of the second is that you should. Deducing these

individual consequences is similar to reasoning in mathematics. However, at this stage, in order

to come to a conclusion, we need to choose which principle to follow given the context.



In this context, many people will conclude that killing her is the morally right thing to do. However, suppose we change the situation to somebody who stubbed a toe and you could kill them without any pain. In that situation, the first moral principle would hopefully prevail.

The important point to see here is that there are relationships between reasoning in ethics and reasoning in mathematics. Ethical principles are analogous to axioms even though they are different in many ways. The deductive process of arriving at consequences from those principles also bears a similarity with mathematics even though there are significant differences as shown above.

**Mindsets.**

One of the goals of this course, and of transdisciplinary courses of this nature, is to develop in students the ability to construct and evaluate knowledge. The deeper, societal goal implicit in that is to have a more reasoned society which is capable of coming to conclusions for the good of humanity and the planet. One thing required for this goal which not discussed in this thesis is an ethical sense. The tools on knowledge construction and evaluation can be used for positive or negative ends – without an ethical sense (as distinct from the ability to reason about ethics) arming students with these tools could have potentially harmful outcomes in the long run. Since this is not something I have worked on, I will leave this to others to discuss.

Apart from this, to achieve these goals, certain mindsets need to be developed. It is quite clear from the psychology literature that all of us engage in motivated reasoning (Kunda, 1990; Epley & Gilovich, 2016). We also don't change our minds as often as we should. We also accept justifications and conclusions from authority figures, whether nation states, corporations, politicians, religious leaders, academics or secular gurus, without subjecting them to the level of critical thinking we ought to.



These seem to be natural properties of our minds and we may not be able to rid ourselves of these traits. The best we potentially could do would be to develop certain mindsets which remind us to push back against these properties of our minds.

For instance, the mindset of being okay with being wrong – in fact even reveling in admitting being wrong. Alongside that, the mindset of wanting to put our own beliefs through critical evaluation, both by ourselves and by others. Closely connected to that is the mindset of intellectual skepticism. Intellectual skepticism here refers to the habit of not blindly accepting claims made by authority. This does not involve blindly rejecting claims made by authority as is the wont of many a conspiracy theorist. Rather, subjecting the belief to an adequate level of skepticism and critical evaluation.

Unpacking these mindsets and developing learning materials to address them is outside the scope of this thesis. However, without developing these mindsets, along with the ethical sense mentioned above, the work in this thesis is largely pointless.

**Existing curricula and theory building.**

In this section, I will focus on India's national curriculum and the US Common Core state standards.

India's National Educational Policy, 2020 (NEP) is a document which is intended to feed into curricular design and eventually into a new National Curriculum Framework (NCF). It states, "... education must develop not only cognitive capacities - both the 'foundational capacities' of literacy and numeracy and 'higher-order' cognitive capacities, such as critical thinking and problem solving – but also social, ethical, and emotional capacities and dispositions" (NEP, 2020). These 'higher order' cognitive capacities are not explicitly defined in the document.



However, assuming that the new NCF is going to build upon NCF 2005, it is useful to look at the goals of mathematics education in mathematics position paper for that curriculum framework. It mentions the need for 'abstraction, structuration and generalization' (NCERT, 2006). It also mentions argumentation, reasoning and perceiving relationships as important goals of mathematics education along with students discussing and communicating mathematics. Given these goals, as well as the broader goals of critical thinking mentioned in NEP-2020, a course on theory building like the one presented in this thesis is consistent with them.

Added to that, NEP-2020 also mentions the need for 'holistic' education across disciplines 'in order to ensure the unity and integrity of all knowledge.' A transdisciplinary course on theory building built upon this course could help in achieving that.

The mathematical practices mentioned in the US Common Core state standards include making use of structure, constructing and critiquing arguments, searching for regularity and attending to precision (Standards for mathematical practice, n.d.). These practices are an important part of the theory building course presented in this thesis.

In both the cases of the Common Core and the NEP/NCF, there are many goals in common with this Theory Building course. However, there are also practice related goals unpacked in Chapter 3 which are not explicitly mentioned in either of these and could be valuable learning outcomes to pursue.

**Equity and theory building.**

Mathematics education, as it exists in many places, excludes many people. For instance, the way exams in India are designed, the people who tend to do well in mathematics exams are those who can withstand un-creative, repetitive, problem solving, following pre-set rules they have to memorize. Apart from that, those who can afford expensive tuitions have a significant



advantage over the rest of society, especially those from historically marginalized groups who don't tend to have a history of formal education in their families.

Education in India is focused on those who do well in the system, largely ignoring those who do not. It is also focused on those who are looking to do medicine or engineering. To change this, a significant part of education needs to be valuable to everybody, no matter what they later choose as their specialization. There is no reason calculus needs to be introduced in high school to everybody who does mathematics. Rather, the focus ought to be on thinking abilities which would be useful to everybody in their personal, professional and public lives even if they choose not to pursue a career in a STEM related discipline.

**Summary and future research**

This thesis can be thought of as having three parts. Firstly, unpacking the learning outcomes associated with theory building including the questions of theory building. Secondly, the design of a course aiming at these learning outcomes. Finally, a qualitative description of the implementation of the course with two groups of students.

In unpacking the learning outcomes associated with theory building, I made use of the education literature as well as the philosophy of mathematics literature on mathematical practice. As I mentioned, this is a relatively new area of research stretching back only a few decades. Hence, there is a lot of work to be done in thinking about these aspects of mathematical practice, both in their own right and for the purposes of education. This is a research program which could benefit from close collaboration between educators and philosophers.

Designing the course and unpacking the learning outcomes went hand in hand. There were insights gained from unpacking the learning outcomes which fed into the design, and insights from thinking about the mathematics which helped shed light on the learning outcomes.



This course can be improved further using some of the insights from the course implementation and from thinking more carefully about mathematical practice. As I mentioned while discussing transdisciplinarity, this course can also be extended to include areas outside of mathematics while aiming at a similar set of goals. This requires collaboration not just with academics exploring mathematical practice, but also those exploring practice in other disciplines.

The implementation of the course was intended to be descriptive. Since this was an implementation in two schools, part of the same school group, by a single facilitator, there are probably not many general empirical claims which can be made. In terms of conclusions we can draw from this, at best this serves as an existence proof to show that it is possible to run a course of this type achieving some of these learning outcomes. Even that depends on the reader's judgment of whether they see any valuable learning in the interactions. So, in order to learn more about the value of this particular course, it needs to be tried out in various situations with other educators teaching it.

To return to the goals of education, if we want to achieve Dewey's vision of education I started the thesis with, we need a system of education which has value to students no matter what they decide to do with their lives. Hence, the focus needs to be on ways of thinking which are generalizable and transferrable outside of their domain of learning. The thinking tools associated with mathematical theory building, and mathematical practice more generally, can form an important part of that endeavor.



**Appendix A: Lesson plans**

The following are the lesson plans I used while conducting the course. I have included a plan titled 'Straight Lines and Intersections,' which I had to cut from the course due to unforeseen circumstances. However, the rest of the plans were used. While conducting the course, the plans were mostly bullet points I took into the class. I have expanded on them in this appendix.

Apart from the first plan (on definition guessing), the rest of the plans are written as dialogues. These dialogues are intended to simulate a classroom discussion. Of course, no actual classroom discussion will ever look anything remotely like these plans. However, such plans hopefully give a sense of how such a class can be conducted.

The plans are:

- Definition of Podgon

- Triangle Theory Building

- Discrete Geometry



**Definition of Podgon**

In this activity, students guess definitions based on examples. This session can either be done with two instructors or a single instructor playing the role of two people. In this document, I will assume there are two instructors. The word podgon is not something which you would find in the dictionary or in a math textbook. It is a word I made up. Let us refer to the two instructors as A and B. A has the following definition of the word podgon, which I will refer to as definition A:

*A podgon is an equilateral polygon with an odd number of sides.*

B's definition (referred to as definition B) is:

*A podgon is a closed shape with an odd number of straight line sides. (it could have curved sides, but the number of straight line sides is odd)*

We start with telling the students that a podgon is a type of closed shape, and we give them the following information (Podgon A is the object defined by definition A and the same for Podgon B):

| Shape | Podgon A | Podgon B |
|---|---|---|
| 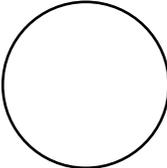 | 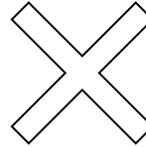 | 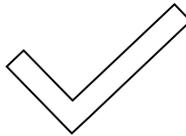 |
| 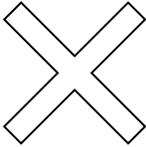 *All sides equal* | 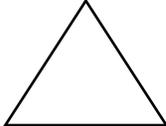 | 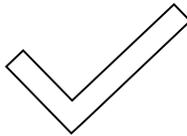 |



The task for students now is to come up with a definition of podgon which matches this table. The instructors are not allowed to say that the definition the students have given is the same one they have. Once they have definitions in place, students need to come up with examples which test whether the definition they have come up with is the same as those of the instructors. Some possible pictures they may come up with:

| Shape | Podgon A | Podgon B |
|---|---|---|
| 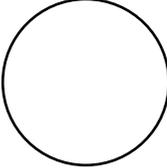 | 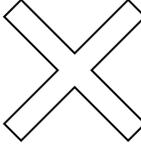 | 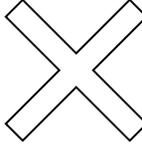 |
| 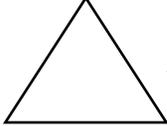 *All sides equal* | 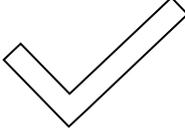 | 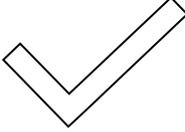 |
| 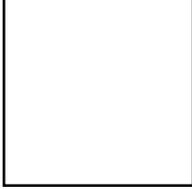 | 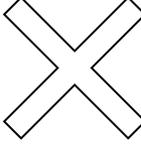 | 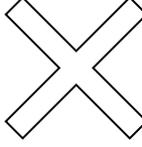 |
| 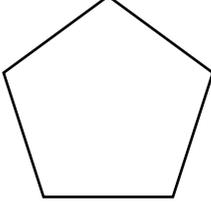 | 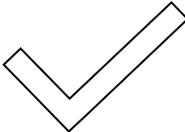 | 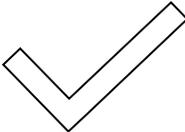 |
| 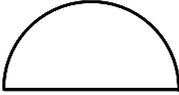 | 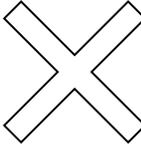 | 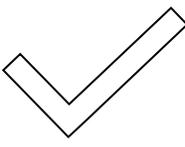 |



| 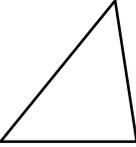 | 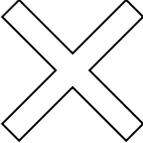 | 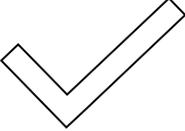 |
|---|---|---|
| 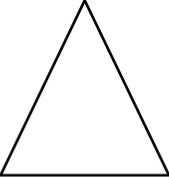 | 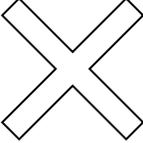 | 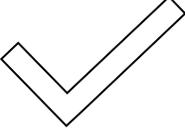 |
| 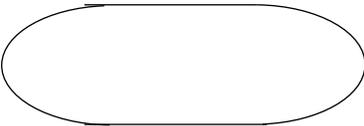 | 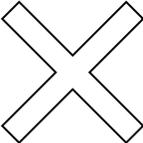 | 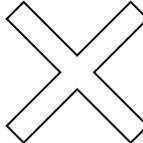 |
| 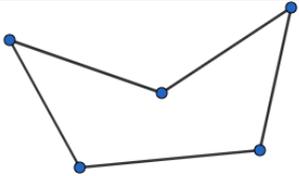<br>*Equal Sides* | 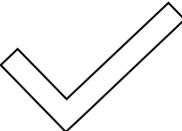 | 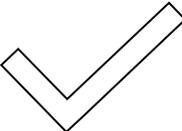 |
| 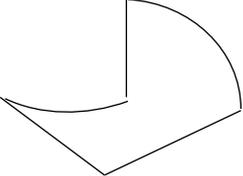 | 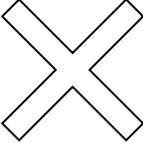 | 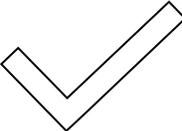 |



**Triangle theory building**

The goal of this module is to help students develop a theory of Triangles in Euclidean Geometry. The idea is that students list out various things they believe to be true about triangles, and then construct a theory by linking claims to one another through proofs. This module could be something done over a few hours or a few months - the shorter the time frame, the more selective you will have to be. In this module plan, I will be working through various directions this can be taken in. However, this is just a small subset of possible directions to take this. Also, the plan will be written as a series of dialogs between a teacher and students. No actual class will look like the dialog - the goal here is to give the teacher a sense of what is involved so that they can chart their own path.

I will be using the letter T to refer to the teacher. Any other letter refers to students.

**Listing out what we know.**

The first thing is to ask students to list things they believe to be true about triangles. At the end of this exercise, we should have 10-15 claims on the board. The teacher will have to clarify what students mean when they are making claims so that everybody understands what is being said. Here is a possible list:

1. The sum of angles of a triangle is 180 degrees

2. There are many types of triangles: isosceles, equilateral, right angled, scalene, obtuse and acute

3. The angles opposite equal sides of a triangle are equal

4. The angles of an equilateral triangle are equal

5. Any triangle can be circumscribed

6. The altitudes of a triangle meet at a point



7.  The medians of a triangle meet at a point

8.  The perpendicular bisectors of a triangle meet at a point

9.  The area of a triangle is 1/2 base x height

The sum of the lengths of two sides of a triangle is greater than the length of the third side

**Equilateral Triangles.**

T: Let us pick up 4: The angles of an equilateral triangle are equal. Why should I believe that? Can you try and show that is true using other things written on the board?

*Time for group discussion*

A: It is because the angles opposite equal sides are equal.

T: Okay. You have mentioned one of the other things we wrote. But, I do not see the connection between that and what we are trying to show. You need to show me the steps of reasoning involved starting from what you are assuming to what you are trying to show.

*A goes to the board and starts writing*

The angles opposite equal sides of a triangle are equal (3)

An equilateral triangle has equal sides

So, the angles of an equilateral triangle are equal

T: Okay. That makes sense. Notice that the second thing you said is not currently in our inventory: that equilateral triangles have equal sides. Let us add it in and explore it. Why should we believe that equilateral triangles have equal sides?

B: Isn't that what equilateral triangles are?

T: Are you saying that is how you want to define equilateral triangles?

B: Yes.



T: Great. In that case, if we assume 3, we can conclude that the angles of an equilateral triangle are equal. Since we defined equilateral triangles as triangles with equal sides, we concluded that triangles with equal sides have equal angles. An obvious question to ask is whether triangles with equal angles have equal sides. Is that always true or can you create a triangle with equal angles which does not have equal sides.

C: It is always true. Triangles with equal angles have equal sides.

T: Why should I believe that?

C: Can't we use the same argument as above?

T: Can we? The assumption we made for that argument was that angles opposite equal sides of a triangle are equal, not that sides opposite equal angles are equal.

C: Why don't we add that assumption in?

T: Sure. We can do that. Once we have done that, we get:

    The sides opposite equal angles of a triangle are equal

    So, the sides of a triangle with equal angles are equal

I want to be a little careful here. You say that the sides opposite equal angles are equal and from that you have concluded that all the sides are equal. Implicit in that is that all the sides are opposite to the equal angles. Is that true?

B: Isn't that obvious?

T: It may seem obvious. However, what we are trying to do here is to be explicit about our assumptions.

B: Since there are three sides and three angles, each side is opposite an angle.

T: Triangles have three sides and three angles?

B: Yes..



T: We haven't written that on the board yet. Let us add those statements in. Even assuming those two things, I'm not sure whether your argument works. Just because there are three sides and three angles does not mean that each side is opposite an angle. Maybe there is a side which is not opposite any angle.

B: That isn't possible

T: Why?

B: Maybe we can say that the side of a triangle which does not make up the angle is opposite it.

T: Why would such a side exist?

B: There are three sides and two sides make an angle. So, there will always be a third side left over.

T: Two sides make an angle? We need to add that in to our inventory. However, given that, I accept your argument. You have shown that there is always exactly one side opposite to an angle of a triangle, and that given a side, it can only be opposite to one angle. Hence, since there are three sides, each must be opposite an angle.

We have now shown that triangles with equal sides have equal angles and triangles with equal angles have equal sides. So, we can now define equilateral triangle in two ways:

A.  Triangles with equal sides

B.  Triangles with equal angles

These are called 'equivalent definitions.'

The obvious thing to ask you now is to justify some of the premises used in the two proofs you have given so far. However, I want to move on to something else and we will get back to that.



**Types of triangles.**

T: You mentioned various different types of triangles. You said triangles are equilateral, isosceles, right, obtuse, acute and scalene. If we were to classify triangles, would you classify them like this?

(Tree with all the types directly underneath triangle)

B: No, there are two different thing happening here. One classification is based on sides while the other is based on angles. Equilateral, isosceles and scalene go together while right, obtuse and acute go together.

T: So, let's start with equilateral-isosceles-scalene. Is this the classificatory system you are looking for?

(All below triangle directly)

C: Yes.

A: No, equilateral triangles are types of isosceles triangles so, this is the system we should use:

(Add tree)

C: That's not true.

A: It is.

T: Well, we can keep going back and forth asserting one or the other claim. I want you to give reasons.

C: Equilateral triangles have all sides equal while isosceles triangles only have two sides equal.

A: I agree that isosceles triangles have two sides equal. However, equilateral triangles also have two sides equal. They also happen to have a third side equal to the other two.



T: C, look at this definition of isosceles triangles:

Triangles with two sides equal

If we accept this definition, do you see that equilateral triangles are types of isosceles triangles since equilateral triangles do fit this definition – there is no mention that the third side cannot be equal.

C: I accept that.

T: So, if we want equilateral triangles to not be types of isosceles triangles, we have to change the definition to something like:

Triangles with two equal sides and the third side not equal to the other two.

If we use this definition, it follows that equilateral triangles are not types of isosceles triangles. So, why should we choose one definition over the other?

B: I don't know.

T: An advantage of using C's definition is when you want to describe some object you see in the real world. It would be weird to refer to a triangular shape with equal sides as a type of isosceles triangle. If there were a set of tables shaped like triangles, some of which were equilateral while the others were isosceles, you would probably say something like: there are some tables there which are equilateral and others which are isosceles.

However, the advantage of A's definition is that if we show something to be true about isosceles triangles, we automatically know it is true for equilateral triangles since equilateral triangles are types of isosceles triangles. Under C's system we would have to justify the claims separately for both types of shapes. So, assuming you want to minimize the amount of work you need to do, it would make sense to pick A's definition and classificatory system if we are trying



so show things to be true about types of triangles. Since that is what we are exploring in these sessions, we pick A's classificatory system over C's.

B: Does that mean that we would want isosceles triangles to be types of scalene triangles?

T: That's a great question. Let's get back to that once we have gone through the other type of classificatory system – the one based on angle measures. Is this acceptable?

A: Yes. In this case, right triangles are not types of acute or obtuse triangles, and acute and obtuse triangles are definitely different.

T: Okay. What are the definitions of these three?

A: Acute triangles are triangles where all angles are acute.

Right triangles are triangles with one right angle.

Obtuse triangles are triangles with one obtuse angle.

T: Okay. Is this a complete classification. By that I mean, does every triangle fit into one of these categories?

A: Yes.

T: On what basis do you say that. Let me give you the argument for the previous classification. The previous one had isosceles triangles and scalene triangles at the top level. Since every triangle has either at least two equal sides or no equal sides, every triangle must fit into one of those classes. In this case, how do we know that every triangle fits into one of these?

A: Triangles either have three acute angles, one right angle or one obtuse angle.

T: How do you know that there are no other types of triangles. For example, why is there no triangle with two right angles?



A: You cannot have a triangle with two right angles since the sum of the angles of a triangle is 180 degrees. If you have two right angles, the third angle will be zero degrees and triangles cannot have zero degree angles.

T: Let me add that in to our inventory – triangles cannot have zero degree angles.

A: Similarly, you can have at most one obtuse angle in a triangle, and you cannot have a right angle with an obtuse angle.

T: Your argument seems to be reliant on the claim that the sum of the angles of a triangle is 180 degrees. If this wasn't the case, your conclusion would not follow, and the classificatory system might not be complete. I just wanted to point that out.

I have another question about this classification. Why not replace it with the following classification:

(Add 101 degree classification)

This is also a complete classification. Why are we using right angles rather that 101 degree angles?

B: Our textbooks use right angled triangle and we have heard about right angled triangles since we were kids.

T: In what ways have you heard about right angled triangles?

B: We learnt about Pythagoras's theorem.

T: Ah.. you learnt about theorems about right angled triangles, but you haven't learnt any theorems about 101 degree angled triangles, right?

B: Yes.



T: That is a good reason for why you choose right angled triangles as an important type of triangle rather than 101 degree triangles – there are interesting theorems which are specifically about right angled triangle, which there are not about 101 degree triangles.

Another reason is that there are theorems about acute angles, but I can't seem to think of interesting theorems about triangles with all angles less than 101 degrees.

A: Okay. So, that means that when we define something, we need to make sure it is important – that there are theorems about the things we define.

T: Precisely. So, let me go back to the earlier classification and ask you whether you know interesting things about scalene triangles.

A: Aren't all the angles of a scalene triangle different?

T: Yes.

A: I can't think of others.

T: Notice there are some things we can say about 101 degree triangles as well – we can say that they are not equilateral and that the sum of the other two angles is 79 degrees. However, these are not anywhere close to as valuable as the things we can say about right angled triangles. I'm going to leave to you as to whether scalene triangles are just a useful term to refer to non-isosceles triangles or whether it is actually an interesting concept with useful theorems.

Let us go back now to justifying the assumptions we used when talking about equilateral triangles.

**Isosceles Triangle Theorem.**

T: We had assumed the following claims in our earlier arguments:

1.  Angles opposite equal sides of a triangle are equal

2.  Sides opposite equal angles of a triangle are equal



Let's start with the first one. Why should we believe it?

A: We learnt a proof for this in school using congruent triangles. Let the triangle vertices be A, B, and C. Assume AB = AC in length. Then, drop a perpendicular from A to BC – assume it intersect BC at D. This creates two triangles, ABD and ACD. Angle ADC = angle ADB = 90 degrees. AB = AC and AD = AD in length. So, the triangles are congruent by RHS. Hence, by CPCT, you get that angle ABC = angle ACB.

T: Okay! There are a lot of things you said there – let us take it step by step. You are starting with what we were given – that two sides of the triangle, AB and AC, are equal. Then, you drop a perpendicular from A to BC. Why is it possible to do that?

B: It is possible

T: But, we don't have that written anywhere. So, let us add in that given a point and a line, you can always draw another line perpendicular to the first line, passing through the given point. You get two right angled triangles – ADC and ADB. Then, when you say RHS, you mean the following:

Two right angled triangles are congruent if their hypotenuses are equal and one of the other sides is equal.

Since one of the sides is common to both the triangles, it is equal to itself. Why are the hypotenuses equal?

B: The hypotenuse is the side opposite the right angle. In these two triangles, this would mean that AB and AC are the hypotenuses, which we assumed were equal.

T: Great! So, if we assume the RHS claim, it follows that the two triangles are congruent. So, how does that give us that the angles opposite the equal sides are equal?

B: CPCT



T: What does that mean?

B: Corresponding parts of congruent triangles are equal

T: You are saying that if two triangles are congruent, then the parts of them are equal? What exactly does that mean? If ADB is congruent to ADC, which parts are equal to which parts?

B: We already showed that AD and AD are equal, and AB and AC are equal. What we can say from congruence is that BD = CD.

T: Okay. But, what about the angles?

B: It is confusing – if we look at a picture, it is obvious but it is hard to state in words.

T: Let me help you. We know that BD and CD are equal, and that AB and AC are equal, so the angle between BD and AB must equal the angle between CD and AC, which means that the angle ABD = the angle ACD. But, we wanted to show that the angle ABC = the angle ACB. How do you get that?

B: Isn't ABD the same as ABC?

T: Seems like it but why?

B: Can we say that the angle between two line segments is the same if they are both part of the same line?

T: That seems fair and we can arrive at the conclusion we want to. I would like to point out that in doing this, we assumed an understanding of congruence. We haven't really defined congruence or shown that RHS results in congruence. However, we will get back to that at some later point. For now, let us try to prove the other statement:

Sides opposite equal angles of a triangle are equal.

A: Maybe we can use congruence again?



T: How would we do that?

A: What if we drop the same perpendicular?

T: Okay. But, this time we do not know that the hypotenuses are equal – in fact, that is what we will need to show.

A: We can use ASA instead. The angle ADB = the angle ADC = 90 degrees, the side AD is common and the angle DAC = the angle DAB.

T: I accept the first two. Why is the last one true?

A: We know that the sum of the angles of a triangle is 180 degrees, and two of the angles of the two triangles are equal – the right angle and the angle we are given is equal. So, the third must be equal as well.

T: I accept that. You are assuming that if a side and the two angles that side is a part of are equal, then the triangles are congruent. Since, they are congruent, the respective sides must be equal and hence AB = AC. Great! As I said earlier, we will get back to congruence. However, notice that you needed that the sum of the angles of a triangle is 180 degrees. We have needed to appeal to that before as well. So, let us move on to asking about that. Why should we believe that the sum of the angles of a triangle is 180 degrees?

**Sum of angles of a triangle.**

T: Why should we believe that the sum of the angles of a triangle is 180 degrees?

C: We learnt a proof of this: Draw a line parallel to one of the sides of the triangle passing through the opposite vertex. Let the vertex be A, the opposite side be BC, and pick two points on the line on either side of A, D and E. Then, angle ABC = angle DAB, and angle EAC = angle ACB. The sum of the angles DAB, EAC and BAC is 180 degrees, and hence the sum of ABC, ACB and BAC is 180 degrees.



T: Let me start by asking you why you believe that the sum of angles DAB, EAC and BAC is 180 degrees.

C: Because DE is a straight line and the angles make up the straight line.

T: What does DE being a straight line have to do with the sum of the angles being 180 degrees?

C: The angle of a straight line is 180 degrees – the angle DAE is 180 degrees.

T: Okay. So, what does that have to do with the sum of the angles being 180 degrees?

C: The three angles together make up the angles of the straight line.

T: I could ask you why that is, but I'm going to leave that for now. However, I will add to our list that the angle of a straight line is 180 degrees. Let me ask you about the other claim: that angle ABC = angle DAB, and angle EAC = angle ACB. Why is that true?

A: Since BC and DE are parallel, AB is a transversal which intersects both these lines. Since alternate interior angles are equal in the case of parallel lines and their transversals, we get that the angles are equal. The same is true for the other pair of angles.

T: What about the lines being parallel makes the alternate interior angles equal?

A: For any pair of parallel lines and a transversal, alternate interior angles are equal.

T: Okay. We can add that to our list of statements.

**Other areas to explore.**

Moving forward, there are many other areas which can be explored. Some of them are:

- Centers of a triangle.

- Defining straight lines.



## Discrete geometry

The following is a lesson plan for the discrete geometry module. In this, we start with a question which is not well specified. From there, we clarify the question by defining terms and setting up rules for the world(s) we are in.

### Bisecting a line in a world with six points.

T: In a world with six points, can every straight line be bisected?

A: We learnt in geometry that every straight line can be bisected.

T: That is the case in Euclidean Geometry, but does Euclidean Geometry have exactly 6 points?

A: No.

T: So, in a world with 6 points, can every straight line be bisected like it can in Euclidean Geometry?

A: I don't know.

T: Why don't you know?

A: I've never learnt about worlds with exactly 6 points.

T: Neither have I. I'm asking you to use you imagination.

B: But, what does a world with 6 points look like?

T: Good question. Let us take a look at the question a little more carefully. We know what words like bisect, straight line and points mean in Euclidean Geometry. But, what do they mean in a world with exactly six points? Without answering that question, we cannot proceed.

B: So, what do they mean?

T: You are probably used to doing mathematics where the teacher or textbook tells you what things mean and you have to answer questions given those definitions and rules. One part



of real mathematics you have probably never done is setting up the objects and rules of the world we are operating in.

B: How do we start?

T: Let me give you a possible direction. Have you seen science fiction films or TV shows like Star Trek which have teleportation? Or floo powder in Harry Potter?

B: Yes.

T: Think about points as places you can teleport to and from. You can think of two points as neighbors if you can teleport between them. Let's say there was a teleportation machine on Mars which was connected to teleportation machines in Delhi and Pune – you can go back and forth. However, there is no teleportation connection between Delhi and Pune. We can represent this as follows:

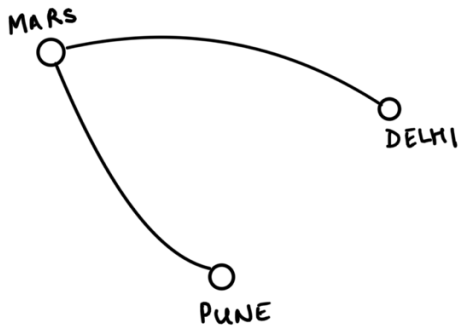

*Figure A 1 Teleportation*

Say, you wanted to travel by teleportation from Delhi to Pune. How would you do it?

C: Travel from Delhi to Mars and then Mars to Pune.



T: There might be many other ways to go about doing this, but for now let us stick to the analogy with teleportation. Now, we need to know what straight lines are. What is a straight line? Not just in this world but generally.

C: We could say that a straight line is made of some points in the same direction

T: Yes, that concept of straight line seems to fit with our notion of straight line in Euclidean Geometry. However, it is less clear what this could mean in a world with six points. If we want a notion of straight line in worlds with 6 points, it might be useful to think of other notions of straight line from Euclidean Geometry. We might be able to define some concept of distance with six points. Is there a concept of straight line to do with distance?

C: You can think of a straight line as the shortest path between two points.

T: So, in our earlier example, what is the shortest path between Delhi and Mars?

C: Just the direct teleportation

T: What about between Delhi and Pune?

B: Delhi to Mars and then Mars to Pune.

T: What if we now add a path from Delhi to Pune directly

B: Wouldn't it depend now on how long each path is?

T: Good point. If we think of it in terms of time, maybe the time it takes to get from Delhi to Pune via Mars is shorter than the direct route. Maybe the direct route is made using older technology. For now let us assume that each direct route between two neighboring locations takes the same amount of time or, in other words, has the same length.

B: In that case, the shortest path from Delhi to Pune now becomes the direct path without going through Mars.

T: So, that is now the straight line path, right?



B: Yes. If straight line path is the shortest path.

T: Okay, so now that we have a way of understanding what a world with six points is and what a straight line in that world is, can every straight line be bisected?

A: Doesn't it depend on what points we consider to be neighbors?

T: Before I address that, let us discuss how we can represent such worlds. Here is my suggestion of a way to do that: We can represent points by a dot and two points have an arc connecting them if they are neighbors. Now, let us explore your suggestion. Maybe it does matter points are neighbors in order to determine whether every straight line can be bisected. So, rather than me asking that question generally, let us take a specific example. Try seeing whether every straight line in this world can be bisected:

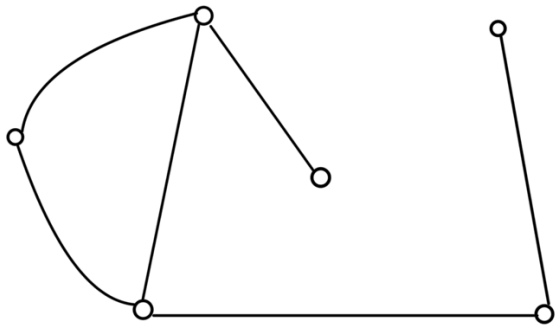

*Figure A 2 A world*

*After some time*

A: Not all lines can be bisected

B: I agree

T: Can you give me an example of lines which cannot be bisected?

A: Here (the dark path):



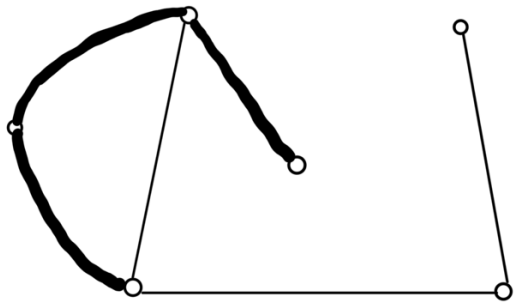

*Figure A 3 A path*

B: No! That line can be bisected. Here is an example of a line which cannot be bisected:

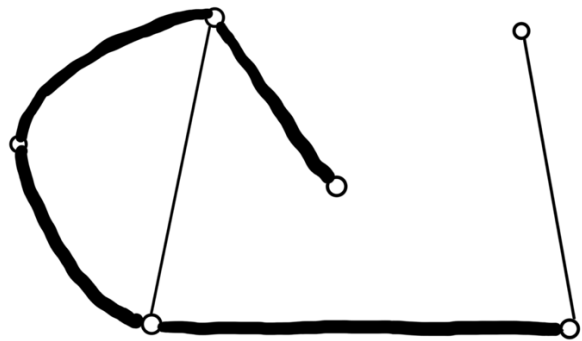

*Figure A 4 Another path*

A: But that line can be bisected. The line I made cannot.

T: This is interesting. Both of you seem to be pretty convinced you are right but you seem to be saying opposite things. Do you think one of you is wrong or do you think something else is going on?

C: I think both of them mean different things by the word bisect.

T: Expand on that.



C: It seems like A thinks a line broken into two equal parts at a point is bisection while B thinks a line broken into two equal parts at a connection between points is bisection.

T: Okay! So, let us say A-Bisection is the one at a point and B-Bisection is the one at a connection between points. Clearly, as they have pointed out, A's line cannot be A-bisected while B's line cannot be B-bisected. So, B, can the line you drew be A-Bisected?

B: Yes

T: A, can the line you drew be B-Bisected?

A: Yes

T: Okay, so as you can see the conclusion as to whether a particular line can be bisected depends on which of these two definitions you choose. However, as you have shown, for both the definitions there are some lines which cannot be bisected in this particular world with six points and these connections.

Now that we know that there is at least one world made of six points where every straight line cannot be bisected by both definitions of bisection we have so far, a natural question to ask is: are there worlds we can create where every straight line can be bisected? Let us try answering this question for both concepts of bisection starting with A-Bisection.

Can you create a world with six points, by changing the connections, such that all straight lines can be A-bisected?

A: How about this:



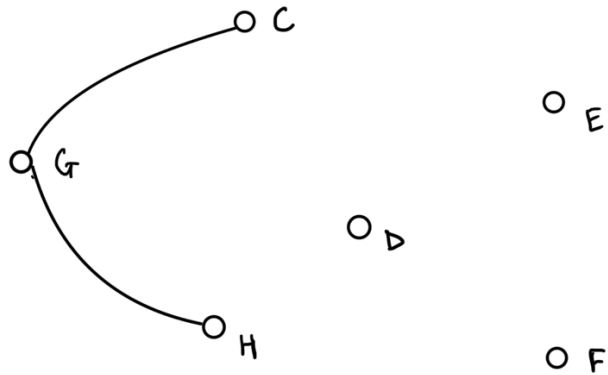

*Figure A 5 World with all lines bisectable?*

B: I agree that the straight line containing C, G and H can be A-bisected. However, the straight line containing just G and H cannot be A-Bisected.

A: Yeah. You are right. Maybe it isn't possible to make such a world.

T: That is an interesting thought. How would you prove that?

A: Well, any world with even one connection in it would have two points connected. The line connecting those two neighbors would be a straight line. That straight line cannot be A-bisected.

T: That is interesting. How about a world with no connections?

A: Then there are no straight lines.

T: Okay. But, since there are no straight lines, every straight line (of which there are none) which is there can be A-bisected. This is known as a vacuous truth. Look it up on the internet if it interests you. In any case, what you have shown is that in any world with at least one connection, there exists at least one straight line which cannot be bisected.

Let us move on to B-bisection. Can you create worlds in which every straight line can be B-bisected?

B: Well, the world with no connections!



T: Good point! Let me rephrase my question: Can you create a world with at least one connection where every straight line can be bisected?

B: Here is one:

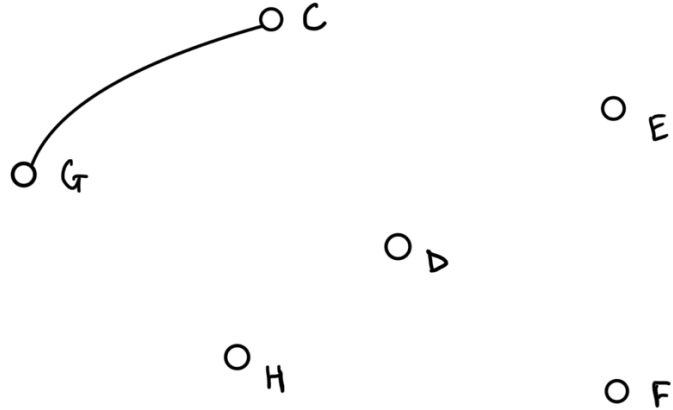

*Figure A 6 World with all lines bisectable?*

T: That looks like it works. Are there any others?

A: Wait! Does that actually work? What are we breaking the line containing C and G into?

B: The points G and C.

A: When we bisect a line in Euclidean Geometry, we get two lines. Is it okay to get two points?

T: That is a great observation. When we say bisection, we usually mean breaking a line into two lines. If we want this bisection into points to count as bisection, we have two options:

1.  We can allow bisection to result in points or lines

2.  We can say a point is a straight line

A: It seems weird to call a point a line.

T: Why? Does the point itself not represent the shortest path from the point to itself?



A: I guess it does.

B: But, if we consider a point to be a line, then in the world I have created, not every line can be B-bisected. For instance, if you pick any point, that cannot be B-bisected.

T: Great! That is very interesting. Let us use the term P-straight lines to mean straight lines which include points and N-straight lines to mean straight lines which do not include points. Let us also divide B-bisection into two types: P-B-Bisection which allows for bisecting lines into non-lines and N-B-Bisection which does not. In the world B has created, all N-straight lines can be P-B-Bisected. However, not all N-straight lines can be N-B-Bisected. P-straight lines can neither be P-B-Bisected nor N-B-Bisected.

Can you think of other worlds with six points where all N-straight lines can be P-B-Bisected?

B: Well, you can create as many pairs of connected points such that the pairs are not connected to each other.

T: Great! Are there any other possibilities?

A: I don't think there are since when we have at least three points connected to each other in any way, then there will always be a straight line which cannot be B-bisected – a straight line containing three points.

T: That is interesting. I want to come back to this, but implicit in your answer is another question. Rather than asking what worlds can be created where all straight lines can be bisected, we can ask: what sorts of lines can be A and B bisected in any world?

A: I think I've figured that out. A straight line can be A-bisected if it has an odd number of points while a straight line can be B-bisected if it has an even number of points.

T: On what basis can you conclude that.



A: Well, B-bisection requires us to break a line into two parts such that there are no common points in the two parts. The length of each part must be the same and hence the number of points must be the same. So, the number of points must be even. For A-bisection, the two parts have to share exactly one point. Hence, there are an odd number of points.

T: Great! Now getting back to the claim you made that if in any world there are three points which are connected, there will always be a N-straight line which cannot be P-B-bisected. Consider this world:

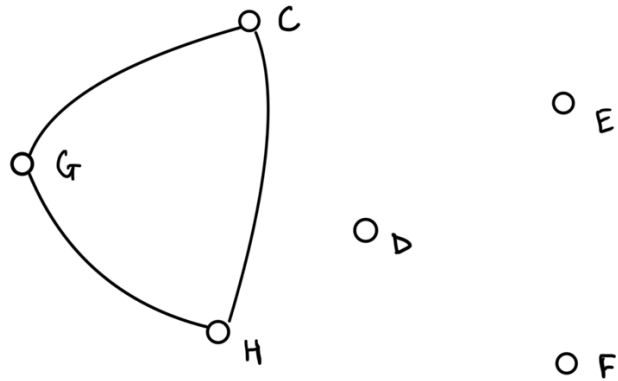

*Figure A 7 Another world with all lines bisectable?*

Can you find one of them here?

A: What about the straight line containing H, G and C, starting at H and ending at C?

T: Is that a straight line?

A: I think so.

T: What did we say a straight line was?

B: The shortest path between two points.

T: Is the line containing H, G and C, starting at H and ending at C the shortest path between H and C?



B: No it isn't! The direct path connecting H and C is.

T: So, in this world, are there N-straight lines which cannot be P-B-Bisected.

B: It doesn't seem like since even though there are paths which have an odd number of points, no straight line has an odd number of points.

T: Precisely. What we need to look for are not worlds where all paths have an even number of points, but all straight lines have an even number of points.

B: I can think of another – just create a similar triangular connection with the remaining three points.

T: Yes that works. The question to ask now is if there are worlds apart from the ones we have created so far where this is possible. I am going to leave this as an open question for you.

***Exploring some types of worlds.***

T: Let us now move on to exploring particular worlds with six points.

Consider this world.

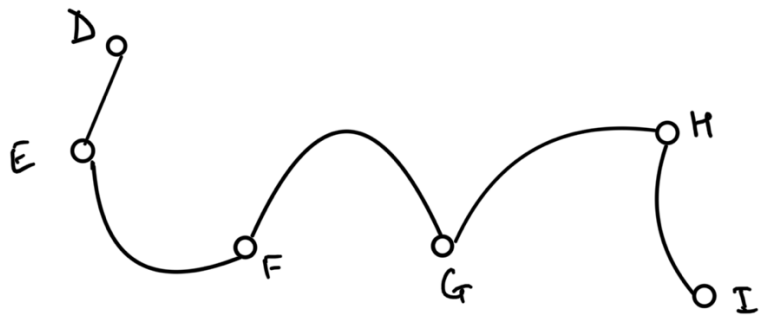

*Figure A 8 Simple world with six points*

What do straight lines in this world look like?

A: You can just pick any connected set of points. Those will be straight lines.

T: Great! What about circles?



A: Circles?! How can there be circles here?

T: What is a circle? Give me the definition you used in Euclidean Geometry

B: A circle is all the points which are equidistant from the center.

T: Yes. To create a circle you need a center and a radius. All points which are exactly radius-distance away from the center are part of the circle. Let us say that the length of the path connecting any two neighbors is 1. Is there a circle of radius 1 centered at G?

B: Yes. It contains the points F and H.

A: But is that a circle? Aren't circles supposed to be closed shapes?

T: Interesting. Clearly, if we allow for non-closed shapes to be circles, then the object B described is a circle. Let us call circles which include only closed shapes C-Circles. Let us call the more general notion NC-Circles. In order to understand what C-Circles are, we need to understand what 'closed shape' means in six point worlds. Let us return to that later and focus on NC-Circles which for now we will call circles for convenience.

The circle that B described contains two points. Are there circles which contain 1 point?

A: Yes. If we pick point D and a radius of anything less than 6, you get a circle with 1 point.

B: We can even get a circle with zero points – take any point and make the radius more than 6.

T: That is very interesting. The question is: do we want to call that a circle or not? We can side-step that question for now, but it is important to raise such questions. Can you create circles with three points.



A: No. If you pick any point as the center and pick any radius there are at most two points which are that radius away from the center. So, not just are three point circles impossible. Only circles with 2 or fewer points are possible.

T: Great! Let us now ask some more questions about circles. What do two-point and one point circles look like?

A: What do you mean by that?

T: I agree that the question I asked isn't completely clear. Maybe think about it as follows: what points can form one-point circles and what circles are those points a part of? What points can form two point circles and what circles are those a part of.

A: Any point can be a one-point circle. You can just take one of the end points and vary the radius. Then, you can take the other end point with radius 5 to get the last circle.

T: Great! Now, notice you can get a one-point circle containing the point E either as a radius 1 circle with center D or a radius 4 circle with center I. A question you could ask is: for each of the points, how many one-point circles are there containing that point? I'm going to leave that question for you to think about.

Now, does any pair of two points make a circle?

B: No. Take any two neighbors. There is no circle which contains both of them.

T: So, what types of pairs of points have a circle containing them and which do not?

B: Neighbors cannot. I think all other pairs can.

T: What about the points D and G?

B: Yeah. That cannot be made into a circle.

T: So, how would we characterize such points.



B: Okay. So, if the distance between two points is odd, then you cannot make a circle. Otherwise you can.

T: Great! Can you connect this to bisection?

B: Yes! If the straight line connecting two points can be A-bisected, you can make the two points into a circle with the bisection point as the center. Otherwise you cannot.

T: Interesting! We will come back to circles later in other worlds, but for now I want to move on to triangles. Are there triangles in this world?

A: I don't think so.

T: What are triangles?

A: They are closed shapes with three sides and three angles.

T: It is hard to see what we would mean by angle in this world, so let us think of triangles as closed shapes with three straight line sides. Do such triangles exist?

A: What would a closed shape be?

T: Good point. What does a closed shape mean in Euclidean Geometry?

A: It has some area inside of it and some area outside of it.

T: Let us try and think of an alternative way of thinking about it – I'm not sure what area means in this world. What can you say about the end points of the sides?

B: The end points of the sides have to connect to the end points of two other sides.

T: Let me suggest a slight change to that definition: each end points of a side has to connect to exactly one other side. So, now are there triangles in this world.

B: I don't think so.

T: Let me suggest a candidate: The shape with the sides EFG, EF and EG.

B: How is that a triangle?



T: Let us work through the definition. All of EFG, EF, and EG are straight lines. Each of their end points is shared by exactly one other side and hence the shape is closed.

B: But that shouldn't be a triangle.

T: Why not?

B: If we look at Euclidean Geometry and think about triangles there, we do not think of the following as a triangle:

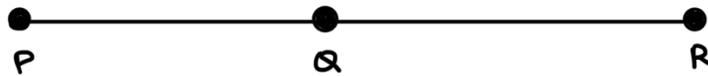

*Figure A 9 Collinear triangle*

It has three straight line sides: PQ, QR and PR and is closed in the way we have thought about closed shapes.

T: Okay, so if we want to exclude that possibility, how we do so?

A: We can just say that the three points should not be collinear – they should not all lie on a straight line.

T: Great! I'm going to do what we usually do. Let us have two types of triangles: C-Triangles and NC-Triangles. C-Triangles include both collinear and non-collinear triangles while NC-triangles include only non-collinear triangles. Are there NC-Triangles in this world?

A: No there are not.

T: Why?

A: If we take any three points in this world, there is always a straight line they are on. Hence, every three points are collinear.

T: Great! There are no NC-triangles, but there are examples of C-Triangles. Let us now move on to exploring another type of world – what I'm going to call a necklace worlds



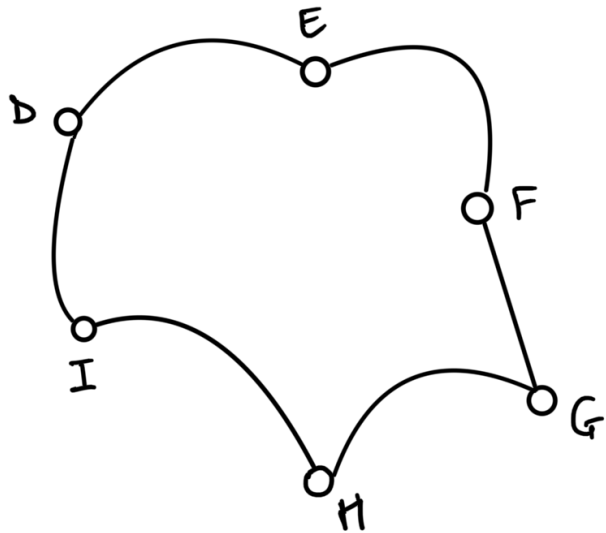

*Figure A 10 Necklace world with six points*

What do circles look like in this world?

A: Isn't the whole thing a circle?

T: Is it? What is the center?

A: Oh yeah. I guess it isn't – there is no point equidistant from all the other points.

T: I think what you did was to use the mental pictures of circles that you have from your past experience and use those to come to a conclusion here. Sometimes pictures can be misleading and we have to go by the definitions and rules of the world we are working in. Pick a point as the center and vary the radius to explore circles.

A: Let us pick D as the center. The circle of radius 1 has the points I and E. The circle of radius 2 has the points F and H. The circle of radius 3 has only the point G. The circle of radius 4 contains the points F and H.

T: Hold on a minute. I want you to spend some time thinking about the circle of radius 4. Does it contain the points F and H?

A: Yes. The points F and H are 4 away from D.



T: Are they 4 away from D? You said they are part of a circle of radius 2.

A: But there is a path of length 3 connecting D to F and H.

T: But, when we are talking about distance, what are we talking about. The length of any path or the length of the shortest path? Take the following example:

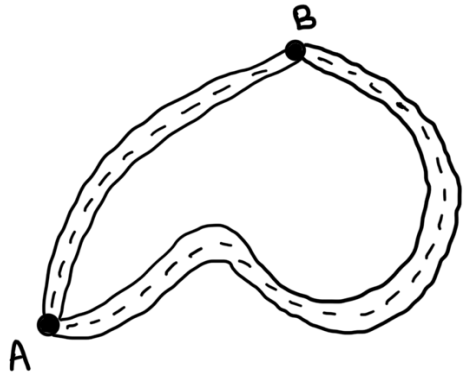

*Figure A 11 Roads between cities*

Say the above represents two paths between two cities. When we say the distance between two cities is 100km, which of the paths are we talking about.

A: The shortest path. So, there is not circle of radius 4 since there are no points at distance 4 away from the center.

T: Great! Let us try and characterize the circles by the number of points. What are the possibilities for the number of points?

A: 1 and 2, the same as in the other world.

T: What sort of single points can be 1-point circles?

B: Any of them. You just take the point opposite and radius 3.

T: What pairs of two point can be circles?

B: I think it is the same – the distance between them has to be even.

T: Is there a unique way to create such circles?



B: No. If you can create a circle using a point as the center, then you can create a circle with the same points using the opposite point as the center.

T: Great! Let us change the world slightly and continue exploring circles. Let us say rather than six points, we have five points.

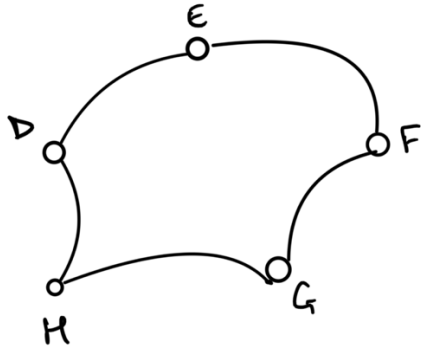

*Figure A 12 Necklace world with five points*

B: There are now no circles with single points.

A: Yes there are.

B: Which ones?

A: If we take a circle with radius 0, then each point is a circle with itself as the center.

T: Do we want to consider that to be a circle. Notice that in Euclidean geometry if we were to allow radius 0 circles, then every point there will be a circle. We usually do not allow that.

A: Why don't we take two types of circles: Z-Circles which allow 0 radius circles and NZ-Circles which do not.

T: Good! I want you to explore the consequences of Z-circles on your own. For now, let us call NZ-circles circles for convenience and continue with that concept.

A: In that case, there are no 1-point circles.



T: What about 2-point circles?

A: Once again, would it not be the case that those points even distance apart form two point circles and those non-even distance apart do not.

B: I don't think that is true here. If you take two neighbors, their distance is odd, but you can create a circle which includes both of them. For example, if you take the points G and F, you can create a circle of radius 2 with center G containing both of them. In fact, I think you can create a circle for every pair of points.

T: Can you justify that?

B: If you pick two points, there are two paths connecting them. At least one of those paths is of even length. Pick that path and think of it as a straight line. The point of bisection of that path is the center and the radius is half the length of that path.

T: Good job! How do you think this will generalize to 7-point and 8-point necklace worlds?

A: I think it has to do with even and odd number of points in the world.

T: In what way.

A: In even point worlds, there are 1-point circles. In fact every point can be a one-point circle by taking the point opposite and radius half of the number of points. In odd point worlds, 1-point circles do not exist since there is no point opposite to a particular point.

T: What about 2-point circles?

B: In odd-point worlds, any two points can be a circle. You just take the path which is even length connecting the two points and take the point of bisection of that path as the center of the circle. In even-point worlds, only points with even distance between them can be circles with two possible centers – the bisection points of either of the paths connecting them.



T: Great! Now, let us explore triangles in necklace worlds. Are there triangles in six point worlds?

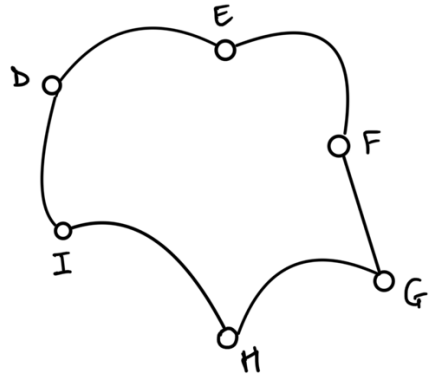

*Figure A 13 Necklace world with six points*

B: Yes. For example, DFH is a triangle – there are three straight line sides each end point of which is shared by another side. Also, these points are not collinear. So, we have NC-Triangles.

T: Of course, we still do have C-triangles which are not NC-triangles. For the purposes of moving forward, let us stick to NC-triangles and call them triangles. For a triangle we need three points. Let us call those vertices. Can any three points form the vertices of a triangle?

A: I think so.

T: Take the points I, D and E. Do those form a triangle?

A: Yes. The triangle with the sides ID, DE and EFGHI

T: Is that a triangle? What was required for a shape to be a triangle?

A: It had to be closed and the sides have to be straight lines.

T: Is that the case here?

A: It is clearly closed.

T: Yes. Also, ID and DE are straight lines. What about EFGHI?



A: Oh! It is not a straight line. The straight line connecting E and I is EDI.

T: So, it isn't a triangle. Not all vertices can form a triangle, assuming we stick to NC-triangles. Great! Let me ask you another question: The first triangle you mentioned was DHF. It is equilateral, right?

B: Yes. All the sides are of length 2. If we mean equal sided by equilateral, then DHF is equilateral.

T: Do 5-point worlds have equilateral triangles?

B: No they cannot.

T: Why?

B: Since the vertices of triangles cannot be collinear, every triangle contains all the points in the world. Hence, the length of the entire necklace must be divisible by 3 to get an equilateral triangle. Since 5 is not divisible by 3, you cannot have an equilateral triangle.

T: I would like you to explore isosceles triangles on your own time. You may enjoy doing that. For now, let us move on.

**Varying rules and exploring.**

T: There are many ways we can explore more worlds. We have explored changing the connections between points. We have also explored some examples of changing the number of points. Another possibility is to have directed paths. Think about a teleportation system where you can teleport from Delhi to Mars but not from Mars to Delhi. We could also change the length of paths between two neighbors. There could be different lengths for different connections. Maybe it takes more time to teleport between some locations than between others. Here is a world where we have some of these variations:



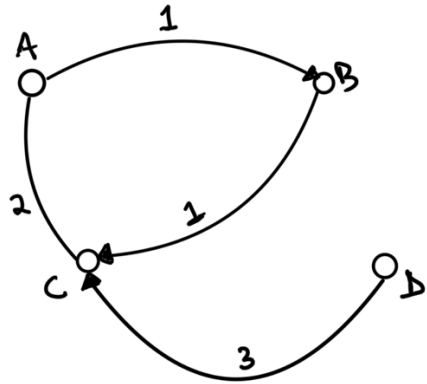

*Figure A 14 World with directed paths of different lengths*

The arrows represent directed paths while the paths with no arrows represent two-way paths. What is the straight line from A to C.

B: Just the path connecting them of length 2.

T: Is that the only straight line?

B: No. The path ABC has the same length. Hence, it is also a straight line.

T: What about the straight lines connecting C to A?

B: This time there is only one since the other path is directed.

T: What about A and D.

A: The only straight line path from D to A is DCA. There is no straight line path from A to D.

T: So, if you allow directed edges, the straight line from one point to another does not necessarily need to be the same as that from the second to the first point. That is interesting. There are many other variations you can explore. We are going to leave things here, but if you wish to explore further, think about classes of worlds rather that individual worlds. For example,



we explored necklace worlds where every point had exactly two neighbors. Maybe you can

explore worlds where every point has exactly three neighbors or four neighbors. You can explore

necklace worlds with directed paths in one direction. The possibilities are endless and you now

have the background to explore these things on your own.



## Appendix B: Pre and post-tests

The next few pages contain the pre and post tests administered to students on the first and last days of the course respectively.



## Pre-test

**Question 1:** The word 'toilider' is something I've made up and given a meaning. I've decided that the following are examples of toiliders:

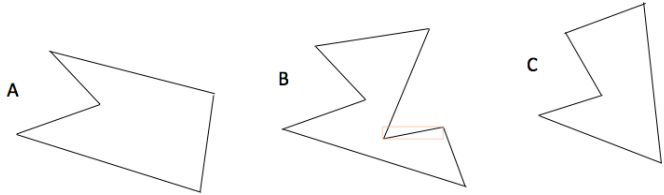

And the following are examples of non-toiliders:

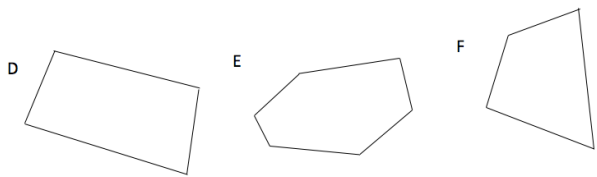

- Notice that D and F have four corners each, A and C have five corners each, E has six corners, and B has seven corners. Two corners are *adjacent* if a side connects them.

- We will use the term *diagonal* to define/describe a straight line connecting any non-adjacent corners.

The first column in the table below contains some definitions. Can each of these be a definition of a toilider? Circle your answer in the second column. If you say no, choose an example from the figures given above that shows that it cannot be a definition for toilider. Indicate the letter of the relevant figure in the second column.

| | (1) Possible Definition of Toilider? | (2) Example Number if Not |
|---|---|---|
| Any polygon | | |
| A polygon with 5 vertices | | |
| A polygon such that at least one diagonal lies outside of the polygon | | |
| A polygon with an odd number of vertices | | |
| A polygon with an even number of vertices | | |
| A polygon such that all the diagonals lie within the polygon | | |



**Question 2**: On a planet far, far away, we have heard that there are creatures called bleeks, plooks, glooks, pleeks, flooks, and blooks. A creature from that planet visited the Earth and it told us the following:

a. all blooks are bleeks

b. all bleeks are glooks

c. all plooks are glooks

d. all plooks are pleeks

e. all pleeks eat flooks

f. only pleeks eat flooks

g. no blook is a pleek

If we believe the creature, which of the following statements can we be sure are true, which are we sure are false, and which can we not say anything about? For example, from a) and b), we can conclude that the statement 'all blooks are glooks' is true. So we can place a tick mark under true in the first row below.

|  | True | False | Can't Say |
|---|---|---|---|
| all blooks are glooks |  |  |  |
| all plooks eat flooks |  |  |  |
| all plooks are bleeks |  |  |  |
| all bleeks eat flooks |  |  |  |
| only bleeks eat flooks |  |  |  |
| only plooks eat flooks |  |  |  |
| all pleeks are glooks |  |  |  |
| there is a blook who eats flooks |  |  |  |



## Post Test

**Question 1:** The word 'hobrub' is something I've made up and given a meaning. I've decided that the

following are examples of hobrub:

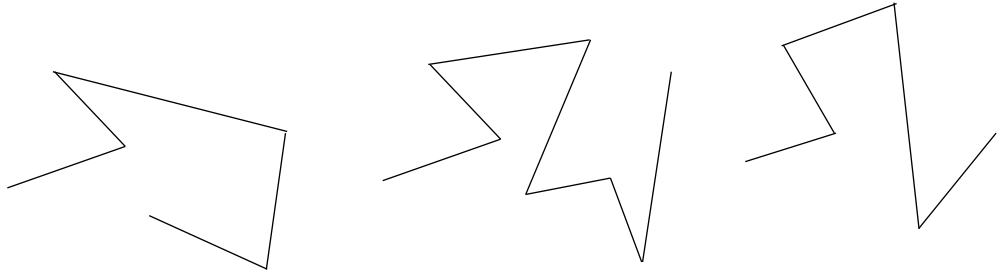

And the following are examples of non-hobrubs:

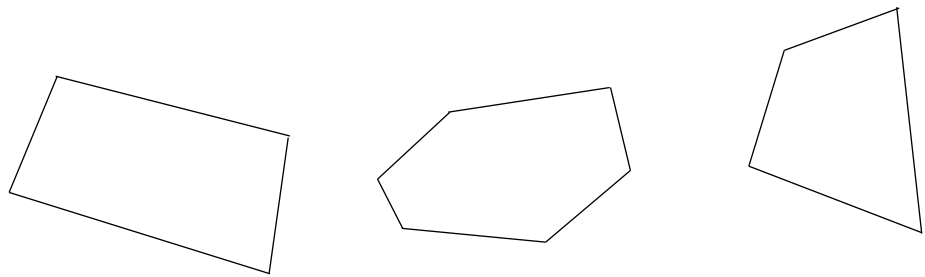

For each specification in the first column in the following table, indicate, in the second column, whether

can be the definition of a hobrub. Say Yes or No. If you say no, give an example that shows it cannot be a definition

for hobrub. Indicate the letter (A, B, C, etc.) of the relevant figure.

| | Possible Definition? | Example Number if Not |
|---|---|---|
| Any polygon | | |
| Any non-closed figure | | |
| A figure made only of straight line segment sides | | |
| A figure made of an odd number of straight line segments | | |
| A figure made of an even number of straight line segments | | |
| A figure made of straight line segments such that the total number of intersection points are even (meeting points count as intersection points) | | |



## Question 2

You are an investigator who has been asked by the local police to come to the island of Gadook. In this place, there are some criminal gangs: Hadooks, Padooks, Ladooks, Tadooks, Nadooks, and Badooks. It is possible for a person to be in many gangs at a time. Here is some information given to you about these gangs by the police:

a)   all Hadooks are Padooks
b)   all Padooks are Ladooks
c)   all Tadooks are Ladooks
d)   all Tadooks are Nadooks
e)   all Nadooks hate Badooks
f)   only Nadooks hate Badooks
g)   no Hadook is a Nadook

Assuming that the police has correct information, which of the following can you conclude?

For example, notice that from a) and b), we can conclude that 'all Hadooks are Ladooks' is true. This allows us to place a tick mark under true in the first row below.

| | True | False | Can't Say |
|---|---|---|---|
| all Hadooks are Ladooks | ✓ | | |
| all Tadooks hate Badooks | | | |
| all Tadooks are Padooks | | | |
| all Padooks hate Badooks | | | |
| only Padooks hate Badooks | | | |
| only Tadooks hate Badooks | | | |
| all Nadooks are Ladooks | | | |
| there is a Hadook who hates Badooks | | | |



**Appendix C: Initial and final survey**

The following pages contain the initial and final surveys. The initial survey focused on demographic information while the final survey focused on what students got out of the course.



**Initial Survey**

Adapted from NAEP Mathematics Student Questionnaire

**General Information**

Name:

Age:

Class:

Gender:

School:

**Home and School**

About how many books are there in your home?

Few (0–10)

Enough to fill one shelf (11–25)

Enough to fill one bookcase (26–100)

Enough to fill several bookcases (more than 100)

Select which of the following you have access to at home:

Access to the Internet

Your own bedroom



A desktop or laptop computer (including Chromebooks) that you can use

A tablet (for example, Surface Pro, iPad, Kindle Fire) that you can use

A smartphone (for example, iPhone, Samsung Galaxy, HTC One) that you can use

How often do you use the Internet for homework at home?

Never

Once or twice a month

Once or twice a week

Almost every day

Every day

In this school year, how often have you felt any of the following ways about your school?

Select one answer choice on each row.

|  | Never or hardly ever | Less than half the time | About half the time | More than half the time | All or almost all the time |
|---|---|---|---|---|---|
| I felt awkward and out of place at school |  |  |  |  |  |
| I felt happy at school |  |  |  |  |  |
| I felt that I learned something that I can use in my daily life |  |  |  |  |  |

How far in school did your mother go?

Did not finish 12th Grade

Finished 12th Grade but did not go to College



Graduated from College

I don't know

How far in school did your father go?

Did not finish 12th Grade

Finished 12th Grade but did not go to College

Graduated from College

I don't know

Does your mother work?

Yes

No

I don't know

Does your father work?

Yes

No

I don't know

**Mathematics**

How often do you receive help or tutoring with maths outside of school or after school?

Never

About once or twice a year



About once or twice a month

About once or twice a week

Every day or almost every day

For school this year, how often have you been asked to write long answers (several sentences or paragraphs) to questions on tests or assignments that involved math?

Never

Once

Two or three times

Four or five times

More than five times

How much does each of the following statements describe you? Select one answer choice on each row.

|  | Not like me at all | A bit like me | Somewhat like me | Quite a bit like me | Exactly like me |
|---|---|---|---|---|---|
| I enjoy doing maths. |  |  |  |  |  |
| I look forward to my maths class. |  |  |  |  |  |
| I think maths will help me even when I am not in school. |  |  |  |  |  |
| I am interested in the things I learn in maths. |  |  |  |  |  |
| I think making an effort in |  |  |  |  |  |
| maths is worthwhile. | | | | | |
| I think it is important to do well in maths. | | | | | |

How much does each of the following statements describe you? Select one answer choice on each row.

| | Not like me at all | A bit like me | Somewhat like me | Quite a bit like me | Exactly like me |
|---|---|---|---|---|---|
| I like complex problems more than easy problems | | | | | |
| I like activities that challenge my thinking abilities. | | | | | |
| I enjoy situations where I will have to think about something. | | | | | |
| I enjoy thinking about new solutions to problems | | | | | |

How much do you enjoy each of the following types of math activities? Select one answer choice on each row.

| | Not at all | A little bit | Somewhat | Quite a bit | A lot | Never done this |
|---|---|---|---|---|---|---|
| Addition, subtraction, multiplication, and division | | | | | | |



| | | | | | |
|---|---|---|---|---|---|
| Finding areas of shapes and figures | | | | | |
| Solving for probabilities and events (for example, card, coin, marble, and spinner problems) | | | | | |
| Solving equations or simplifying expressions | | | | | |
| Constructing and building different types of graphs (for example, bar graph, line graph, or box and whisker plots) | | | | | |
| Working with geometric shapes | | | | | |

To you, what is Mathematics?

**Experience with Knowledge Creation and Critical Thinking**

Were you a part of the ThinQ course in the 6th?

Yes

No

What workshops (if any) have you attended with me (Madhav)?

Summer 2017

Summer 2018

Summer 2019

Other (please specify)

What workshops have you had from other ThinQ members (including Mohanan, Aarthi, etc.)



**Final Survey**

1. What was your favorite part of the course? Why?

2. What was your least favorite part of the course? Why?

3. What did you find valuable in this course? Why?

4. What are one to three specific things about the course that could be improved?